\documentclass{article}

    \usepackage[final, nonatbib]{neurips_2022}

\usepackage[utf8]{inputenc} 
\usepackage[T1]{fontenc}    
\usepackage{hyperref}       
\usepackage{url}           
\usepackage{booktabs}       
\usepackage{amsfonts}       
\usepackage{nicefrac}       
\usepackage{microtype}      
\usepackage{xcolor}         
\usepackage[ruled, lined]{algorithm2e}
\usepackage{aligned-overset}
\usepackage{bbm}
\usepackage{dsfont}
\usepackage{enumitem}
\usepackage{wrapfig}
\usepackage{colortbl,booktabs}
\usepackage{xcolor}
\usepackage{caption}
\usepackage{subcaption}
\usepackage{wrapfig}

\newtheorem{lemma}{Lemma}

\newtheorem{remark}{Remark}

\newtheorem{theorem}{Theorem}

\def\Eset#1{\mathcal{E}(#1)}  
\newcommand\ntil{\tilde{n}}
\newcommand\vtil{\tilde{v}}
\newcommand\algl{\; \vspace{1mm} }

\def\modk#1#2{\mathrm{mod}(#1, #2  )}   
\def\floor#1{\left\lfloor #1 \right\rfloor} 
\def\ceil#1{\left\lceil #1 \right\rceil} 

\def\pr#1{\left( #1 \right ) }
\def\br#1{\left[ #1 \right ] }
\def\dr#1{\left\{#1\right\}} 

\def\eq#1{\begin{align*}#1\end{align*}}
\def\eql#1#2{\begin{equation}{#1}\begin{split}#2\end{split}\end{equation}}

\def\pleq{\preceq}

\newcommand{\norm}[1]{\left\|#1\right\|}
\newcommand{\abs}[1]{\left|#1\right|}
\newcommand{\E}{\mathbb{E}}
\newcommand{\R}{\mathbb{R}}
\newcommand{\Prb}{\mathbb{P}}
\newcommand{\Pin}{\boldsymbol{\Pi}}

\newcommand{\J}{\boldsymbol{\mathit{J}}}
\newcommand{\A}{\boldsymbol{\mathit{A}}}
\newcommand{\B}{\boldsymbol{\mathit{B}}}

\newcommand{\cG}{\mathcal{G}}

\newcommand{\W}{\boldsymbol{\mathit{W}}}

\newcommand{\I}{\boldsymbol{\mathit{I}}}
\newcommand{\x}{\boldsymbol{\mathit{x}}}
\newcommand{\y}{\boldsymbol{\mathit{y}}}
\newcommand{\bb}{\boldsymbol{\mathit{b}}}
\newcommand{\bs}{\boldsymbol{\mathit{s}}}
\newcommand{\bn}{\boldsymbol{\mathit{n}}}
\newcommand{\bh}{\boldsymbol{\mathit{h}}}
\newcommand{\bv}{\boldsymbol{\mathit{v}}}
\newcommand{\g}{\boldsymbol{\mathit{g}}}
\newcommand{\Om}{\mathcal{O}}

\renewcommand\AA{\boldsymbol{\mathit{A}}}
\newcommand\XX{\boldsymbol{\mathit{X}}}

\newcommand\Atil{\widetilde{\A}} 
\newcommand\II{\I} 
\newcommand\tp{^T}
\newcommand{\zero}{\mathbf{0}} 
\newcommand\WW{\W} 
\newcommand\Real{\R}

\newcommand\xx{\x} 
\newcommand\yy{\y}

\newcommand\Con{\boldsymbol{\Pi}}
\def\mx#1{\begin{pmatrix} #1 \end{pmatrix}}
\def\Xl#1{\boldsymbol{\mathit{X}}^{(#1)}}
\def\Yl#1{\boldsymbol{\mathit{Y}}^{
(#1)}}
\def\Gl#1{\boldsymbol{\mathit{G}}^{(#1)}}
\newcommand\gam{\gamma}
\def\gF#1{\nabla \boldsymbol{\mathit{F}}\big(\Xl{#1}
\big)}
\def\Wl#1{\WW^{(#1)}}
\def\nt#1{\left\| #1 \right\|}
\def\mt#1{\left\| #1 \right\|}

\newcommand\pa{c_1}
\def\nf#1{\left\| #1 \right\|_{\rm F}}
\def\Ly#1{\Phi^{(#1)}}

\def\Et#1{\E\br{\nf{#1}^2}}
\def\Ql#1{\boldsymbol{\mathit{Q}}^{(#1)}}
\newcommand\pb{c_2}

\def\Xa#1{\bar{\XX}^{(#1)}}
\def\xa#1{\bar{\xx}^{(#1)}}
\def\xl#1{\xx^{(#1)}}
\def\gl#1{\g^{(#1)}}
\def\yl#1{\yy^{(#1)}}
\def\gFa#1{\nabla \boldsymbol{\mathit{F}}\big(\bar{\XX}^{(#1)}\big)}
\def\jr#1{\left< #1 \right>} 

\def\Eb#1{\E\br{ #1  }} 
\def\gf#1{\nabla f\big(#1\big)} 
\def\Ot#1{\mathcal{O}\pr{#1}} 
\newcommand\bet{\beta}

\newcommand\BB{\boldsymbol{\mathit{B}}}

\newcommand\dffFinit{F_0}
\newcommand\ConXinit{C_0}
\newcommand\dFinit{D_0}

\def\nmod#1{#1 \ \mathrm{mod} \ n}

\def\normB#1{\Big\| #1 \Big\|}
\newcommand\Ctil{\widetilde{C}} 
\newcommand\Aol{\overline{\AA}^{(v_t)}} 
 
\newcommand\Col{\overline{C}} 
 
\newcommand\MM{\boldsymbol{\mathit{M}}}

\title{Communication-Efficient Topologies for Decentralized Learning with $\Om(1)$ Consensus Rate}

\author{
Zhuoqing Song$^{1}$\thanks{Equal Contribution. Corresponding Author: Kun Yuan},~ Weijian Li$^{2*}$, Kexin Jin$^{3*}$, Lei Shi$^{1,7}$, Ming Yan$^{4,5}$, Wotao Yin$^2$, Kun Yuan$^{2,6}$ \vspace{1mm}\\ 
    $^1$Fudan University\ \ $^2$Alibaba DAMO Academy~~ $^3$Princeton University\\ $^4$The Chinese University of Hong Kong, Shenzhen~~$^5$Michigan State University\\
    $^6$Peking University~~$^7$Shanghai Artificial Intelligence Laboratory\vspace{1mm}\\
{\small\texttt{zqsong19@fudan.edu.cn,\ weijian.li@alibaba-inc.com,\ kexinj@math.princeton.edu}}\\{\small\texttt{yanming@cuhk.edu.cn, \ leishi@fudan.edu.cn}}\\
{\small\texttt{wotao.yin@alibaba-inc.com,\quad kun.yuan@alibaba-inc.com}} 
}

\begin{document}

	\maketitle

	\begin{abstract}

 Decentralized optimization is an emerging paradigm in distributed learning in which agents achieve network-wide solutions by peer-to-peer communication without the central server. Since communication tends to be slower than computation, when each agent communicates with only a few neighboring agents per iteration, they can complete iterations faster than with more agents or a central server. However, the total number of iterations to reach a network-wide solution is affected by the speed at which the agents' information is ``mixed'' by communication. We found that popular communication topologies either have large maximum degrees (such as stars and complete graphs) or are ineffective at mixing information (such as rings and grids). To address this problem, we propose a new family of topologies, EquiTopo, which has an (almost) constant degree and a network-size-independent consensus rate that is used to measure the mixing efficiency.

 In the proposed family, EquiStatic has a degree of $\Theta(\ln(n))$, where $n$ is the network size, and a series of time-dependent one-peer topologies, EquiDyn, has a constant degree of 1. We generate EquiDyn through a certain random sampling procedure. Both of them achieve an $n$-independent consensus rate. We apply them to decentralized SGD and decentralized gradient tracking and obtain faster communication and better convergence,  theoretically and empirically. Our code is implemented through BlueFog and available at \url{https://github.com/kexinjinnn/EquiTopo}.

	\end{abstract}
	
	\section{Introduction}
	
	Modern optimization and machine learning typically involve tremendous data samples and model parameters. The scale of these problems calls for efficient distributed algorithms across multiple computing nodes. Traditional distributed approaches usually follow a centralized setup, where each node needs to communicate with a (virtually) central server. This communication pattern
	incurs significant communication overheads and long latency.
	
	Decentralized learning is an emerging paradigm to save communications in large-scale optimization and learning. In decentralized learning, all computing nodes are connected with some network topology (e.g., ring, grid, hypercube, etc.) in which each node averages/communicates locally with its immediate neighbors. This decentralized setup allows each node to communicate with fewer neighbors and hence has a much lower overhead in per-iteration communication. However, local averaging is less effective in ``mixing'' information, making decentralized algorithms converge slower than their centralized counterparts. Therefore, seeking a balance between communication efficiency and convergence rate in decentralized learning is critical.  
	
	The network topology (or graph) determines decentralized algorithms' per-iteration communication and convergence rate. The maximum graph degree controls the communication cost, whereas the connectivity influences the convergence rate. Intuitively speaking, a densely-connected topology enables decentralized methods to converge faster but results in less efficient communication since each node needs to average with more neighbors. Selecting an appropriate network topology is key to achieving light communication and fast convergence in decentralized learning.

	\begin{table}[t]
	    \centering 
		\caption{\small Comparison between different commonly-used topologies. ``Static Exp.'': static exponential graph; ``O.-P. Exp.'': one-peer exponential graph; ``E.-R. Rand'': Erdos-Renyi random graph $G(n,p)$ with probability $p=(1+a)\ln(n)/n$ for some $a>0$; ``Geo. Rand'': geometric random graph $G(n,r)$ with radius $r^2 = (1+a)\ln(n)/n$ for some $a > 0$. Undirected graphs can admit symmetric gossip matrices. If some graph has a  dynamic pattern, its associated gossip matrix will vary at each iteration.}
		\begin{tabular}{rccllc}
			\toprule
			\textbf{Topology} & \textbf{Connection} & \textbf{Pattern} & \textbf{Degree} & \hspace{-2mm}\textbf{Consensus Rate} & \textbf{size $n$}\\ \midrule
			Ring \cite{nedic2018network}              & undirect.          &                static                         &            $\Theta(1)$                &\hspace{-2mm}   $1 \hspace{-0.5mm}-\hspace{-0.5mm} \Theta(1/n^{2})$                 &  arbitrary  \\ 
			Grid \cite{nedic2018network}              & undirect.          &                static                         &            $\Theta(1)$                &\hspace{-2mm}   $1 \hspace{-0.5mm}-\hspace{-0.5mm} \Theta(\hspace{-0.5mm}1\hspace{-0.5mm}/\hspace{-0.5mm}(n\ln(n))\hspace{-0.5mm})$                 &   arbitrary \\ 
			Torus \cite{nedic2018network}              & undirect.          &          static                               &         $\Theta(1)$                   &\hspace{-2mm}     $1 \hspace{-0.5mm}-\hspace{-0.5mm} \Theta(1/n)$                &   arbitrary   \\ 
			Hypercube \cite{trevisan2017lecture}              &undirect.           &            static                             &       $\Theta(\ln(n))$                   &\hspace{-2mm}       $1 \hspace{-0.5mm}-\hspace{-0.5mm} \Theta(1/\ln(n))$         &      power of $2$     \\ 
			Static Exp.\cite{ying2021exponential} & directed &  static      & $\Theta(\ln(n))$         &\hspace{-2mm} $1 \hspace{-0.5mm}-\hspace{-0.5mm} \Theta(1/\ln(n))$   &  arbitrary\\
			O.-P. Exp.\cite{ying2021exponential} & directed &  dynamic      & $1 $         &\hspace{-2mm} finite-time conv.$^\dagger$ & power of $2$ \\
			E.-R. Rand \cite{nedic2018network} & undirect. & static &$\Theta(\ln(n))^{\diamond}$ &\hspace{-2mm} {$\Theta(1) $}  &arbitrary\\
			Geo. Rand \cite{boyd2005mixing} & undirect. & static&$\Theta(\ln(n))$ &\hspace{-2mm} $1 \hspace{-0.5mm}-\hspace{-0.5mm} \Theta(\ln(n)/n)$  & arbitrary\\
			{\color{blue}D-EquiStatic} & {\color{blue}directed} & {\color{blue}static}& {\color{blue}$\Theta(\ln(n))$}&\hspace{-2mm}    {\color{blue}$\rho \in (0,1)^\ddagger$}      & {\color{blue}arbitrary}\\
			{\color{blue}U-EquiStatic} & {\color{blue}undirect.} & {\color{blue}static} & {\color{blue}$\Theta(\ln(n))$}         &\hspace{-2mm} {\color{blue}$\rho \in (0,1)^\ddagger$}  & {\color{blue}arbitrary}\\
			{\color{blue}OD-EquiDyn} & {\color{blue}directed}& {\color{blue}dynamic}& {\color{blue}$1 $}         &\hspace{-2mm}  {\color{blue}$\sqrt{(1+\rho)/2}$} & {\color{blue}arbitrary}\\
			{\color{blue}OU-EquiDyn} & {\color{blue}undirect.} & {\color{blue}dynamic}& {\color{blue}$1 $}         &\hspace{-2mm} {\color{blue}$\sqrt{(2+\rho)/3}$}  & {\color{blue}arbitrary}\\
			\bottomrule
			\multicolumn{6}{l}{$^\dagger$\,\footnotesize{One-peer exponential graph has finite-time exact convergence only when $n$ is the power of $2$.}} \\
              \multicolumn{6}{l}{$^\diamond$ \footnotesize{$\Theta(\ln(n))$ is the averaged degree; its maximum degree can be $O(n)$ with a non-zero probability.}} \\
			\multicolumn{6}{l}{$^\ddagger$ \footnotesize{Constant $\rho=\Theta\pr{1}$  is independent of network-size $n$.}} \\
		\end{tabular}
		\label{Table:Summary}
	\end{table}
	
	\subsection{Prior arts in topology selections}
	\label{sec:topology-selection}
	
	\textbf{Gossip matrix and consensus rate.} Given a connected network of size $n$ and its associated doubly-stochastic gossip matrix $\W \in \mathbb{R}^{n\times n}$ (see the definition in \S~\ref{sec:notation}), 
 its consensus rate $\beta $ determines how effective the gossip operation $ \W \x$ is to mix information (see more explanations in \S~\ref{sec:notation}). 
 It is a long-standing topic in decentralized learning to seek topologies with both a small maximum degree and a fast consensus rate (i.e., a small $\beta$ as close to $0$ as possible). 
	
	\textbf{Static graphs.} Static topologies maintain the same graph connections throughout all iterations. The directedness, degree, and consensus rate of various common topologies are summarized in Table~\ref{Table:Summary}. The ring, grid, and torus graphs \cite{nedic2018network} are the simplest sparse topologies with $\Theta(1)$ maximum degree. However, their consensus rates quickly approach $1$ as network size $n$ increases, which leads to inefficient local averaging. The hypercube graph~\cite{trevisan2017lecture} maintains a nice balance between degree and consensus rate since $\ln(n)$ varies slowly with $n$. However, this graph cannot be formed when size $n$ is not the power of $2$.  The static exponential graph extends hypercubes to graphs with  any size~$n$, but its directed communications cannot enable symmetric gossip matrices required in well-known decentralized algorithms such as EXTRA \cite{shi2015extra}, Exact-Diffusion \cite{yuan2017exact1}, NIDS \cite{li2017decentralized}, D$^2$~\cite{tang2018d}. Two widely-used random topologies, i.e.,  the Erdos-Renyi graph \cite{nachmias2008critical,benjamini2014mixing} and the geometric random graph~\cite{beveridge2016best, boyd2005mixing}, are also listed in Table \ref{Table:Summary}. It is observed that the Erdos-Renyi graph achieves a network-size-independent consensus rate with a $\Theta(\ln(n))$ {\em averaged} degree in expectation. However, it is worth noting that the communication overhead in network topology is determined by the {\em maximum} degree. Since some nodes may have much more neighbors than others in a random realization, the maximum degree in the Erdos-Renyi graph can be $O(n)$ with a non-zero probability. Moreover, the random graphs listed in Table \ref{Table:Summary} are undirected. The may not be used in scenarios where directed graphs are preferred.

\textbf{Dynamic graphs.} Dynamic graphs allow time-varying topologies between iterations. When an exponential graph allows each node to cycle through all its neighbors and communicate only to a single node per iteration, we achieve the time-varying one-peer exponential graph \cite{assran2019stochastic,ying2021exponential}. When the network size $n$ is the power of $2$, a sequence of one-peer exponential graphs can together achieve periodic global averaging. However, its consensus rate is unknown for other values of $n$. A closely related work Matcha \cite{wang2019matcha} proposed a disjoint matching decomposition sampling strategy when training learning models. While it decomposes a static dense graph into a series of sparse graphs with small degrees, the consensus rates of these dynamic graphs are not established. 

Finally, it is worth noting that the consensus rates of all graphs (except for the Erdos-Renyi graph) discussed above are either unknown or dependent on size $n$. Their efficiency in mixing information gets less effective as $n$ goes large.

	\subsection{Main results}
	
	\textbf{Motivation.} Since existing network topologies suffer from several limitations, we ask the following questions. 
 {\em Can we develop topologies that have (almost) constant degrees and network-size-independent consensus rates that admit both symmetric and asymmetric matrices of any size? Can these topologies allow one-peer dynamic variants?}
 This paper provides affirmative answers.
	
	\textbf{Main results and contributions.} This paper develops several novel graphs built upon a set of basis graphs in which the label difference between any pair of connected nodes are {\em equivalent}. With a general name EquiTopo, these new graphs can achieve network-size-independent consensus rates while maintaining (almost) constant graph degrees. Our contributions are: 
	\begin{itemize}[leftmargin=1.5em]
		\item We construct a \underline{directed} graph named D-EquiStatic that has a network-size-independent consensus rate $\rho$ with a degree 
		{$\Theta(\ln(n))$}. 
  Furthermore, we develop a \underline{one-peer} time-varying variant named OD-EquiDyn to achieve a network-size-independent consensus rate with degree $1$. 
		
\item We construct a \underline{undirected} graph U-EquiStatic, which has a network-size-independent consensus rate $\rho$ with degree {$\Theta(\ln(n))$}. 
  It admits symmetric gossip matrices that are required by various important algorithms. We also develop a \underline{one-peer} time-varying and {undirected} variant named OU-EquiDyn to achieve a network-size-independent consensus rate with degree $1$. 

		\item We apply the EquiTopo graphs to two well-known decentralized algorithms, i.e., decentralized stochastic gradient descent (SGD) \cite{chen2012diffusion,lian2017can,koloskova2020unified} and stochastic gradient tracking (SGT) \cite{nedic2017achieving,xu2015augmented,di2016next,xin2020improved}, to achieve the state-of-the-art convergence rate while maintaining $\Theta(\ln(n))$ (with D/U-EquiStatic) or $1$ (with OD/OU-EquiDyn) degree in per-iteration communication. 
	\end{itemize}
	
	The comparison between EquiTopo and other common topologies in Table \ref{Table:Summary} shows that the EquiTopo family (especially the one-peer variants) has achieved the best balance between maximum graph degree and consensus rate. The comparison between EquiTopo and other common topologies when applying to DSGD and DSGT are listed in Tables \ref{Table:dsgd} and
	\ref{Table:dsgdstro}
	in Appendix \ref{Ap:TranIter} and Table \ref{Table:dsgta}.

	\textbf{Note.} This paper considers scenarios in which any two nodes can be connected when necessary. The high-performance data center cluster is one such scenario in which all GPUs are connected with high-bandwidth channels, and the network topology can be fully controlled. EquiTopo may not be applied to wireless network settings where two remote nodes cannot be connected directly. 
	
	\subsection{Other related works}
	In decentralized optimization, decentralized gradient descent \cite{nedic2009distributed,chen2012diffusion,sayed2014adaptive,yuan2016convergence} and dual averaging~\cite{duchi2011dual} are well-known approaches. While simple and widely used, their solutions are sensitive to heterogeneous data distributions. Advanced algorithms that can overcome this drawback include explicit bias-correction \cite{shi2015extra, yuan2017exact1, li2017decentralized},  gradient tracking \cite{nedic2017achieving,di2016next,qu2018harnessing,xu2015augmented}, and dual acceleration \cite{scaman2017optimal,uribe2020dual}. Decentralize SGD is extensively studied in \cite{chen2012diffusion,lian2017can,assran2019stochastic} to solve stochastic problems. It has been extended to directed \cite{assran2019stochastic,ying2021exponential}  or time-varying topologies \cite{koloskova2020unified,nedic2014distributed,wang2019matcha,ying2021exponential}, asynchronous settings \cite{lian2018asynchronous}, and data-heterogeneous scenarios \cite{tang2018d,xin2020improved,alghunaim2021unified,koloskova2021improved,lu2019gnsd} to achieve better performances.

	\section{Notations and Preliminaries}\label{sec:notation}

	\textbf{Notations.} 
	We let $\mathds{1}_n\in \R^n$ be the all-ones vector and $\I \in \mathbb{R}^{n\times n}$ be the identity matrix. Furthermore, we define  $\J = \frac 1 n\mathds{1}_n \mathds{1}_n^T$ and $\Pin = \I - \J$. 
	A matrix $\AA = [a_{ij}]\in \Real^{n\times n}$ is nonnegative if $a_{ij} \geq 0$ for all $1\leq i, j \leq n$. A nonnegative matrix $\AA$ is doubly stochastic if $\AA\mathds{1}_n = \AA^T\mathds{1}_n = \mathds{1}_n$.
	Given a matrix $\A \in \R^{m \times n}$, 	$\norm{\A}_2$ is its  spectral norm. 
	For $\x \in \R^n$, $\norm{\x}$ is its Euclidean norm. We let $[n] = \{1, \cdots, n\}$. 
	Throughout the paper, we define a $\mathrm{mod}$ operation that returns a value in $[n]$ as 
	\begin{align}
	i\ \mathrm{mod}\ n = 
	\left\{
	\begin{array}{ll}
	\ell & \mbox{if $i = kn + \ell$ for some $k \in \mathbb{Z}$ and $\ell \in [n-1]$}, \\
	n & \mbox{if $i = kn $ for some $k \in \mathbb{Z}$}.
	\end{array}
	\right.
	\end{align}
	\textbf{Network.}
	Given a graph $\mathcal{G}(\mathcal{V},\mathcal{E})$ with a set of $n$ nodes $\mathcal{V}$ and a set of directed edges $\mathcal{E}$. 
	An edge $(j,i) \in \mathcal{E}$ means node $j$ can directly send information to node $i$. For undirected graphs, $(j,i) \in \mathcal{E}$ if and only if $(i,j)\in \mathcal{E}$. Node $i$'s degree is the number of its in-neighbors $\abs{\{j|(j,i)\in \mathcal{E}\}}$. A one-peer graph means that the degree for each node is {\em at most} 1. 

\textbf{Weight matrices.}
 To facilitate the local averaging step in decentralized algorithms, each graph is associated with a nonnegative weight matrix $\W=[w_{ij}]\in\R^{n\times n}$, whose element $w_{ij}$ is non-zero only if $(j,i)\in\mathcal{E}$ or $i=j$. One benefit of an undirected graph is that it can be associated with a symmetric matrix. Given a nonnegative weight matrix $\W\in\R^{n\times n}$, we let $\mathcal{G}(\W)$ be its associated graph such that $(j,i)\in \mathcal{E}$ if $w_{ij}>0$ and $i\neq j$. 
 
 	\textbf{Consensus rate.} 
 	For weight matrices $\{\Wl{t}\}_{t\geq 0} \subseteq \Real^{n\times n}$, the consensus rate $\beta$ is the minimum  nonnegative number such that  for any $t \geq 0$ and vector $\xx \in \Real^n $ with the average $\bar{x} = \frac{1}{n}\sum_{i=1}^n x_i$, 
 	\begin{align*}
 	   \E\br{\nt{\Wl{t}\xx - \bar{x} \cdot\mathds{1}_n}^2 } \leq  \beta^{2}\nt{\xx - \bar{x} \cdot \mathds{1}_n}^2,  
 	   \end{align*}
 	   or equivalently, 
 	   $\E\br{\|\Con\Wl{t}\xx \|^2 } \leq  \beta^{2}\|\Con\xx \|^2.  $
 	   For $\Wl{t}\equiv \WW$,  $\beta$ essentially equals $\mt{\Con\WW }_2$.

	\section{Directed EquiTopo Graphs} \label{sec-D-EquiRand}

	\subsection{Basis weight matrices and basis graphs}\label{subsec-basis-graph}
	Given a graph of size $n$, we introduce a set of doubly stochastic {\em basis matrices} $\{\A^{(u, n)}\}_{u=1}^{n-1}$, where $\A^{(u, n)} = [a_{ij}^{(u, n)}] \in \R^{n \times n}$ with 
	\begin{align}\label{Au}
	a_{ij}^{(u, n)} = 
	\left\{
	\begin{array}{cl}
	\frac{n-1}{n}, & \mbox{if $i = (j + u )\ \bmod\ n $}, \\
	\frac{1}{n}, & \mbox{if $i = j$} ,\\
	0, & \mbox{otherwise}.
	\end{array}
	\right.
	\end{align}
	
\begin{figure}[!t]
		\centering
		\includegraphics[scale=0.22]{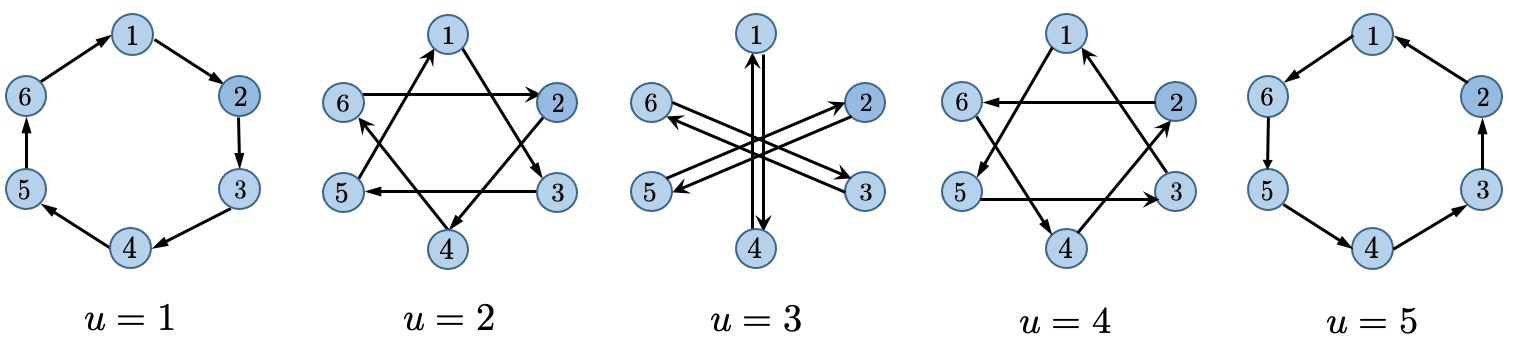}
		\caption{\small The set of the basis graphs $\{\cG(\A^{(u)})\}_{u=1}^5$ for $n=6$.}
		\label{fig:basis-lattice}
	\end{figure}

	Their associated graphs $\{\mathcal{G}(\A^{(u, n)} )\}_{u=1}^{n-1}$ are called {\it basis graphs}.  A basis graph $\mathcal{G}(\A^{(u, n)})$ has degree one and the same {\it label difference} $(i-j) \bmod n$ for all edges $(j,i)$.
	The set of five basis graphs $\{\mathcal{G}(\A^{(u, 6)})\}_{u=1}^5$ for $n=6$ is shown in Fig.~\ref{fig:basis-lattice}. 
	When $n$ is clear from the context, we omit it and write $\A^{(u)}$ instead. 

	\subsection{Directed static EquiTopo graphs (D-EquiStatic)}
	
	Our directed graphs are built on the above basis graphs, 
    and a weight matrix has the form
	\begin{equation}
	\label{eq:defW}
	\textstyle \W	= \frac{1}{M}\sum_{i=1}^M \A^{(u_i)},
	\end{equation}	
	where $u_i\in[n-1]$ and $M>0$ is the number of basis graphs we will sample. 
	Throughout this paper,  the multiset $\{u_i\}_{i=1}^M$ are called {\it basis index}.
	It is possible that $u_i = u_j$ for some $i\neq j$.
	Since each $\A^{(u)}$ has the form~\eqref{Au}, the matrix $\W$ is doubly stochastic, and all nodes of the directed graph $\cG(\W)$ have the same degree that is no more than $M $.

	Since $\cG(\W)$ is a directed static graph and built with $M$ basis graphs, we name it D-EquiStatic. The following theorem shows that we can construct a weight matrix $\W$ such that its consensus rate is independent of the network size $n$ by setting $M $ properly. The proofs of all theorems are in the Appendix.

	\begin{theorem}
		\label{lm-staconstr}
		Let $\A^{(u)}$ be defined by \eqref{Au} for any $u\in [n-1]$. 
		For any constant $\rho \in (0, 1)$,  
		we can 
		choose a sequence of $u_1, \cdots, u_M $  from $[n-1] $ with $M = \Theta\pr{\ln\pr{n}/\rho^2} $
		and construct the D-EquiStatic weight matrix $\W$ as in \eqref{eq:defW} such that  the consensus rate of $\WW$ is $\rho$, i.e.,  
		\begin{equation}
            \nt{\Con\WW\xx} \leq \rho\nt{\Con\xx},\ \forall \xx\in \Real^n  
		\end{equation}	
	\end{theorem}
	
	The graph $\cG(\W)$ has degree at most $M$. In the following, we will just say that the degree is $\Theta(\ln(n))$. As $\rho$ is tunable, we choose $\rho$ as a constant, e.g., $\rho = 0.5$. A method of constructing D-EquiStatic weight matrix $\W$ can be found in Appendix \ref{Ap:AlgD-EquiStati}.

	\begin{remark}
		Compared to all common topologies listed in Table \ref{Table:Summary}, D-EquiStatic achieves a better balance between degree and consensus rate. Moreover, D-EquiStatic works for any size $n \ge 2$. 
		Different from the Erdos-Renyi random graph and the geometric random graph, whose degree cannot be predefined before the implementation, we can easily specify the degree $M$ for D-EquiStatic. 
	\end{remark}

	\subsection{One-peer directed EquiTopo graphs (OD-EquiDyn)}

	While D-EquiStatic achieves a size-independent consensus rate with $\Theta(\ln(n))$ degree,  we develop a one-peer dynamic variant
 to further reduce its degree. Given a weight matrix $\W$ of form~\eqref{eq:defW} and its associated basis matrix $\{\A^{(u_i)}\}_{i=1}^M$, the one-peer directed variant, or OD-EquiDyn for short, samples a random $\A^{(u)}$ per iteration and utilizes it as the one-peer weight matrix, see Alg.~\ref{Alg:OD-EquiRand}. Since each node in $\cG(\A^{(u)})$ has exactly one neighbor, the graph $\cG(\W^{(t)})$ has degree one for every iteration $t$. Note that $\W^{(t)}$ is a random time-varying weight matrix. Its consensus rate (in expectation) can be characterized as below. 

\begin{algorithm}[H]{
\label{Alg:OD-EquiRand}
\caption{OD-EquiDyn weight matrix generation at iteration $t$}
\KwIn{constant $\eta \in (0,1)$;
	basis index 
	$\{u_1, u_2,\dots, u_M\}$ from a weight matrix $\W$ of form~\eqref{eq:defW};
	}

Pick $v_{t}$ from uniform distribution over the basis index $\{u_1, u_2,\dots, u_M\}$; 

Produce basis matrix $\A^{(v_t)}$ according to \eqref{Au};

\KwOut{$\W^{(t)}=(1-\eta)\I + \eta \A^{(v_{t})}$}
}
\end{algorithm}

\begin{theorem}
	\label{the:od-EquiRand}	
	Let the one-peer directed weight matrix $\W^{(t)}$ be generated by Alg.~\ref{Alg:OD-EquiRand}. It holds that 
	\begin{align*}\E\Big[\norm{ \Pin \W^{(t)} \x}^2\Big] \leq \big(1-2\eta(1-\eta)(1-\rho)\big)\norm{\Pin \x}^2,
	~~\forall \x \in \R^n \end{align*}
	where $\rho $ is the consensus rate of the weight matrix $\W$ (which can be tuned freely as in Theorem \ref{lm-staconstr}).
\end{theorem}

\begin{remark}
 The OD-EquiDyn graph has a degree of $1$ no matter how dense the input matrix $\W$ is. 
	When $\eta = 1/2$, which is used in our implementations, it holds that $\E \|{ \Pin \W^{(t)} \x}\|^2 \leq (1+\rho)/2 \norm{\Pin \x}^2$ for any $\x \in \R^n$. Thus, the OD-EquiDyn graph maintains the same $\Theta(1)$ degree as ring, grid, and the one-peer exponential graph but with a faster size-independent consensus rate.  
\end{remark}		

\begin{remark} \label{rmk-complete-od}
	For basis index $\{1,\cdots, n-1\}$ ($\W = \J$), Alg.~\ref{Alg:OD-EquiRand} returns a sequence of OD-EquiDyn graphs with $\E \|{ \Pin \W^{(t)} \x}\|^2 \leq \frac{1}{2} \norm{\Pin \x}^2$ when $\eta = \frac{1}{2}$ because $ \|\Pin \J \x\| = 0$ implies $\rho = 0$.   
\end{remark}
{
\begin{remark}
Although the Erdos-Renyi random graph also enjoys $\Theta(1)$ consensus rate and $\Theta(\ln(n))$ average degree (see Table \ref{Table:Summary}), its maximum degree could be as large as $\Theta(n)$ which implies that Erdos-Renyi random graphs could be highly unbalanced. 
Moreover, Erdos-Renyi random graphs are undirected graphs, while EquiStatic graphs can be both directed (Section~\ref{sec-D-EquiRand}) and undirected (Section~\ref{sec-U-EquiRand}).  In addition, the structure of EquiStatic allows simple construction of one-peer random graphs which preserve $\Theta(1)$ consensus, while it is still an open problem on whether Erdos-Renyi random graphs admit one-peer variants with $\Theta(1)$ consensus rate.   
\end{remark}
}
	
\section{Undirected EquiTopo Graphs}\label{sec-U-EquiRand}

The implementation of many important algorithms such as EXTRA \cite{shi2015extra}, Exact-Diffusion \cite{yuan2017exact1}, NIDS~\cite{li2017decentralized}, decentralized ADMM \cite{shi2014linear}, and the dual-based optimal algorithms \cite{scaman2017optimal,uribe2020dual,kovalev2021adom} rely on symmetric weight matrices. 
Moreover, devices in full-duplex communication systems can communicate with one another in both directions, and undirected networks are natural to be utilized.
These motivate us to study undirected graphs. 

\subsection{Undirected static EquiTopo graphs (U-EquiStatic)}

Given a D-EquiStatic weight matrix $\W$ and its associated basis matrices $\{\A^{(u_i)}\}_{i=1}^M$, we directly construct an undirected weight matrix name U-EquiStatic by  
\begin{equation}
\label{undirsta}
\textstyle \widetilde{\W} =\frac{1}{2}(\W + \W^T)= \frac{1}{2M}\sum_{i=1}^M (\A^{(u_i)} + [\A^{(u_i)}]^T),  
\end{equation}	
whose basis index are $\dr{u_i, - u_i }_{i=1}^M $ because 
$\A^{(- u)} = [\A^{(u)}]^T $. 

Since $\widetilde{\W}$ is built upon $\W$ and $\Con\W = \W\Con$,  the following theorem follows directly from 
\begin{align*}\mt{\Con\widetilde{\W }}_2 = \frac{1}{2}\mt{\Con\W + (\Con\W)\tp}_2 \leq \frac{1}{2}\pr{\mt{\Con\W}_2 + \mt{\pr{\Con\W}\tp}_2} = \mt{\Con\W}_2.\end{align*}

\begin{theorem}\label{thm-u-equirand}
	Let $\W$ be a D-EquiStatic matrix  with consensus rate $\rho $ and $\widetilde{\W}$ be the U-EquiStatic matrix defined  by \eqref{undirsta}. It holds that 
	\eql{\label{eq:tildeW1}}{
		\norm{ \Pin \widetilde{\W} \x}
		\leq \rho \norm{\Pin \x}, ~~\forall \x \in \R^n.
	}

\end{theorem}

\begin{figure}[!t]
	\centering
	\includegraphics[scale=0.22]{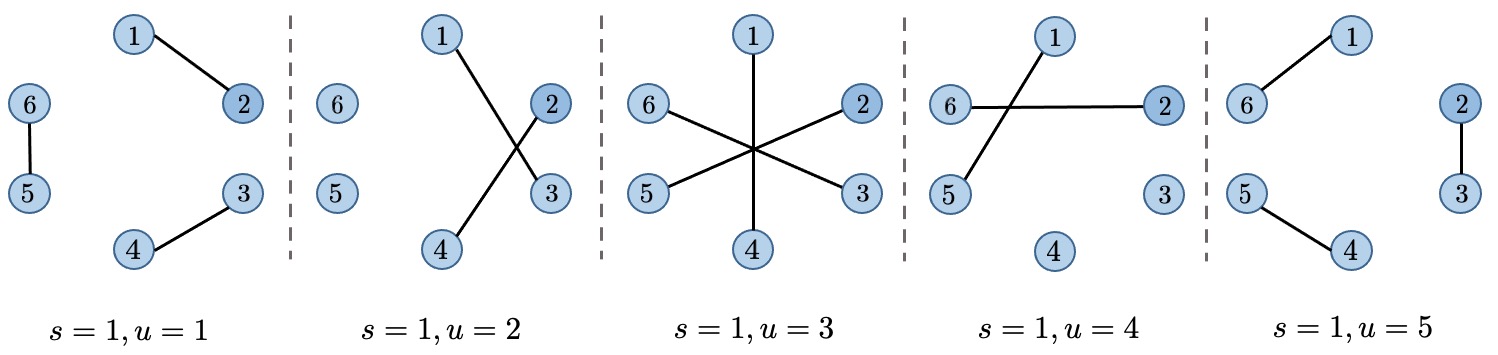}
	\caption{\small A few realizations of the OU-EquiDyn graphs for $n=6$, $s=1$, and $u\in\{1,2,3,4,5\}$.}
	\label{fig:ou-equirand}
\end{figure}

\subsection{One-peer undirected EquiTopo graphs (OU-EquiDyn)}
\label{OU-EquiDyn}

Constructing a one-peer undirected graph OU-EquiDyn is not as direct as U-EquiStatic because $\frac{1}{2}(\A^{(u)} + [\A^{(u)}]^T)$ admits a graph with degree~2, see Appendix \ref{Ap:basisouequidyn} for an illustration.

Alg.~\ref{Alg:CuMN} shows a method to construct a series of OU-EquiDyn  matrices with degree 1.
Starting from node $s$, we connect a node with the $u^{\mbox{th}}$ node after it, as long as both of them have not been connected to any other nodes. 
Fig.~\ref{fig:ou-equirand} illustrates the process when $s=1$ and $u = 1,\cdots, 5$ for a network of size $6$. Some nodes have no neighbors at realizations $u=2$ and $u=4$. This phenomenon is caused by the restriction that each node has no more than one neighbor.  For instance, when $u= 2$, node $5$ wants to connect with node $1$ but node $1$ has already been connected to node $3$. 
Thus, there exist node pairs that are never connected when $s=1$. 
To resolve this issue, we let the starting index $s$ be sampled randomly from $[n]$. Fig.~{\ref{fig:ou-oneundirsu}} in Appendix {\ref{Ap:basisouequidyn}} illustrates the scenarios when $s=3$. It is observed that the pairs $\{3, 5\}$ and $\{4,6\}$ are now connected to each other.
The node version of OU-EquiDyn is illustrated in Alg.~\ref{Alg:CuMNagentpersp1} of 
Appendix. \ref{Ap:EqAgetAlg}.

\begin{algorithm}[H]
\caption{OU-EquiDyn weight matrix generation at iteration $t$}
\label{Alg:CuMN}
\KwIn{$\eta \in (0,1)$;
	basis index  
	$\{u_i, - u_i\}_{i=1}^M$ from a symmetric weight matrix $\widetilde{\W}\in\R^{n\times n}$ of form~\eqref{undirsta};} 

Pick $v_t \in \{u_i, -u_i\}_{i=1}^M$ and  $s_t \in [n]$ uniformly at random\;
Initialize  
$ \A = 
[ a_{ij}]= \I$ and $b_i = 0$, $\forall i \in [n]$\; 
\For {$j = (s_t:s_t+n-1 \mod n)$ }{ 
	
	$i = (j + v_t) \mod n$\; 
	
	\If{$ b_i = 0$ and $b_j = 0 $}{ 
		$ a_{ij} =  a _{ji} = (n-1)/n$\;
		$ a_{ii}  =  a_{jj} = 1/n$\;		
		$b_i = 1 $, $b_j = 1$; 
}}
\KwOut{$\widetilde \W^{(t)} = (1-\eta)I +\eta  \A$}
\end{algorithm}

\begin{theorem}
	\label{the:ou_equirand}	
	Let $\widetilde{\W}$ be a U-EquiStatic matrix  with consensus rate $\rho $, and $\widetilde{\W}^{(t)}$ be an OU-EquiDyn matrix generated by Alg.~\ref{Alg:CuMN}, it holds that
	\begin{align*}\E\br{\big\| \Pin \widetilde{\W}^{(t)} \x\big\|^2} \leq \big(1-\frac{4}{3}\eta(1-\eta)(1-\rho)\big)\norm{\Pin \x}^2,
	~~\forall \x \in \R^n.\end{align*}
\end{theorem}

\begin{remark}
	Theorem \ref{the:ou_equirand} implies that the OU-EquiDyn graph can achieve a size-independent consensus rate with a degree at most $1$. When $\eta = 1/2$, it holds that $\E\|\Pin \W^{(t)} \x\|^2 \leq [(2+\rho)/3] \|\Pin \x\|^2$. 
\end{remark}

When $\widetilde{\W} = \J$ and the basis index $\{1,\cdots, n-1\}$ are input to Alg.~\ref{Alg:CuMN}, we obtain an OU-EquiDyn sequence $\widetilde{\W}^{(t)}$ such that $\E\|\Pin \widetilde{\W}^{(t)} \x\|^2 \leq (2/3) \|\Pin \x\|^2$.

\begin{remark}
An alternative OU-EquiDyn matrix construction that relies on the Euclidean algorithm  is in Appendix \ref{Ap:EOUA}.
\end{remark}

\section{Applying EquiTopo Matrices to Decentralized Learning}

We consider the following distributed problem over a network of $n$ computing nodes:
\begin{equation}
\label{form}
\min_{\x \in \R^d}~~
f(\x) = \frac 1n \sum_{i=1}^n
f_i(\x)
\end{equation}
where $f_i(\x) :=
\E_{\xi_i \sim \mathcal D_i} [F(\x; \xi_i) \big]$.
The function $f_i(\x)$ is kept at node $i$, and $\xi_i$ denotes the local data that follows the local distribution $\mathcal D_i$. Data heterogeneity exists if local distributions $\{\mathcal D_i\}_{i=1}^n$ are not identical. Throughout this section, we let $\x_i^{(t)}$ be  node $i$'s local model at iteration $t$, and $\bar{\x}^{(t)} = \frac 1n\sum_{i=1}^ n \x^{(t)}_i$.

\textbf{Assumptions.} We make the following standard assumptions to facilitate analysis.

\textbf{A.1} \textit{Each local cost function $f_i(x)$ is differentiable, and there exists a constant $L > 0$ such that
$\norm{\nabla f_i(\x)- \nabla f_i(\y)} \le L \norm{\x-\y}$ for all $\x, \y \in \R^d$.}

\textbf{A.2} 
\textit{Let $\g_i^{(t)}=\nabla F(\x_i^{(t)}; \xi_i^{(t)})$. There exists $\sigma^2 >0$ such that for any $t$ and $i$
\begin{equation*}
\E_{\xi_i^{(t)}\sim \mathcal D_i}{\g_i^{(t)}}=\nabla f_i (\x_i^{(t)}),~~
{\rm and}~~
\E_{\xi_i^{(t)}\sim \mathcal D_i}\Big[\norm{\g_i^{(t)} -\nabla f_i(\x_i^{(t)})}^2\Big]
\le \sigma^2.
\end{equation*}}

\textbf{A.3} 
\textit{ (For DSGD only) There exists $b^2$ such that
$\frac 1n \sum_{i=1}^n \norm{\nabla f_i(\x) - \nabla f(\x)} \le b^2$
for all $\x \in \R^d$.}

\subsection{Decentralized stochastic gradient descent}
The decentralized stochastic gradient descent (DSGD) \cite{chen2012diffusion,lian2017can,koloskova2020unified} is given by
\begin{equation}
\label{dsgd}
\textstyle \x_i^{(t+1)}=
\sum_{j=1}^n w_{ij}^{(t)}\big(\x_j^{(t)}-\gamma \g_j^{(t)}\big), 
\end{equation}
where the weight matrix $\W^{(t)} = [w^{(t)}_{ij}]$ can be time-varying and random. Applying the EquiTopo matrices discussed in \S~\ref{sec-D-EquiRand}-\ref{sec-U-EquiRand}, we achieve the following convergence results, whose proof follows Theorem $2$ in \cite{koloskova2020unified} directly, and is omitted here.
More results for DSGD are given in Appendix \ref{dsgd:stron}.

\begin{theorem}
	\label{the:dsgd}
	Consider the DSGD algorithm \eqref{dsgd}. Under Assumptions A.1-A.3, it holds that
	\begin{equation*}
	\begin{aligned}
	\frac 1{T+1} \sum_{t=0}^{T} \E \Big[ \norm{\nabla f(\bar{\x}^{(t)})}^2\Big]=
	\Om \Big( \frac {\sigma}{\sqrt{nT}}
	+\frac {\beta^{\frac 23} \sigma^{\frac 23}}{T^{\frac 23}(1-\beta)^{\frac 13}}
	+\frac {\beta^{\frac 23} b^{\frac 23}}{T^{\frac 23}(1-\beta)^{\frac 23}}
	+\frac {\beta}{T(1- \beta)}
	\Big),
	\end{aligned}
	\end{equation*}	
	\begin{itemize}[leftmargin = 1.5em]
		\item where $\beta = \rho$ with D-EquiStatic $\W$ or U-EquiStatic $\widetilde{\W}$; 
		\item where $\beta = \sqrt{(1+\rho)/2}$ for OD-EquiDyn $\W^{(t)}$ (Alg.~\ref{Alg:OD-EquiRand} with $\eta = 1/2$), and $\beta = \sqrt{(2+\rho)/3}$ for OU-EquiDyn $\widetilde{\W}^{(t)}$ (Alg.~\ref{Alg:CuMN} with $\eta = 1/2$). 
	\end{itemize}
\end{theorem}

For a sufficiently large $T$, the term $\Om(1/\sqrt{nT})$ dominates the rate, and we say the algorithm reaches the linear speedup stage. The transient iterations are referred to as those iterations before an algorithm reaches the linear-speedup stage. 
We compare the per-iteration communication, convergence rate, and transient iterations of DSGD over various topologies in Tables \ref{Table:dsgd} and \ref{Table:dsgdstro} of the Appendix. It is observed that OD/OU-EquiDyn endows DSGD with the lightest communication, fastest convergence rate, and smallest transient iteration complexity.

\subsection{Decentralized stochastic gradient tracking algorithm}

The decentralized  stochastic gradient tracking algorithm (DSGT)~\cite{nedic2017achieving,di2016next,qu2018harnessing,xu2015augmented,xin2020improved}  is given by
\begin{equation}
\begin{aligned}
\label{dsgta}
\x_i^{(t+1)} =&\textstyle \sum_{j=1}^n w_{ij}^{(t)} \big(\x_j^{(t)}  - \gam\y_j^{(t)}\big);\\
\y_i^{(t+1)} =&\textstyle \sum_{j=1}^n w_{ij}^{(t)} \y_j^{(t)}  + \g_i^{(t+1)}- \g_i^{(t)},~~\y_i^{(0)}=\g_i^{(0)}.\\
\end{aligned}
\end{equation}		
The following result of DSGT does not appear in the literature since it admits an improved convergence rate for stochastic decentralized optimization over asymmetric or time-varying weight matrices. Existing works on DSGT assume weight matrix to be either symmetric \cite{koloskova2021improved,alghunaim2021unified} or static \cite{xin2020improved}.
\begin{theorem}
\label{the:dsgt}
	Consider the DSGT algorithm in~\eqref{dsgta}. 
	If $\{\Wl{t}\}_{t\geq 0}$ have consensus rate $\beta$, then
	under Assumptions A.1-A.2, it holds for $T \ge \frac{1}{1 - \beta}$ that
	\begin{equation*}
	\begin{aligned}
	        &\frac{1}{T + 1  }\sum_{t=0 }^{T} \E\br{\|\nabla f(\bar{\xx}^{(t)})\|^2} 
        =
        \mathcal{O}\Big( \frac{ \sigma  }{\sqrt{n T}} + \frac{\sigma^{\frac 23}}{(1-\beta)T^{\frac 23}}
        + \frac{1  }{\pr{1 - \bet }^2 T}\Big).  
	\end{aligned}
	\end{equation*}	
	 
\end{theorem}	
When utilizing the EquiTopo matrices, the corresponding $\beta$ is specified in Theorem \ref{the:dsgd}. Note that DSGT achieves linear speedup for large $T$. The per-iteration communication and convergence rate comparison of DSGT over different topologies is in Table \ref{Table:dsgta}. OD/OU-EquiDyn endows DSGT with the lightest communication, fastest convergence rate, and smallest transient iteration complexity. 

\begin{table}[t]
	\caption{\small Per-iteration communication and computation complexity of the DSGT under different topologies.}
	\begin{tabular}{rccllc}
		\toprule
		& \textbf{Topology} & \textbf{Per-iter Comm.} &  \textbf{Convergence Rate}  & \textbf{Trans. Iters.} \\ \midrule
		&Ring                   &$\Theta(1)$                        &$\Ot{\frac{ \sigma  }{\sqrt{n T}} +
		\frac{n^2 \sigma^{\frac{2}{3}}}
		{T^{\frac{2}{3}}} + \frac{n^4 }{T}}$     
		&$\Om(n^{15})$              \\ 
		&Torus                  &$\Theta(1)$                      & $\Ot{\frac{ \sigma  }{\sqrt{n T}} + \frac{n \sigma^{\frac{2}{3}}}
		{T^{\frac{2}{3}}} + \frac{n^2 }{T}}$       
		&$\Om(n^9)$                 \\
		&Static Exp.    &$\Theta(\ln(n))$      &$\Ot{\frac{\sigma}{\sqrt{n T}} +
		\frac{\ln(n)\sigma^{\frac{2}{3}}}
		{T^{\frac{2}{3}}}
	  + \frac{\ln^2(n)}{T}}$     
		&$\Om(n^3 \ln^6(n))$        			\\
		&O.-P. Exp.		 & $1$      & $\Ot{\frac{ \sigma  }{\sqrt{n T}} +  	\frac{\ln(n)\sigma^{\frac{2}{3}}}
		{T^{\frac{2}{3}}} + \frac{\ln^2(n)}{T}}$    
		&$\Om(n^3 \ln^6(n))$ 		   \\
		&{\color{blue}D(U)-EquiStatic} & {\color{blue}$\Theta(\ln(n))$}  & {\color{blue}$\Ot{\frac{ \sigma  }{\sqrt{n T}} +\pr{\frac{ \sigma}{T}}^{\frac{2}{3}} + \frac{1  }{T}}$} 
		&{\color{blue}$\Om(n^3)$}				\\
		&{\color{blue}OD (OU)-EquiDyn} & {\color{blue}$1$ }
		&{\color{blue}$\Ot{\frac{ \sigma  }{\sqrt{n T}} +\pr{\frac{ \sigma}{T}}^{\frac{2}{3}} + \frac{1  }{T}}$} 
		& {\color{blue}$\Om(n^3)$}  \\
		\bottomrule
	\end{tabular}
	\label{Table:dsgta}
\end{table}

\begin{figure}[!ht]
\centering
\begin{minipage}{.45\textwidth}
  \centering
  \includegraphics[width=0.9\textwidth]{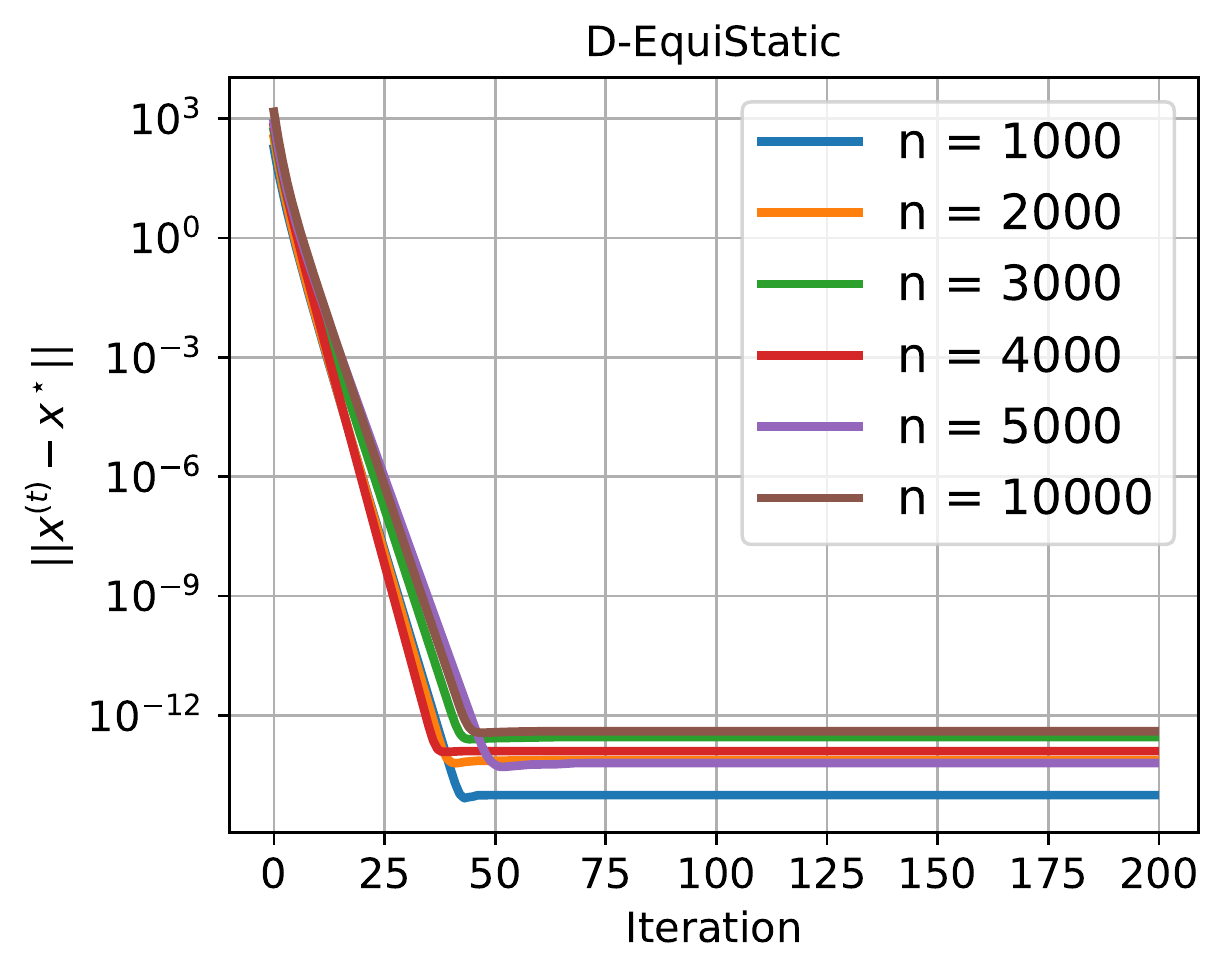}
  \caption{\small The D-EquiStatic topology can achieve network-size independent consensus rate.}
	\label{fig:n-independent-partial}
\end{minipage}%
\hspace{1cm}
\begin{minipage}{.45\textwidth}
  \centering
  \includegraphics[width=0.9\textwidth]{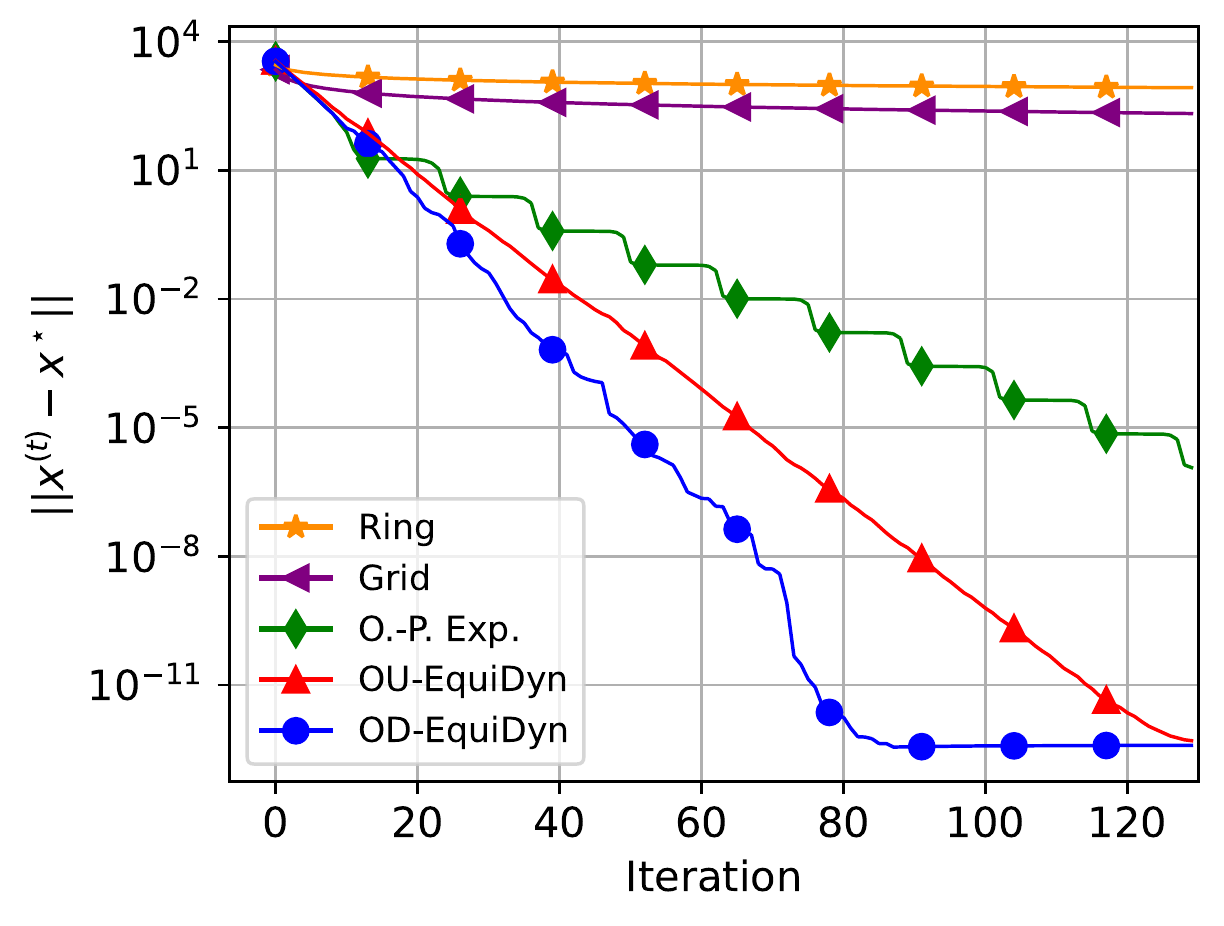}
	\caption{\small OD/OU-EquiDyn is faster than other topologies (i.e., ring, grid, and one-peer exponential graph) with $\Theta(1)$ degree in consensus rate.}
	\label{fig:average-consensus}
\end{minipage}
\end{figure}

\section{Numerical Experiments}
This section presents experimental results to validate EquiTopo's network-size-independent consensus rate and its comparison with other commonly-used topologies in DSGD on both strongly-convex problems and non-convex deep learning tasks. More experiments for EquiTopo in DSGT and all implementation details are referred to Appendix \ref{Ap:Exp}.

\textbf{Network-size independent consensus rate.} 
This simulation examines the consensus rates of all four EquiTopo graphs. We recursively run the gossip averaging $\x^{(t+1)} = \W^{(t)} \x^{(t)}$ with $\x^{(0)} \in \R^{n}$ initialized arbitrarily and $\W^{(t)} \in \R^{n \times n}$ generated as D-EquiStatic (Eq.~\eqref{eq:defW}), OD-EquiDyn (Alg.~\ref{Alg:OD-EquiRand}), U-EquiStatic (Eq.~\eqref{undirsta}), and OU-EquiDyn (Alg.~\ref{Alg:CuMN}), respectively. Fig.~\ref{fig:n-independent-partial} depicts how the quantity $\|\x^{(t)} - \J\x^{(0)}\|$ evolves when $n$ ranges from $1000$ to $10,000$ with D-EquiStatic topology. See Appendix \ref{Ap:Exp} for all other EquiTopo graphs, which also achieve network-size independent consensus rates when $n$ varies. These results are consistent with Theorems~\ref{lm-staconstr} - \ref{the:ou_equirand}.

\textbf{Comparison with other topologies.} We now compare EquiTopo's consensus rate with other commonly-used topologies. Fig.~\ref{fig:average-consensus} illustrates the performance of several graphs with $\Theta(1)$ degree when running gossip averaging. We set $n=4900$ so that the grid graph can be organized as $70 \times 70$. OD/OU-EquiDyn is much faster than other topologies. Note that each node in OD-EquiDyn, OU-EquiDyn, and O.-P. Exp. has exactly one neighbor per iteration. More experiments on graphs with $\Theta(\ln(n))$ degrees and  on scenarios with smaller network sizes are in Appendix \ref{Ap:Exp}.

\begin{figure}[!ht]
	\begin{center}
	\includegraphics[width=0.45\textwidth]{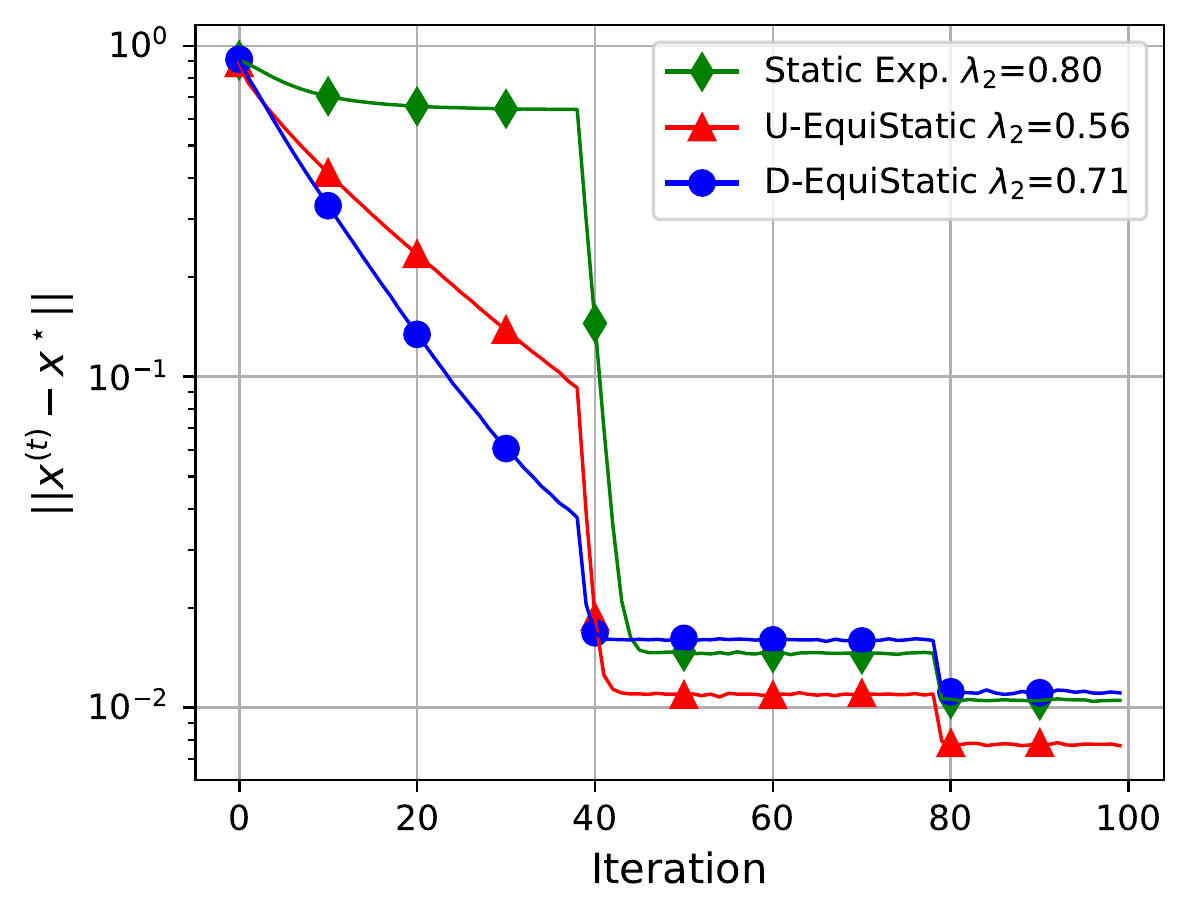}
	\end{center}
	\caption{\small D/U-EquiStatic in DSGD. $\lambda_2$ is the second largest eigenvalue.}
	\label{fig:dsgd}
\end{figure}

\textbf{DSGD with EquiTopo: least-square.} We next apply D/U-EquiStatic graphs to DSGD when solving the distributed least square problems. In the experiment, we let $n=300$ and set $M = 9$ so that D/U-EquiStatic has the same degree as the exponential graph.  Fig.~\ref{fig:dsgd} depicts that D/U-EquiStatic converges much faster than a static exponential graph, especially in the initial stages when the learning rate is large. The U-EquiStatic performs slightly better than D-EquiStatic since its bi-directional communication enables the graph with better connectivity.

\begin{figure}[h!]
	\centering
	\includegraphics[scale=0.5]{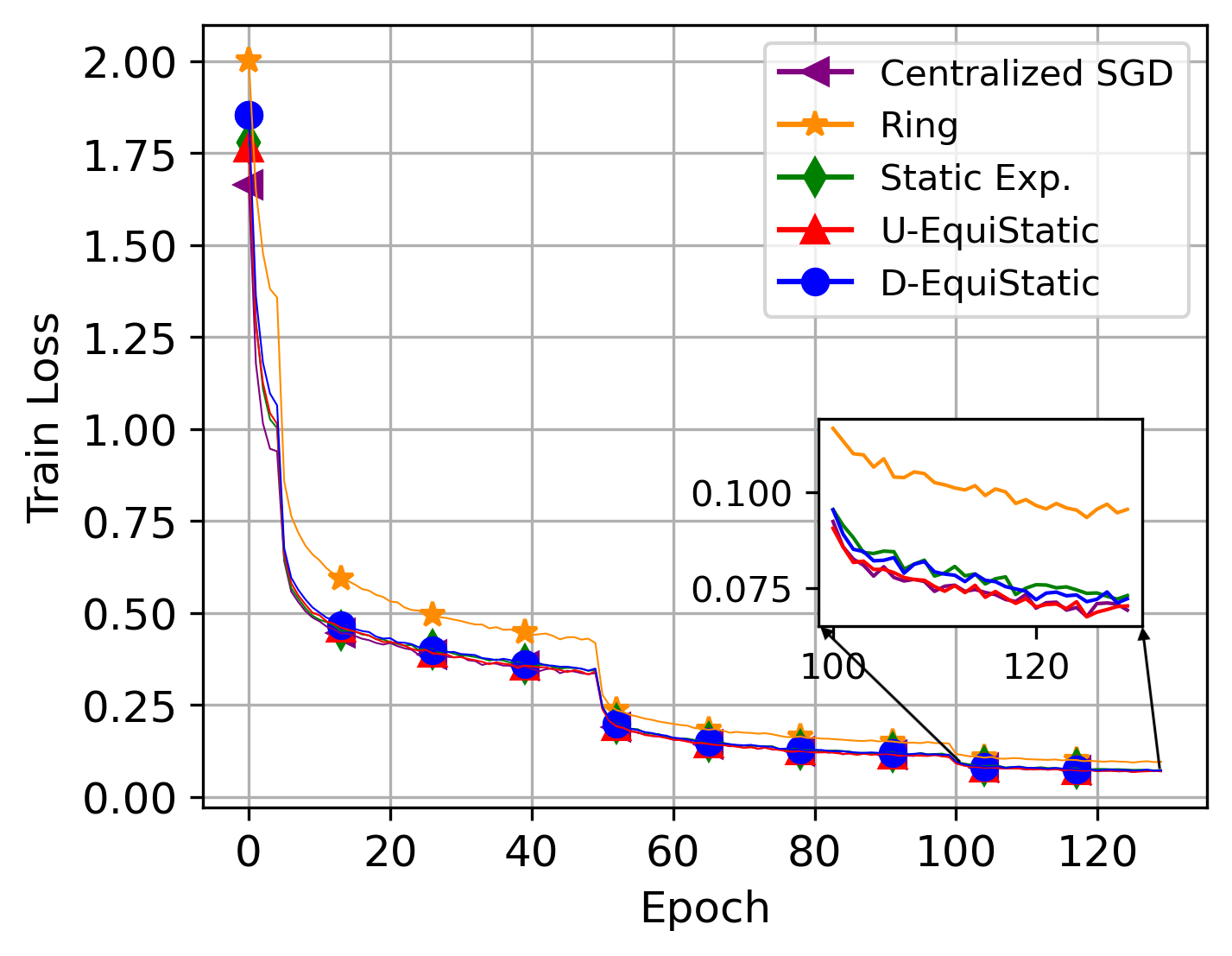}
	\quad 
	\includegraphics[scale=0.5]{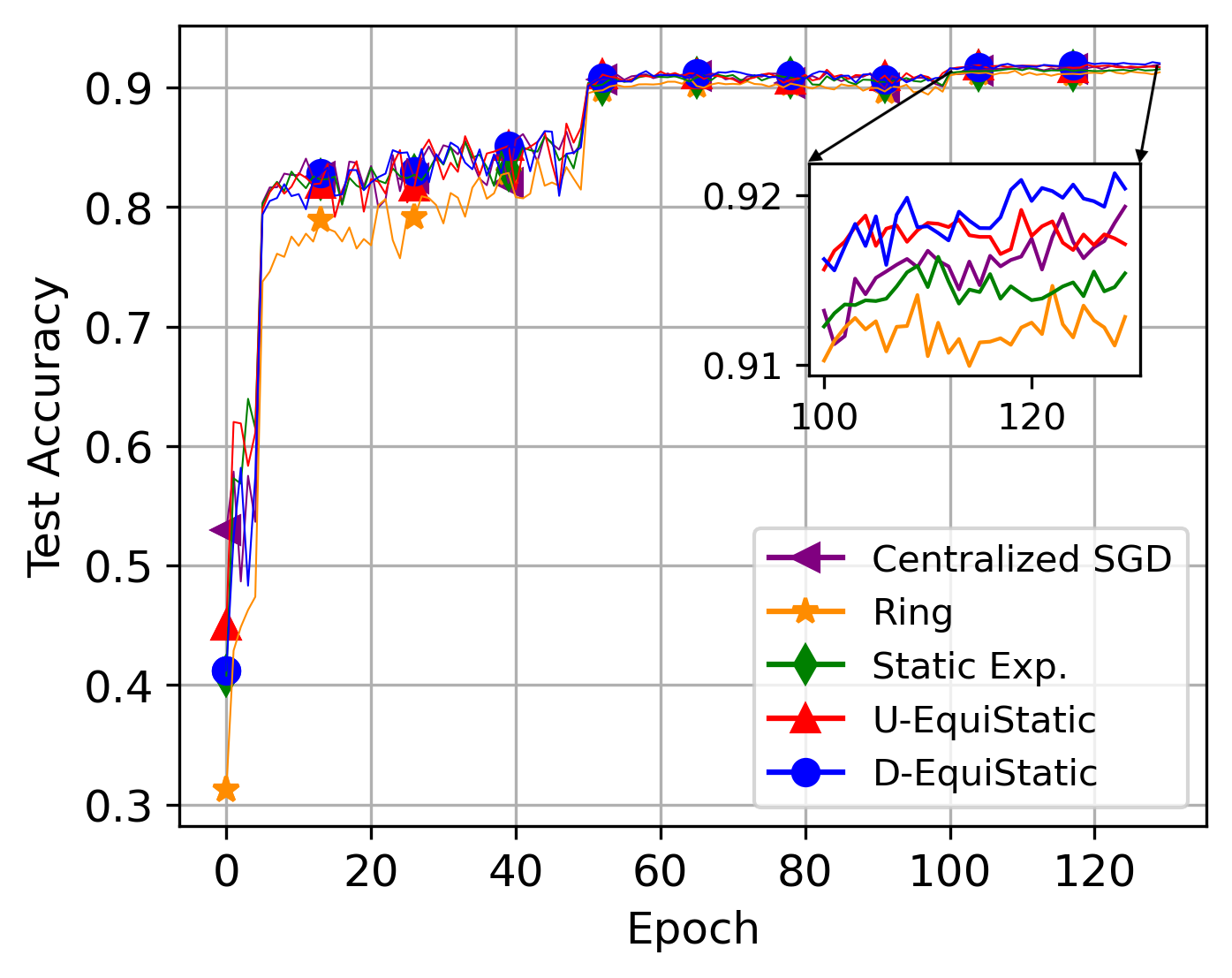}
	\caption{\small Train loss and test accuracy comparisons among different topologies for ResNet-20 on CIFAR-10.}
	\label{fig:cifar-10}
\end{figure}

\textbf{DSGD with EquiTopo: deep learning.} We consider the image classification task with ResNet-20 model \cite{he2016deep} over the CIFAR-10 dataset \cite{krizhevsky2009learning}. We utilize BlueFog \cite{ying2021bluefog} to support decentralized communication and topology setting in a cluster of 17 Tesla P100 GPUs. Fig.~\ref{fig:cifar-10} illustrates how D/U-EquiStatic compares with static exp., ring, and centralized SGD in training loss and test accuracy. It is observed that D/U-EquiStatic has strong performance. They achieve 
competitive training losses to centralized SGD but slightly better test accuracy. Meanwhile, they also outperform static exponential graphs in test accuracy by a visible margin (D-EquiStatic: 92$\%$, U-EquiStatic: 91.7$\%$, Static Exp.: 91.5$\%$). Experiments with EquiDyn topologies 
and results on MNIST dataset \cite{lecun2010mnist} are in Appendix~\ref{Ap:Exp}.

\section{Conclusion}
This paper proposes EquiTopo graphs that achieve a  state-of-the-art balance between the maximum degree and consensus rate. The EquiStatic graphs are with $\Theta(\ln(n))$ degrees and $n$-independent consensus rates, while their one-peer variants, EquiDyn, maintain roughly the same consensus rates with a degree at most $1$. EquiTopo enables decentralized learning with light communication and fast convergence.

\section*{Acknowledgement}
We thank the anonymous reviewers for suggesting updated literature on E.-R. graphs.  
Ming Yan was partially supported by the NSF award DMS-2012439. 
Lei Shi was partially supported by Shanghai Science and Technology Program under Project No. 21JC1400600 and No. 20JC1412700
and National Natural Science Foundation of China (NSFC) under Grant No. 12171093. 
Kexin Jin was supported by Alibaba Research Internship Program.

{
\bibliographystyle{ieee_fullname}
\bibliography{references}
}

\newpage
\appendix

\title{Appendix for Communication-Efficient Topologies for Decentralized Learning with $\Om(1)$ Consensus Rate}

\maketitle

\section{Directed EquiTopo Graphs}
\label{Ap:DirEqu}

\subsection{Construction of a  D-EquiStatic graph}
\label{Ap:AlgD-EquiStati}

A practical method to construct a D-EquiStatic weight matrix $\W$ is provided in Alg.~\ref{Alg:D-EquiRand}.
We should mention that the ``while" loop in the algorithm is adopted to guarantee $\norm{\Pin \W}_2 \le \rho$.

\begin{algorithm} 
\caption{A practical method for D-EquiStatic weight matrix generation}
\label{Alg:D-EquiRand}
\KwIn{Network size $n$; desired consensus rate $\rho \in (0,1)$; probability $p $ }	 \vspace{1mm}

Set $M = \ceil{\frac{8}{3\rho^2}\ln(2n/p)} $ and initialize $\WW = \II$\; \vspace{1mm}

\While {$\norm{\Pin \W}_2 > \rho$}{\vspace{1mm}

Sample $M $ i.i.d random variables $u_1, u_2, \dots, u_M$ uniformly  from $[n-1]$;

Generate basis weight matrices $\{A^{(u_i)}\}_{i=1}^M$ according to \eqref{Au};  \vspace{1mm}

Construct $\W$ by (\ref{eq:defW});
}

\KwOut{The D-EquiStatic weight matrix $\W$ and its associated basis indices $\{u_t\}_{t=1}^M$} \vspace{1mm}
\end{algorithm}

\subsection{Proof of Theorem \ref{lm-staconstr}}
\label{Ap:pflm-staconstr}

Before showing properties of $\W$ defined by (\ref{eq:defW}), 
we provide two lemmas as follows.

Referring to Theorem 1.6 of \cite{tropp2012user}, we have the following result for a sequence of random matrices.
\begin{lemma}[Matrix Bernstein]
\label{lem:MatrixBern}	
Consider a sequence of $K$ independent random $n\times n$ matrices $\{\MM_i\}_{i=1}^K$.
Assume that each random matrix satisfies
    \begin{align*} \E[\MM_i] =0, \quad {\rm and} \quad \Vert \MM_i \Vert_2 \le R~~ {\rm almost ~surely}.\end{align*}
Define 
\begin{align*}\sigma^2 := \max\Big\{\Big\Vert \sum\nolimits_{i =1}^K \E[ \MM_i \MM_i^T]\Big\Vert_2,
\Big\Vert \sum\nolimits_{i =1}^K \E[ \MM_i^T \MM_i]\Big\Vert_2 \Big\}.
\end{align*}
It holds that 
\begin{equation*}
\Prb\Big(\Big\Vert \sum\nolimits^K_{i=1} \MM_i \Big\Vert_2 \ge \delta \Big) \le 2n \exp\Big(- \frac{\delta^2/2}{\sigma^2 + R\delta/3}\Big),\quad \forall \delta \ge 0.
\end{equation*}	
\end{lemma}

\begin{lemma}
\label{norm_comp}
For any matrix $\B = [b_{ij}] \in \R^{n \times n}$, it holds that	
\begin{align*}\norm{\B}_2 \le \max \{\norm{\B}_1, \norm{\B}_{\infty}\}.\end{align*}
\end{lemma}	

\emph{Proof.}
By definition,
\begin{align*}\norm{\B}_2  = \sup_{\norm{\x}\leq 1, \norm{\y} \le 1} \x^T \B  \y.\end{align*}
Moverover, for all $\norm{\x}\leq 1,~\norm{\y} \le 1$, we have
\begin{equation*}
\begin{aligned}
(\x^T \B \y)^2 \le \Big(\sum_{i,j} |b_{ij}| x_i^2 \Big) \Big(\sum_{i,j} |b_{ij}| y_j^2 \Big) \le \norm{\B}_{1}\norm{\B}_{\infty}.
\end{aligned}
\end{equation*}		
Thus, Lemma \ref{norm_comp} holds.
$\hfill\square$

\textbf{Theorem~\ref{lm-staconstr}} 
\emph{
(Formal restatement of Theorem~\ref{lm-staconstr})
Let $\A^{(u)}$ be defined by \eqref{Au} for any $u\in [n-1]$ and the D-EquiStatic weight matrix $\W$ be constructed by \eqref{eq:defW} with $\{u_i\}_{i=1}^M$ following an independent and identical uniform distribution from $[n-1]$. For any size-independent consensus rate $\rho \in (0,1)$ and probability $p \in (0,1)$, if $M \geq \frac{8}{3\rho^2}\ln\frac{2n}{p}$ it holds with probability at least $1 - p$ that   
\begin{equation}
\nt{\Con\WW\xx} \leq \rho \nt{\Con\xx},\ \forall \xx\in \Real^n.  
\end{equation}}	
\emph{Proof.}	
Notice that
\begin{equation*}
\E[\A^{(u_i)}] = \J, ~\ \forall i \in \{1,2,\dots,M\}.
\end{equation*}	
Since each $\A^{(u_i)}$ is doubly stochastic, it follows from Lemma \ref{norm_comp} that
\begin{align*}\norm{\A^{(u_i)} - \J}_2= \norm{\Pin \A^{(u_i)}}_2 \le   \norm{\Pin}_2  \norm{\A^{(u_i)}}_2 \le 1.\end{align*} 
Consequently,
\begin{align*}
\sum_{i=1}^M \E \Big[\norm{(\A^{(u_i)} - \J)(\A^{(u_i)} - \J)^T}_2 \Big] \le
\sum_{i=1}^M  \E \Big[\norm {\Pin \A^{(u_i)}}_2^2 \Big] \le M.
\end{align*}

Analogously,
$\sum_{i=1}^M \E  \big[\norm{(\A^{(u_i)} - \J)^T(\A^{(u_i)} - \J)}_2 \big]\leq M.$
By Lemma \ref{lem:MatrixBern},
\begin{align*}
\Prb \big(\norm{\W-\J}_2 \ge \rho \big)
=&\Prb\Big(
\norm{\sum\nolimits_{i=1}^M (\A^{(u_i)} -\J)  }_2 \ge M\rho
\Big) \\  
\leq& 2n\exp\Big(- \frac{M^2\rho^2  / 2}{M + M\rho/3} \Big) 
\leq 2n\exp\Big(- \frac{M^2\rho^2  / 2}{M + M/3} \Big) \leq  p,
\end{align*}	
i.e.,
\begin{equation*}
\Prb(\norm{ \Pin \W}_2 \le \rho) \ge 1-p.
\end{equation*}	
Note that $\Pin \W  = \Pin \W \Pin $. If $\nt{\Con\WW}_2\leq \rho $, then
\begin{equation*}
\begin{aligned}
\norm{ \Pin \W \x}^2 =  \norm{ \Pin \W \Pin \x}^2 
\le  \norm{ \Pin \W}_2^2 \norm{\Pin \x}^2\leq \rho^2 \norm{\Pin \x}^2.
\end{aligned}
\end{equation*}
Therefore, the conclusion holds.
$\hfill\square$

The relation $M \geq \frac{8}{3\rho^2}\ln(2n/p)$ is required for theoretical analysis, and it is very conservative. In practice, we can set $M$ to be far less than $\frac{8}{3\rho^2}\ln(2n/p)$ and repeat the process described in the formal version of Theorem~\ref{lm-staconstr} until we find a desirable $\WW$ (see Alg.~\ref{Alg:D-EquiRand}). 
In addition, 
the verification condition $\norm{\Pin \W}_2 \leq \rho$ in Alg.~\ref{Alg:D-EquiRand} can also be dropped in implementations so that we only conduct the ``while'' loop once. 
We find that these relaxations can still achieve $\W$ with an empirically fast consensus rate (see the illustration in the experiments). 

{
\begin{remark}
    We have much flexibility in the choice of $p$, such as $p =  1/2$ or $p = 1/n$. 
    If $p \in (1/{\rm poly\pr{n}}, 1)$, the corresponding value of $M$ will only differ in constants. 
\end{remark}
}

\subsection{Proof of Theorem \ref{the:od-EquiRand}}
\label{Ap:OD-EquiRand}

Due to $\Pin \A^{(v_t)} = \A^{(v_t)}  \Pin = \Pin  \A^{(v_t)}  \Pin$,
it follows from Alg. \ref{Alg:OD-EquiRand} that
\begin{equation*}
\begin{aligned}
\E\Big[ &\norm{\Pin \W^{(t)} \x}^2 \Big] = \E \Big[\norm{\Pin \Big((1 - \eta) \I + \eta \A^{(v_t)} \Big) \x}^2 \Big] \\
&\le (1 - \eta)^2 \norm{\Pin \x}^2 + 2\eta(1 - \eta)(\Pin \x)^T \E \big[\Pin\A^{(v_t)} \big] \Pin \x
+ \eta^2  \E \Big[\norm{\A^{(v_t)}}_2^2 \Big] \norm{\Pin \x}^2 \\
\end{aligned}
\end{equation*}	
Notice that $ \E \big[\A^{(v_t)} \big] =  \W$
and $\norm{\Pin \W \y} \le \rho \norm{\Pin \y}$, $\forall \y \in \R^n$. Therefore,
\begin{equation*}
\begin{aligned}
\E\Big[\norm{\Pin \W^{(t)} \x}^2 \Big]
&\le \big( (1-\eta)^2 + 2\eta \rho (1-\eta) +\eta^2 \big) \norm{\Pin \x}^2\\
&=\big(1-2\eta(1-\eta)(1-\rho)\big) \norm{\Pin \x}^2.
\end{aligned}
\end{equation*}	
The proof is completed.
$\hfill\square$

\section{One-Peer Undirected EquiTopo Graphs (OU-EquiDyn)}
\label{Ap:OU-EquiRand}

\subsection{Illustration for basis graphs of OU-EquiDyn}
\label{Ap:basisouequidyn}

\begin{figure}[h]
\centering
\includegraphics[scale=0.39]{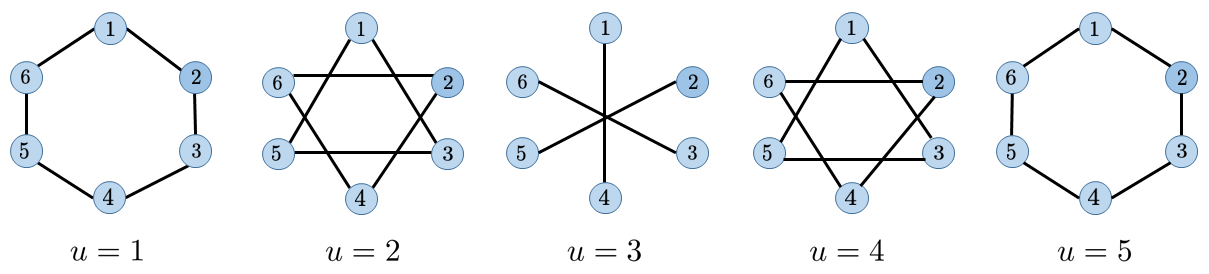}
\caption{\small  Undirected graphs generated by  $\widehat{\A}^{(u)} = \frac{1}{2}(\A^{(u)} + [\A^{(u)}]^T)$.}
\label{fig:ou-equirandb}
\end{figure}

The associated graph of 
$\widehat{\A}^{(u)} = \frac{1}{2}(\A^{(u)} + [\A^{(u)}]^T)$ with $n = 6$  
 are given in Fig. \ref{fig:ou-equirandb}. Clearly, there exist non-one-peer graphs.

\begin{figure}[h]
\centering
\includegraphics[scale=0.39]{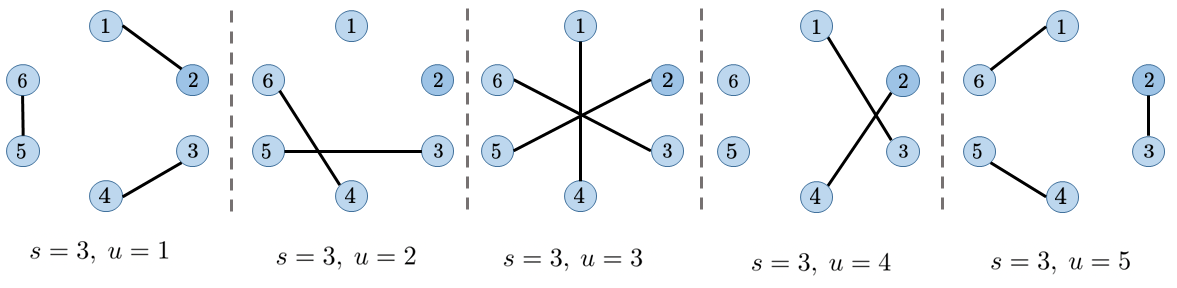}
\caption{\small One-peer undirected graphs generated in Alg. \ref{Alg:CuMN} with  $s=3$ and $u=1,\dots,5$.}
\label{fig:ou-oneundirsu}
\end{figure}

Consider $n = 6$. The  OU-EquiDyn graphs generated by Alg. \ref{Alg:CuMN} (when $s = 3$) are presented in
Fig. \ref{fig:ou-oneundirsu}.    

\subsection{Node version of Alg.~\ref{Alg:CuMN}}
\label{Ap:EqAgetAlg}

From the node's perspective, an equivalent version of 
Alg.~\ref{Alg:CuMN} is presented in 
Alg.~\ref{Alg:CuMNagentpersp1}. 
In the remainder of Appendix~\ref{Ap:OU-EquiRand},  
we denote
 $\floor{x}$ as the largest   integer no greater than $x$.  
We also denote the traditional $\mathrm{mod}$ operator as 
\eq{
    \modk{a}{b} = a - b \cdot \floor{\frac{a}{b}}.   
}

\begin{algorithm}[t]
\caption{OU-EquiDyn weight matrix generation at iteration $t$ (from nodes' perspective) }
\label{Alg:CuMNagentpersp1}  
\KwIn{$\eta \in (0,1)$; basis index 
$\{u_i, -u_i\}_{i=1}^M$ from a weight matrix $\widetilde \W \in \R^{n\times n}$ of form \eqref{undirsta};} 
\vspace{1mm}

Initialize $\AA = [a_{ij}] = \II$\; \vspace{1mm}  

\For{node $i=1$ to $n$ (in parallel) }{\vspace{1mm}  

Pick $v_t$ from $\{u_i, -u_i\}_{i=1}^M$ 
and  $s_t \in [n]$ uniformly at random using the common random seed; \vspace{1mm}

\If{$v_t \leq n/2  $ }{
    $q = v_t $\; \vspace{1mm}
    $k(i) = \modk{i - s_t}{n} $\; \vspace{1mm}  
}\Else{
    $q = n - v_t$\; \vspace{1mm}
    $k(i) = \modk{i - s_t + q  }{n} $\; \vspace{1mm} 
}

$r(i) = \modk{k(i)}{q} $\;\vspace{1mm}
$d(i) = \floor{{k(i)}/q} $\; \vspace{1mm}  

\If{$\floor{(n - 1 - r(i))/q}$ is odd or $d(i) < \floor{(n-1-r(i))/q}$ }{\vspace{1mm}
\If{$d(i)$ is even}{ \vspace{1mm}
 $j = \nmod{(i + q)}$\;  
}\Else{ \vspace{1mm}
 $j = \nmod{(i - q)}$\;\vspace{1mm}  
}
$a_{ij}  = (n-1)/n $\; \vspace{1mm} 

$a_{ii}  = 1/n $\; 
}
}

\KwOut{$\widetilde{\W}^{(t)} = (1-\eta)\I + \eta \A$}

\end{algorithm}

If $v_t = n/2$, every node $i$ is connected to node $\nmod{(i+n/2)}$ for Alg.~\ref{Alg:CuMN},
which is the same as Alg.~\ref{Alg:CuMNagentpersp1}. 

If $v_t < \frac{n}{2}$, then $q = v_t$. 
In Alg.~\ref{Alg:CuMN}, nodes $\dr{s_t, s_t + 1, \cdots, s_t + q - 1} \ \mathrm{mod} \ n $ are connected to $\dr{s_t + q, s_t + q + 1, \cdots, s_t + 2q - 1} \ \mathrm{mod} \ n$, respectively. 
If $j\in \dr{s_t + q, s_t + q + 1, \cdots, s_t + 2q - 1}  \ \mathrm{mod} \ n$, then new edges cannot be added because they have been connected. 
Similar process starts from connecting node $\nmod{(s_t + 2q)}$ with node $\nmod{(s_t + 3q)}$, and as a result, 
for $v_t < \frac{n}{2}$, Alg.~\ref{Alg:CuMN} can be interpreted 
as follows:
divide $[n]$ into $q$ disjoint subsets:  
\eq{
    C_{\ell} = \dr{i \in [n]: i = \nmod{(s_t+\ell + d\cdot q)},\  d\in \mathbb{Z} },\ 0\leq \ell < q.    
}
Equivalently,
\eq{ 
    C_{\ell} = \dr{i \in [n]: i = \nmod{(s_t+\ell + d\cdot q)},\  0\leq d\leq \floor{\frac{n - 1 - \ell }{q}} }. 
}
For node $i$, we define $k(i) = \modk{i - s_t}{n}$, $d(i) = \floor{k(i)/q}$ 
and $r(i) = \modk{k(i)}{q}$.   
Then, $i\in C_{r(i)}$ and $i = \nmod{(s_t + r(i) + d(i) \cdot  q)}$.

In each $C_{\ell}$, $s_t + \ell$ is connected with $s_t + \ell + q$, $s_t + \ell + 2q$ is connected with $s_t + \ell + 3q$, \dots. 
Thus, if $|C_\ell|$ is even, every node in $C_\ell$ has a neighbor. 
If $|C_\ell | $ is odd (equivalently, $\floor{\frac{n - 1 - \ell}{q} }$ is odd), the node $\ell + q\cdot \floor{\frac{n - 1 - \ell}{q}}  $ is idle and the others has neighbors. 

Define 
\begin{equation*}
C'_{\ell} = 
\left\{
\begin{aligned}
& \  C_{\ell},\qquad \mbox{if $\floor{\frac{n - 1 - \ell}{q } } $ is odd} \vspace{1mm} \\        
& \dr{i\in [n]: i = \nmod{\pr{s_t + \ell + d\cdot q}},\  0\leq d < \Big\lfloor \frac{n-1-\ell}{q}\Big\rfloor } , \qquad  \mbox{otherwise.} 
\end{aligned}  
\right.
\end{equation*}
Then,  node $i$ has a neighbor if and only if it is in the set
$
   C'_{r(i)}.    
$
In addition, for node $i$ in $C'_{r(i)} $, it is connected to $\nmod{(i+q)}$  if $d(i) $ is even; and connected to $\nmod{(i - q)}$ if $d(i)$ is odd.   

In Alg.~\ref{Alg:CuMNagentpersp1}, we compute $r(i)$ and $d(i)$ firstly, 
and then check whether node $i$ is in $C'_{r(i)}$. 
If $i\in C'_{r(i)}$ and $d(i)$ is even (odd)
, then it is connected to $\nmod{(i+q)}$ 
($\nmod{(i - q)}$). 
Otherwise, node $i$ is idle, i.e., $a_{ii} = 1$. 
This yields the equivalence between Alg.\ref{Alg:CuMN} and Alg.~\ref{Alg:CuMNagentpersp1} for the case of $v_t < n/2$. 

If  $v_t > n/2$, let $q = n - v_t < n/2$.   
Then the nodes in $\dr{s_t, s_t + 1, \cdots, s_t + q - 1 } \mod n $ are connected to the nodes $\dr{s_t + v_t, s_t + 1 + v_t, \cdots, s_t + q - 1 + v_t } \mod n$ respectively.  
Note that $\nmod{(i + v_t) }$ is equivalent to $\nmod{(i - q)}$. Then, equivalently, nodes $\dr{s_t, s_t + 1, \cdots, s_t + q - 1}  \mod  n $ are  connected with $\dr{s_t - q, s_t - q + 1, \cdots, s_t - 1 }  \mod  n $, respectively.  
Similar process starts from connecting node $\nmod{(s_t + 2q)}$ with node $\nmod{(s_t + q)}$.   
Then, the undirected graph generated with the starting point $s_t $ and the label difference $v_t > n/2$ is equivalent to the undirected graph generated with the starting point $\nmod{(s_t - q)} $ and the label difference $q$. 
So by setting $k(i) = \modk{i - (s_t - q)}{n}$, the proof follows by similar arguments for the case $v_t < n/2$.

\subsection{Proof of Theorem \ref{the:ou_equirand}}
\label{Ap:OUThe}

We first provide the following three lemmas.
In the remainder of Appendix~\ref{Ap:OU-EquiRand}, for any matrix $\AA = [a_{ij}]\in \Real^{n\times n}$, we denote its edge set as
\eq{
\Eset{\AA} = \dr{(i,j)\in [n]\times [n]:\ a_{ij} > 0,\ i\neq j}.  
}

Denote the matrix $\AA $ generated at the $t$-th iteration of Alg.~\ref{Alg:CuMN} by $\widetilde{\AA}^{(v_t)} = [\widetilde{a}_{ij}^{(v_t)}]\in \Real^{n\times n}$.  
Note that $\widetilde{\AA}^{(v_t)}$ is also stochastic even when $v_t$ is given since $s_t$ is randomly chosen from $[n]$.

\begin{lemma}
\label{lem:EsetA1}
For any $n \geq 2$, it holds for 
$\widetilde {\A}^{(v_t)}$ defined by Alg. \ref{Alg:CuMN} that 
it holds  that  
\begin{align*}\E\big[|\Eset{\widetilde{\A}^{(v_t)}}|\big]   \geq \frac{2n}{3}.\end{align*}
\end{lemma}
\emph{Proof.}
Because $|\Eset{\widetilde{\AA}^{(v_t)}}|$ is invariant with respect to $s_t$, it suffices to prove $|\Eset{\widetilde{\A}^{(v_t)}}| \geq 2n/3$ for $s_t=1$.  

For $v_t \leq n/2$, we define  
    $m = \floor{n/(2v_t)}$ and $r = \modk{n}{2v_t} $, then, $m \geq 1$.   
Notice that node $i$ is connected with $i + v_t$ for any $i\in \dr{\ell + 2d v_t: \ 1\leq \ell \leq v_t,\ 0\leq d < m }$.   

If $r \leq v_t $, then the last $r$ nodes are idle. 
As $m\geq 1$, we have $n = 2v_t m + r \geq 2v_t + r \geq 3r. $ 
Thus, 
\eq{
    |\Eset{\widetilde{\A}^{(v_t)}}| \geq n - r \geq n - \frac{1}{3}n = \frac{2}{3}n.  
} 
If $r > v_t$, then node $i$ in $\dr{2 m v_t + \ell: \ 1\leq \ell\leq r - v_t }$ is connected to $i + v_t$. 
As a result, only the nodes in $\dr{2 m v_t + \ell: \  r - v_t + 1\leq \ell \leq v_t}$ are idle, i.e.,  $2v_t - r $ nodes are idle.  
We have  
$ n = 2 m v_t + r > 3v_t  $ from $m\geq 1$ and $r > v_t $.   
Consequently,
\eq{
     |\Eset{\widetilde{\A}^{(v_t)}}| \geq n - (2v_t - r) \geq n - v_t \geq \frac{2}{3}n. 
}

We have shown $|\Eset{\widetilde{\A}^{(v_t)}}| \geq \frac{2}{3}n$ for $v_t \leq n/2$.

For $v_t > n/2$, recall that we have shown in Appendix~\ref{Ap:EqAgetAlg} that the undirected graph generated with label difference $v_t$ and starting point $s_t$  equals the undirected graph generated with label difference $q = n - v_t < n/2$ and starting point $\nmod{(s_t - q)}$.  Since the number of edges is invariant with the starting point, the case $v_t > n/2$ has been reduced to $v_t < n/2$.   
This completes the proof.
$\hfill\square $

\begin{lemma}\label{lem:symA}
For any symmetric matrix $\B = [b_{ij}]\in \R^{n\times n}$, if $\B \mathds{1}_n = \mathbf{0}_n$, then
\begin{align*}\x^T\B \x = - \frac{1}{2} \sum_{i, j} b_{ij}(x_i-x_j)^2,\quad \forall \xx\in \Real^n.  \end{align*}
\end{lemma}

\emph{Proof.}
The $i$-th entry of $\B\x$ is
\begin{align*}[\B\x]_i = \sum_j b_{ij}x_j = b_{ii}x_i + \sum_{j:j\neq i}b_{ij}x_j = \sum_{j}b_{ij}(x_j - x_i).\end{align*}
Hence,
\begin{align*}\x^T\B\x = \sum_{i,j} b_{ij}x_i(x_j - x_i).\end{align*}
Due to $\B^T = \B$, we have
\begin{align*}\x^T\B\x = \sum_{i,j} b_{ij}x_j(x_i - x_j).\end{align*}
Averaging the above equations yields the result.
$\hfill\square$

\begin{lemma}\label{lem:expecW1}
For any $n \geq 2$, it holds for 
$\widetilde {\A}^{(v_t)}$ 
that 
\begin{equation*}
\mathbb E[\widetilde {\A}^{(v_t)}] \preceq  \frac{1}{3} \I + \frac{1}{3} \Big({\A}^{(v_t)}+[{\A}^{(v_t)}]^T)\Big).
\end{equation*}	
\end{lemma}

\emph{Proof.}
If $v_t = n/2$, then $\widetilde{\A}^{(v_t)} = \AA^{(v_t)} = [\AA^{(v_t)}]\tp $. 
By Lemma~\ref{lem:symA}, 
$\xx\tp \pr{\II - \AA^{(v_t)}} \xx \geq 0 $, i.e., $\AA^{(v_t)} \preceq \II$. 
Thus, 
\begin{align*}\widetilde{\AA}^{(v_t)} = \AA^{(v_t)} \preceq \frac{1}{3}\II + \frac{2}{3}\pr{\AA^{(v_t)} + [\AA^{(v_t)}]\tp}.  
\end{align*}

Consider  $v_t \neq n/2$.   
Notice that $\mathcal{E}(\widetilde{\A}^{(v_t)})\subset \mathcal{E}\big(\A^{(v_t)} + [\A^{(v_t)}]^T\big)$ and $|\mathcal{E}(\A^{(v_t)} + [\A^{(v_t)}]^T)| \leq 2n$. 
By Lemma~\ref{lem:EsetA1} and the fact that $s_t $ is from the uniform distribution over $[n]$, it holds 
for any $(i, j)\in \mathcal{E}(\A^{(v_t)}+ [\A^{(v_t)}]^T)$ that
\begin{align*}\Prb{[(i, j)\in \mathcal{E}(\widetilde{\A}^{(v_t)})]} \geq \frac{1}{2n} \E[|\mathcal{E}(\widetilde{\A}^{(v_t)})|] \geq \frac{1}{3}.  
\end{align*}

By the construction of $\widetilde{\A}^{(v_t)}$ in Alg.~\ref{Alg:CuMN} and $\A^{(v_t)}$ in~\eqref{Au}, the non-diagonal and non-zero entries are $\frac{n-1}{n}$. It follows from 
Lemma \ref{lem:symA} that 
for any $\x\in \R^n$, we have
\begin{align*}
\x^T \E[\mathbf{I} - \widetilde{\A}^{(v_t)}]\x &= \frac{n-1}{2n}\E\Bigg[ \sum_{(i, j)\in \mathcal{E}(\tilde{\A}^{(v_t)})} (x_i - x_j)^2 \Bigg]\\
&=\frac{n-1}{2n}\sum_{(i, j)\in \mathcal{E}(\A^{(v_t)} + [\A^{(v_t)}]^T)} \Prb[(i, j)\in \mathcal{E}(\widetilde{\A}^{(v_t)})] (x_i - x_j)^2\\
&\geq \frac{n-1}{6n} \sum_{(i, j)\in \mathcal{E}(\A^{(v_t)} + [\A^{(v_t)}]^T)} (x_i - x_j)^2\\
&= \frac{1}{3} \x^T \big(2\mathbf{I} - \A^{(v_t)} - [\A^{(v_t)}]^T\big)\x.
\end{align*}
Rearranging the terms, we derive 
\eq{
\E\br{\Atil^{\pr{v_t}}} \pleq \frac{1}{3}\II + \frac{1}{3}\pr{\A^{\pr{v_t}} + [\A^{\pr{v_t}}]\tp}.        
}
This completes the proof.
$\hfill\square$

\textbf{Proof of Theorem \ref{the:ou_equirand}}
It follows from Lemma \ref{lem:expecW1} that
\begin{align*}\E\Big[\widetilde{\A}^{(v_t)} \Big]
\preceq \frac{1}{3} \I + \frac{1}{3} 
\E \Big[ \A^{(v_t)}+[\A^{(v_t)}]^T \Big]
=\frac{1}{3} \I + \frac{2}{3} \widetilde{\W}.\end{align*}
Consequently,
\begin{equation*}
\begin{aligned}
\E \Big[&\normB{\Pin \widetilde{\W}^{(t)} 
\x}^2\Big]
= \E \Big[\normB{\Pin \Big((1 - \eta) \I + \eta \widetilde{\A}^{(v_t)} \Big) \x}^2 \Big] \\
&\le (1 - \eta)^2 \norm{\Pin \x}^2 + 2\eta(1 - \eta)(\Pin \x)^T \E \big[\widetilde{\A}^{(v_t)} \big] \Pin \x
+ \eta^2  \E \Big[\normB{\widetilde{\A}^{(v_t)}}_2^2 \Big] \norm{\Pin \x}^2.
\end{aligned}
\end{equation*}	
Combining the above two inequalities, we derive
\begin{equation*}
\begin{aligned}
\E\Big[\normB{\Pin \widetilde{\W}^{(t)} \x}^2 \Big] 
&\le \big( (1-\eta)^2 + \frac{2}{3}\eta (1-\eta) + \frac{4}{3}\eta \rho (1-\eta) +\eta^2 \big) \norm{\Pin \x}^2\\
&=\big(1-\frac{4}{3}\eta(1-\eta)(1-\rho)\big) \norm{\Pin \x}^2.
\end{aligned}
\end{equation*}	
Thus, the proof is completed.
$\hfill\square$

\subsection{An alternative construction of OU-EquiDyn}
\label{Ap:EOUA}
Alg.~\ref{Alg:CuMNagentpersp1EuclideanAlgorithm} provides a different way to construct one-peer undirected graphs which achieve a similar consensus rate
as the graphs generated by Alg.~\ref{Alg:CuMN} but with a different structure. 

\begin{algorithm}[t]
\caption{Alternative OU-EquiDyn weight matrix generation at iteration $t$ (from nodes' perspective)}
\label{Alg:CuMNagentpersp1EuclideanAlgorithm}
\KwIn{$\eta \in (0,1)$;
basis index 
$\{u_i, -u_i\}_{i=1}^M$ from a weight matrix $\widetilde \W \in \R^{n\times n}$ of form \ref{undirsta};} 
\vspace{1mm}

Initialize $\AA = [a_{ij}] = \II$\; \vspace{1mm}  

\For{node $i=1$ to $n$ (in parallel) }{\vspace{1mm}

Pick $v_t$ from $\{u_i, -u_i\}_{i=1}^M$ 
and  $s_t \in [n]$ uniformly at random using the common random seed; \vspace{1mm}

Compute $d = \gcd(v_t, n) $ and find $1\leq b \leq n/d - 1$ such that $\modk{b\cdot v_t}{n} = d$ by Euclidean algorithm\algl

Set $\ntil = n/d$,  and
$m(i)=\modk{\floor{(i-s_t)/d} \cdot b}{\ntil}$; \vspace{1mm}  

\If{$\ntil $ is even or $m(i) < \ntil - 1 $ }{\vspace{1mm}
\If{$m $ is even}{\vspace{1mm}
$j = \nmod{(i + v_t)}   $\;\vspace{1mm}  
}
\Else{\vspace{1mm}
$j = \nmod{(i - v_t)}  $\;\vspace{1mm}  
}
$a_{ij} = (n-1)/n $\; \vspace{1mm} 

$a_{ii} = 1/n $\;\vspace{1mm} 
}
}

\KwOut{$\widetilde{\W}^{(t)} = (1-\eta)\I + \eta \A$}
\end{algorithm}

We denote the matrix $\AA $ generated at the $t$-th iteration of Alg.~\ref{Alg:CuMNagentpersp1EuclideanAlgorithm} by $\overline{\AA}^{(v_t)} = [\overline{a}^{(v_t)}_{ij}]\in \Real^{n\times n}$. Next, we explain the motivation of Alg.~\ref{Alg:CuMNagentpersp1EuclideanAlgorithm}. 

Denote $\gcd(a, b)$ as the greastest common divisor of $a$ and $b$.  
Let $d = \gcd(v_t, n)$ and $\vtil_t = v_t/d$, $\ntil = n/d$. 
Then, $\vtil_t$ and $\ntil$ are coprime.   

Firstly, we divide $[n]$ into $d$ disjoint subsets:
\eq{
    \Col_{\ell} = \dr{i\in [n]:\ \modk{i - s_t}{ d} = \ell },\ 0\leq \ell < d.      
} 
Clearly, $\abs{\Col_{\ell}} = n/d = \ntil $, $\forall 0\leq \ell < d$. 

We claim that 
\eql{\label{eq:ColdefBezout1}}{
    \Col_{\ell} = \dr{\nmod{s_t + \ell + m v_t  }:\ 0\leq m < \ntil },\ 0\leq \ell < d.  
}
To proof the claim, we denote the RHS of~\eqref{eq:ColdefBezout1} by $\Ctil_{\ell}$.
Since $v_t$ can be divided evenly by $d$, for any $i\in \Ctil_{\ell} $, it satisfies $\modk{i}{d} = \modk{(s_t + \ell)}{d} $.
Combining with the fact that $\Ctil_{\ell} \subset [n]$, we have $\Ctil_{\ell} \subset \Col_{\ell}$. 

Then, since $\ntil$ and $\vtil$ are coprime, for any $0\leq m_1 <m_2 <  \ntil $, $\modk{m_1\vtil_t}{\ntil} \neq \modk{m_2\vtil_t}{\ntil}$. 
Then, $\modk{m_1 v_t}{n} \neq \modk{m_2 v_t}{n}$. 
Thus, $\abs{\Ctil_{\ell}} = \ntil = \abs{\Col_{\ell}}$.   
Combining with $\Ctil_{\ell} \subset \Col_{\ell}$, we have $\Col_{\ell} = \Ctil_{\ell}$. 

The above analysis also implies that for each $i\in \Col_{\ell}$, there exists a unique $0\leq m(i) < d$ such that $i = \nmod{(s_t + \ell + m(i)v_t)} $.    

By~\eqref{eq:ColdefBezout1}, a natural way to construct one-peer undirected graphs is: 
in each $\Col_{\ell}$, connect $s_t + \ell $ with $s_t + \ell + v_t$, connect $s_t + \ell + 2v_t$ with $s_t + \ell + 3 v_t $, \dots. 
Equivalently, $i$ is connected with $\nmod{(i + v_t)}$ if $m(i)$ is even and connected with $\nmod{(i - v_t)}$ if $m(i)$ is odd.  

In this way, 
 if $\ntil $ is even, every node in  $C_{\ell}$ has a neighbor, $\forall 0\leq \ell < d$, i.e., every node in $[n]$ has a neighbor.  
Equivalently, a node $i\in \Col_{\ell}$ has a neighbor if it is in the set
\begin{equation}\label{eq:Colellnotidle1}
    \Col_{\ell}' = \left\{
    \begin{aligned}
        & \ C_{\ell},\qquad \mbox{if $\ntil $ is even} \\  
        & \dr{\nmod{(s_t + \ell + mv_t)}:\ 0\leq m < \ntil - 1 },\qquad \mbox{if $\ntil $ is odd}
    \end{aligned}
    \right.
\end{equation} 

To give a practical way of the process described above from each node's perspective, for each node $i\in \Col_{\ell}$, we provide a more efficient way to compute the unique $0\leq m(i) < \ntil  $ such that $i = \nmod{(s_t + \ell + m(i) v_t) } $. 
Firstly, for each $i\in \Col_{\ell}$, since $\ell = \modk{i - s_t}{d}$, we have     
\eq{
    i = \nmod{(s_t + \ell + \floor{(i - s_t)/d} \cdot d)}.      
}
Then, by B$\acute{e}$zout's theorem, since $\ntil$ and $\vtil_t$ are coprime, 
there exist integers $1 \leq b \leq \ntil - 1 $ and $b'  $ such that 
$b \cdot v_t + b' \cdot n = d   $.  
The pair $(b, b')$ can be computed by the Euclidean algorithm in $\mathcal{O}(\ln(n))$ time.   
Then,  
\eq{
     \floor{(i - s_t)/d}\cdot d = \floor{(i - s_t) / d } \cdot (b v_t + b' n).     
}
Define $m(i) = \modk{\floor{(i - s_t)/ d } \cdot b }{\ntil} $, then, we have 
$\nmod{(\floor{(i - s_t)/d}\cdot d)} = \nmod{m(i) v_t},   $ 
Thus,  
$i = \nmod{(s_t + \ell + m(i) v_t)}.  $   

In Alg.~\ref{Alg:CuMNagentpersp1EuclideanAlgorithm}, each node $i$ computes $\ntil$ and its $m(i)$ firstly. 
If $\ntil$ is even, every node has a neighbor. 
If $\ntil$ is odd but $m(i) < \ntil - 1 $, then $i\in \Col_{\ell}'$ where $\ell = \mod(i-s_t, d)$, i.e., $i$ also has a neighbor.   
In this way,  each node can determine whether it has a neighbor in this iteration.  
If node $i$ has a neighbor, as we have described above, it is connected with $\nmod{(i+v_t)}$ if $m(i)$ is even and $\nmod{(i - v_t)}$ otherwise.  

The following lemma is used to prove Lemma~\ref{lem:expecW12}.   

\begin{lemma}
\label{lem:EsetA12}
Let $d = \gcd(v_t, n)$, then  
\begin{align*} \E\big[|\Eset{  \Aol}|\big] = (2d) \cdot\floor{\frac{n}{2d}}  \geq \frac{2n}{3}.   \end{align*}
\end{lemma}
\emph{Proof.}
Define $\Col_{\ell}$ and $\Col_{\ell}'$ as in~\eqref{eq:ColdefBezout1} and~\eqref{eq:Colellnotidle1}.  
If $n$ can be divided by $2d$ evenly, i.e., $\abs{\Col_{\ell}} = \ntil$ is even for any $0\leq \ell < d$. 
Then, every node has a neighbor, i.e., 
$\abs{\Eset{\Aol}} = n = (2d) \cdot\floor{\frac{n}{2d}}. $

If $n$ cannot be divided by $2d$ evenly, 
by~\eqref{eq:Colellnotidle1},  
in each $\Col_{\ell}$, there is one node idle in this iteration.  
Then, we also have  
$\abs{\Eset{\Aol}} = n - d = (2d)\cdot \floor{\frac{n}{2d}}.   $ 
As $d = \gcd(v_t, n)$, we have $n = (2k+1)d$, with $k\in \mathbb{Z}$. 
Since $v_t$ can be divided by $d$ evenly and $v_t \leq n-1 $, we have $d < n$. 
Thus, $k\geq 1$, i.e., $n \geq 3d$.  
Then,  
$\abs{\Eset{\Aol}} = n - d \geq n - \frac{n}{3} = \frac{2n}{3}.  $ 

Since the above analysis holds for arbitrary $s_t\in [n]$, the lemma is proved.  
$\hfill\square $

The following lemma follows by similar arguments with Lemma~\ref{lem:expecW1}.  
\begin{lemma}\label{lem:expecW12}
For any $n \geq 2$, the output matrix $\Aol$ of Algorithm~\ref{Alg:CuMN} satisfies 
\begin{equation*}
\mathbb E[\Aol] \preceq  \frac{1}{3} \I + \frac{1}{3} \Big({\A}^{(v_t)}+[{\A}^{(v_t)}]^T)\Big).
\end{equation*}		
\end{lemma}

The following consensus rate for  Alg.~\ref{Alg:CuMNagentpersp1EuclideanAlgorithm} is proved similarly to Theorem~\ref{the:ou_equirand}.  
\begin{theorem}
\label{the:ou_equirand2}	
Let $\widetilde{\W}$ be a U-EquiStatic matrix  with consensus rate $\rho $, and $\widetilde{\W}^{(t)}$ be an OU-EquiDyn matrix generated by Alg.~\ref{Alg:CuMNagentpersp1EuclideanAlgorithm}, it holds that
\begin{align*}\E\br{\big\| \Pin \widetilde{\W}^{(t)} \x\big\|^2} \leq \big(1-\frac{4}{3}\eta(1-\eta)(1-\rho)\big)\norm{\Pin \x}^2,
~~\forall \x \in \R^n.\end{align*}
\end{theorem}

\section{Applying EquiTopo Matrices to Decentralized Learning}

\subsection{Convergence of DSGD for strongly convex cost functions}
\label{dsgd:stron}

We assume that $f_i(\x)$ is $\mu$-strongly convex for any $i$, i.e., there exists a constant $\mu>0$ such that
\begin{align*}
f_i(\y) \ge f_i(\x) + \langle \nabla f_i(\x), \y - \x \rangle +\frac {\mu}{2} \norm{\y-\x}^2, ~~\forall \x, \y \in \R^d.
\end{align*}

As we have tested the performance of DSGD with EquiTopo matrices for strongly convex cost functions,     
we attach the following convergence result of the algorithm (\ref{dsgd}).
The proof follows by \cite[Theorem $2$]{koloskova2020unified} (or Appendix A.4 therein) and is omitted here.

\begin{theorem}
	Consider the DSGD algorithm \eqref{dsgd} utilizing the EquiTopo matrices, and $f_i$ being $\mu$-strongly convex for all $i$. Under Assumptions A.1-A.3, it holds that
	\begin{equation*}
	\begin{aligned}
	\frac {1}{H_T} \sum_{t=0}^{T}h^{(t)}\E \Big[
	f(\bar \x^{(t)}) - f(\x^*)
	\Big]=
	\tilde{\Om} \Big( \frac {\sigma^2}{nT}
	+\frac {\kappa \beta \sigma^2}
	{(1-\beta)T^2}
	+\frac {\kappa \beta b^2}
	{(1-\beta)^2 T^2}
	+\frac{1}{1-\beta}\exp\big(-\frac{(1-\beta)T}{\kappa}\big)
	\Big)
	\end{aligned}
	\end{equation*}
	where $\kappa = L/\mu$,  $\tilde{\Om}(\cdot)$ hides constants and polylogarithmic factors, positive weights $h^{(t)} = \pr{1 - \frac{\mu \gam  }{2}}^t$, and $H_T=\sum_{t=0}^T h^{(t)}$.
	Furthermore,
	\begin{itemize}[leftmargin = 1.5em]
		\item  $\beta = \rho$ with D-EquiStatic $\W$ or U-EquiStatic $\widetilde{\W}$; 
		\item  $\beta = \sqrt{(1+\rho)/2}$ for OD-EquiDyn $\W^{(t)}$ (Alg.~\ref{Alg:OD-EquiRand} with $\eta = 1/2$), and $\beta = \sqrt{(2+\rho)/3}$ for OU-EquiDyn $\widetilde{\W}^{(t)}$ (Alg.~\ref{Alg:CuMN} with $\eta = 1/2$). 
	\end{itemize}
\end{theorem}	

\subsection{Transient iteration}
\label{Ap:TranIter}

\textbf{The computation of transient iteration.}

For nonconvex cost functions, the convergence rate of (\ref{dsgd}) is given by
\begin{equation*}
\begin{aligned}
\frac 1{T} \sum_{t=0}^{T-1} \E \Big[ \norm{\nabla f(\bar{\x}^{(t)})}^2\Big]=
\Om \Big( \frac {\sigma}{\sqrt{nT}}
+\frac {\beta^{\frac 23} \sigma^{\frac 23}}{T^{\frac 23}(1-\beta)^{\frac 13}}
+\frac {\beta^{\frac 23} b^{\frac 23}}{T^{\frac 23}(1-\beta)^{\frac 23}}
+\frac {\beta}{T(1- \beta)}
\Big)
\end{aligned}
\end{equation*}	
To reach the linear speedup stage, 
the iteration T has to be sufficiently large so that the $\sqrt{nT}$-term
dominates, i.e.,
$
\frac {\sigma}{\sqrt{nT}} \ge
\frac {\beta^{\frac 23} \sigma^{\frac 23}}{T^{\frac 23}(1-\beta)^{\frac 13}},
$
$
\frac {\sigma}{\sqrt{nT}} \ge
\frac {\beta^{\frac 23} b^{\frac 23}}{T^{\frac 23}(1-\beta)^{\frac 23}},
$
and moreover,
$
\frac {\sigma}{\sqrt{nT}} \ge 
\frac {\beta}{T(1- \beta)}.
$  
Then $T\ge \frac {\beta^4 n^3}{(1-\beta)^2 \sigma^2}$,
$T\ge \frac {\beta^4 b^4 n^3}{(1-\beta)^4 \sigma^6}$,
and $T \ge \frac{\beta^2 n}{(1-\beta)^2 \sigma^2}$. 
Substituting $\beta$ into the inequalities, transient iterations under different networks can be computed.
Similar methods can be adopted for the transient iterations of the distributed gradient tracking algorithm.

Under different network topologies, for non-convex and strongly convex cost functions, 
convergence results and transient iterations are shown in Table \ref{Table:dsgd} and \ref{Table:dsgdstro}.
The results indicate that the proposed networks are at faster rates.

\begin{table}[h]
	\caption{\small For non-convex cost functions, per-iteration communication and convergence rate comparison between DSGD over different topologies.
	The smaller the transient iteration complexity is, the faster the algorithm converges.}
	\begin{tabular}{rccllc}
		\toprule
		&\textbf{Topology}        &\textbf{Per-iter Comm.} 			& \textbf{Convergence Rate} 			& \textbf{Trans. Iters.} \\ \midrule
		&Ring 
		&$\Theta(1)$                
		&$\Om \Big( \frac {\sigma}{\sqrt{nT}}
		+\frac {n^{\frac 23}\sigma^{\frac 23}}{T^{\frac 23}}
		+\frac {n^{\frac 43}b^{\frac 23}}{T^{\frac 23}}
		+\frac {n^2}{T}\Big)$   
		&$\Om(n^{11})$          		\vspace{0.5mm} \\ 
		&Torus   
		&$\Theta(1)$                        
		&$\Om \Big( \frac {\sigma}{\sqrt{nT}}
		+\frac {n^{\frac 13} \sigma^{\frac 23}}{T^{\frac 23}}
		+\frac {n^{\frac 23} b^{\frac 23}}{T^{\frac 23}}
		+\frac {n}{T}\Big)$      
		&$\Om(n^7)$         	   \vspace{0.5mm}\\
		&Static Exp.    &$\Theta(\ln(n))$      
		&$\Om \Big( \frac {\sigma}{\sqrt{nT}}
		+\frac {\ln{^{\frac 13}}(n) \sigma^{\frac 23}}{T^{\frac 23}}
		+\frac {\ln{^{\frac 23}}(n)b^{\frac 23}}{T^{\frac 23}}
		+\frac {\ln(n)}{T}\Big)$   
		&$\Om(n^3 \ln^4(n))$        \vspace{0.5mm}\\
		&O.-P. Exp.	
		& $1$      & $\Om \Big( \frac {\sigma}{\sqrt{nT}}
		+\frac {\ln^{{\frac 13}}(n) \sigma^{\frac 23}}{T^{\frac 23}} 
		+\frac {\ln{^{\frac 23}}(n)b^{\frac 23}}{T^{\frac 23}}
		+\frac {\ln(n)}{T}\Big)$   
		&$\Om(n^3 \ln^4(n))$  \vspace{0.5mm}  \\
		&{\color{blue}D(U)-EquiStatic} 				&{\color{blue}$\Theta(\ln(n))$}  
		& {\color{blue}$\Om \Big( \frac {\sigma}{\sqrt{nT}}
			+\frac {\sigma^{\frac 23}} {T^{\frac 23}}
			+\frac {b^{\frac 23}}
			{T^{\frac 23}}
			+\frac {1}{T}\Big)$} 
		&{\color{blue}$\Om(n^3)$}			\vspace{0.5mm}\\
		&{\color{blue}OD (OU)-EquiDyn}
		&{\color{blue}$1$ }
		&{\color{blue}$\Om \Big( \frac {\sigma}{\sqrt{nT}}
			+\frac{\sigma^{\frac 23}} {T^{\frac 23}}
			+\frac{b^{\frac 23}} {T^{\frac 23}}
			+\frac {1}{T}\Big)$} 
		& {\color{blue}$\Om(n^3)$}      \\
		\bottomrule
	\end{tabular}
	\label{Table:dsgd}
\end{table}

\begin{table}[h]
	\caption{\small For strongly convex cost functions, per-iteration communication and convergence rate comparison between DSGD over different topologies. 
	The smaller the transient iteration complexity is, the faster the algorithm converges.}
	\begin{tabular}{rccllc}
		\toprule
		&\textbf{Topology}        &\textbf{Per-iter Comm.} 			& \textbf{Convergence Rate} 			& \textbf{Trans. Iters.} \\ \midrule
		&Ring           			&$\Theta(1)$                
		&$\tilde{\Om} \Big(\frac {\sigma^2}{nT}+\frac {  \kappa n^2\sigma^2}{T^2}
		+\frac {\kappa  n^4 b^2}{T^2}\Big)$  
		&$\tilde\Om(\kappa n^5)$          		\vspace{0.5mm} \\ 
		&Torus                   &$\Theta(1)$                        
		&$\tilde{\Om} \Big(\frac {\sigma^2}{nT}+\frac { \kappa  n\sigma^2}{T^2}
		+\frac {\kappa  n^2 b^2}{T^2}\Big)$    
		&$\tilde\Om(\kappa n^3)$         	   \vspace{0.5mm}\\
		&Static Exp.       &$\Theta(\ln(n))$      
		&$\tilde{\Om} \Big(\frac { \sigma^2}{nT}+\frac {\kappa \ln(n) \sigma^2}{T^2}
		+\frac {\kappa \ln^2(n)b^2}{T^2}
		\Big)$
		&$\tilde\Om(\kappa n\ln^2(n))$
		\vspace{0.5mm}\\
		&O.-P. Exp.			
		& $1$      
		&$\tilde{\Om} \Big(\frac {\sigma^2}{nT}+\frac { \kappa  \ln(n)\sigma^2}{T^2}
		+\frac {\kappa \ln^2(n)b^2}{T^2}
	    \Big)$  
		&$\tilde\Om(\kappa n\ln^2(n))$  \vspace{0.5mm}  \\
		&{\color{blue}D(U)-EquiStatic} 				&{\color{blue}$\Theta(\ln(n))$}  
		& {\color{blue}$\tilde{\Om} \Big( \frac {\sigma^2}{nT}+\frac { \kappa  \sigma^2}{T^2}
		+\frac{\kappa b^2}{T^2}\Big)$} 
		&{\color{blue}$\tilde\Om(\kappa n)$}			\vspace{0.5mm}\\
		&{\color{blue}OD (OU)-EquiDyn} 		&{\color{blue}$1$ }
		&{\color{blue}$\tilde{\Om} \Big( \frac {\sigma^2}{nT}+\frac { \kappa  \sigma^2}{T^2}
		+\frac{\kappa  b^2}{T^2}\Big)$} 
		& {\color{blue}$\tilde\Om(\kappa n)$}      \\
		\bottomrule
	\end{tabular}
	\label{Table:dsgdstro}
\end{table}

\subsection{Decentralized stochastic gradient tracking algorithm}
\label{Ap:DSTGA}

We write local variables compactly into matrix form, for instance
\eq{
  \Xl{t} = \br{\xl{t}_1, \cdots, \xl{t}_n}\tp \in \Real^{n\times d}, \ 
  \gF{t} = \br{\nabla f_1\pr{\xl{t}_1}, \cdots, \nabla f_n\pr{\xl{t}}}\tp \in \Real^{n\times d}.    
}
The matrices $\Yl{t}, \Gl{t}\in \Real^{n\times d}$ are defined analogously. We also denote $\gF{-1} = \Gl{-1} = \zero$ for notational simplicity.

Clearly, 
the DSGT algorithm can be simplified as
\eql{\label{eq:XYstackrel1}}{
    \mx{\Xl{t+1} \\ \Yl{t+1} }
    =
    \mx{
        \Wl{t} & -\gam \Wl{t} \\
        \zero & \Wl{t}
    }
    \mx{
        \Xl{t} \\
        \Yl{t}
    }
    +
    \mx{
        \zero \\
        \Gl{t+1} - \Gl{t}
    }.
}

For simplicity, we define
\eq{
    \Wl{j:k} = \Wl{k} \cdots \Wl{j},
    ~~\forall
    k \ge j \ge 0,
}
and moreover,
$\Wl{j:k} = \II$ for $j > k$.

Notice that
\eq{
    \mx{
     \Wl{k} & -\gam \Wl{k}  \\
     \zero & \Wl{k}
    }
    \cdots
        \mx{
     \Wl{j} & -\gam \Wl{j}  \\
     \zero & \Wl{j}
    }
    =
    \mx{
        \Wl{j:k} & - \gam \pr{k-j+1} \Wl{j:k}  \\
        \zero & \Wl{j:k}
    }.
}
Consequently,
it holds for all 
$t \ge 1$ that
\eql{\label{eq:Xexpand1}}{
    \Xl{t} &= \Wl{0:\pr{t-1}}\Xl{0} - \gam \sum_{j=0}^{t - 1} \pr{t-j}\Wl{j:\pr{t-1}} \pr{\Gl{j} - \Gl{j - 1}}.  
}

Moreover, we have two inequalities as follows.

\begin{lemma}
\label{lem:t^2exp1}
    For any $t\geq 1$ and $\beta \in (0,1)$,
    we have
    \eq{
        t^2\bet^{2\pr{t-1}} \leq \frac{\pa}{\pr{1 - \bet}^2  }\bet^{t - 1  }  , 
        \quad t^2\bet^{t - 1  } \leq \frac{\pb}{\pr{1 - \bet}^2} \pr{\frac{1 + \bet}{2}}^{t - 1},     
    }
    where $\pa = 4 $, $\pb = 16 $.  
\end{lemma}
\emph{Proof.}
  Define $r(x) = x^2\bet^{x - 1 }$, 
  where $x\geq 1$. 
  Then, for the first inequality, it suffices to show that $r(x) \leq \pa/{(1-\beta)^2}$ for $x\geq 1$.  
  Due to $r'(x) = x\bet^{x-1 }\pr{2 + x\ln \bet  } $,  
   $r(x)$ attains its maximum at $x_0 = \max\dr{1,  - \frac{2}{\ln \bet }}   $.

  If $ - \frac{2}{\ln \bet } > 1 $,  
  by combining with the fact that $  \ln \bet \leq \bet - 1 < 0 $ for $\bet\in (0, 1)$, we have      
  \eq{
    x^2\bet^{x-1} \leq r\pr{- \frac{2}{\ln \bet  } } = \frac{4  }{\pr{\ln\bet}^2}\bet^{- \frac{2}{\ln\bet} - 1} 
    \leq \frac{4  }{\pr{\ln\bet}^2} \leq \frac{4}{\pr{1 - \bet}^2}.    
  } 
  
  If $- \frac{2}{\ln \bet } \leq 1 $, then
  \eq{
    x^2\bet^{x-1} = r\pr{1} = 1 \leq \frac{4 }{\pr{1 - \bet}^2}.   
  } 
  
  The second inequality follows by similar arguments and the fact that $\sqrt{\bet} \leq \frac{1 + \bet}{2}$.  
The proof is completed.
$\hfill\square$

The following lemma is a generalization of Cauchy-Schwartz inequality.
Its proof follows by using $\nf{\AA + \BB}^2 \leq \frac{1}{\alpha}\nf{\AA}^2 + \frac{1}{1 - \alpha}\nf{\BB}^2 \ (\alpha \in \pr{0, 1})  $ repeatedly.

\begin{lemma}\label{lem:mx1}

    Consider a sequence of matrices $\dr{\BB_i}_{i=1}^m$. If $\alpha_1, \alpha_2, \cdots, \alpha_m > 0$ and $\sum_{i=1}^m \alpha_i \leq 1$, then 
    \eq{
        \nf{\sum_{i=1}^{m}\BB_i}^2 \leq \sum_{i=1}^{m} \frac{1}{\alpha_i} \nf{\BB_i }^2.
    }

\end{lemma}

We define a potential function as
\eql{\label{eq:Lyz1}}{
   \Ly{t} 
   =&  \frac{4 \bet^t }{\pr{1 - \bet}^2 }\Et{\Pin\Xl{0}} + \frac{16  \pa n \gam^2  }{\pr{1 - \bet}^4} \sigma^2  \\  
     & + \frac{4 \pb \gam^2 }{\pr{1 - \bet}^4}\sum_{j=0}^{t - 1} \pr{\frac{1 + \bet}{2}}^{t - j - 1 }   \Et{ \gF{j} - \gF{j-1}},~~
     \forall t\ge 1,
}
and moreover,
\eql{\label{eq:Ly0!1lzy}}{
    \Ly{0} = \frac{4  }{\pr{1 - \bet}^2} \nf{\Con\Xl{0}}^2 + \frac{16 \pa n \gam^2 }{\pr{1 - \bet}^4} \sigma^2.
}

The following Lemma~\ref{lem:ConXlt1Lyz!1} and Lemma~\ref{lem:gFt+1-gFt} are used to prove Lemma~\ref{lem:sumEtConXlt1}. 
Theorem~\ref{the:dsgt} follows by combining Lemma~\ref{lem:sumEtConXlt1} with the descent lemma (Lemma~\ref{lem:choosegam1}).  
\begin{lemma}\label{lem:ConXlt1Lyz!1}
Consider the DSGT~\eqref{dsgta}.
Let Assumptions A.1 and A.2 hold.
If $\dr{\Wl{t}}_{t\geq 0}$ have convergence rate $\beta$, i.e., $\E\br{\nt{\Con\Wl{t}\yy}^2} \leq \beta^2\nt{\Con\yy}^2$ for any $\yy\in \Real^n$,  
then
    \eq{
        \E\br{\nf{\Con\Xl{t}}^2} \leq \frac{1 - \bet}{2}\Ly{t},~~
        \forall t\ge 0.
    }
\end{lemma}
\emph{Proof.}
  The case $t = 0$ follows by definition directly.
  
  For $t\geq 1$, we
  define \eq{\Ql{t, 1} = \Wl{0:\pr{t-1}}\Xl{0} - \gam \sum_{j=0}^{t - 1} \pr{t-j}\Wl{j:\pr{t-1}} \pr{\gF{j} - \gF{j-1}},} and
  \eq{\Ql{t, 2} = - \gam \sum_{j=0}^{t - 1} \pr{t-j}\Wl{j:\pr{t-1}} \pr{\Gl{j} - \gF{j} - \Gl{j-1}  + \gF{j-1} }.  }

  Recalling \eqref{eq:Xexpand1} gives
  \eq{
    \Xl{t} = \Ql{t, 1} + \Ql{t, 2}.    
  }

  Then
  \eq{
    \E\br{\nf{\Con\Xl{t}}^2} \leq 2\Et{\Con \Ql{t, 1}  } + 2\Et{\Con \Ql{t, 2} }.
  }

  Rearranging $\Ql{t, 2}$ yields
  \eq{
    \Ql{t, 2} =& - \gam \Wl{t-1}\pr{\Gl{t-1} - \gF{t-1}} \\ 
      & - \gam \sum_{j=0}^{t-2} \pr{\pr{t - j}\Wl{j:(t-1)} - \pr{t - j - 1}\Wl{\pr{j+1}:\pr{t-1}}}\pr{\Gl{j} - \gF{j}}.
  }

  By Assumption~A.2 and  the assumption on consensus rate, we have  
  \eq{
     &\Et{\Con\Ql{t, 2}} \\
      =& \gam^2 \E\br{\nf{\Con\Wl{t-1}\pr{\Gl{t-1} - \gF{t-1}}}^2}
       \\
      & + \gam^2 \sum_{j=0}^{t-2} \E\br{\nf{\pr{\pr{t - j}\Con\Wl{j:(t-1)} - \pr{t - j - 1}\Con\Wl{\pr{j+1}:\pr{t-1}}}\pr{\Gl{j} - \gF{j}}}^2} \\
      \leq& \gam^2 \beta^2 \Et{\Gl{t-1} - \gF{t-1}}   
       + \gam^2\sum_{j=0}^{t - 2} 4\pr{t - j}^2 \bet^{2\pr{t - j - 1}}\Et{\Gl{j} - \gF{j}}  \\
      \leq&  4n  \gam^2 \sum_{j=1}^{t} j^2\bet^{2\pr{j-1}} \sigma^2.  
  }

  By Lemma~\ref{lem:t^2exp1},
  we have
  \eq{
    n \gam^2 \sum_{j=1}^{t} j^2 \bet^{2\pr{j - 1} } \sigma^2 \leq \frac{n  \pa \gam^2 }{\pr{1 - \bet}^2}\sum_{j=1}^{t} \bet^{j-1} \sigma^2 \leq \frac{n \pa \gam^2  }{\pr{1 - \bet}^3} \sigma^2.  
  }

 As a result,
 \begin{align*}
 \Et{\Con\Ql{t, 2}} \le 
 \frac{4n \pa \gam^2  }{\pr{1 - \bet}^3} \sigma^2.
 \end{align*}
 
 Moreover,
  \eq{
    &\E\br{\nf{\Con\Ql{t, 1}}^2  } \\
    \leq& \frac{1}{\pr{1 - \bet}\bet^t}\Et{\Con\Wl{0:\pr{t-1}}\Xl{0}} \\
     & +  \sum_{j=0}^{t - 1} \frac{ \gam^2\pr{t-j}^2 }{\pr{1 - \bet} \bet^{t - j - 1 } }\Et{\Con\Wl{j:\pr{t-1}} \pr{\gF{j} - \gF{j-1}}}  \\
    \leq& \frac{\bet^t }{\pr{1 - \bet}  }\Et{\Pin \Xl{0}} +  \frac{\gam^2}{1 - \bet }\sum_{j=0}^{t - 1} \pr{t-j}^2 \bet^{t - j - 1 } \Et{ \gF{j} - \gF{j-1}}       \\
    \leq& \frac{\bet^t }{1 - \bet }\Et{\Pin \Xl{0}} + \frac{\pb \gam^2 }{\pr{1 - \bet}^3}\sum_{j=0}^{t - 1} \pr{\frac{1 + \bet}{2}}^{t - j - 1 }   \Et{ \gF{j} - \gF{j-1}},
  }
  where the first inequality follows by Lemma~\ref{lem:mx1} and the fact that $\pr{1 - \bet}\sum_{j=0}^{t }\bet^j < 1 $,
  the second inequality is by the assumption on consensus rate, 
  and the third inequality is by Lemma~\ref{lem:t^2exp1}.

 Therefore, the conclusion holds by the definition of $\Phi^{(t)}$.
 $\hfill\square$

\begin{lemma}\label{lem:gFt+1-gFt}
 Consider the DSGT~\eqref{dsgta}.  
Let Assumptions A.1 and A.2 hold.
If $\gam \leq \frac{1}{L}$, it holds for $t\geq 0$ that
    \eq{
        &\Et{\gF{t+1} - \gF{t}} \\
        \leq& 6n \gam^2 L^2 \E\br{\nt{\nabla f\pr{\xa{t}}}^2} + 9L^2 \Et{\Con\Xl{t}}  + 3L^2\Et{\Con\Xl{t+1}} + 3  \gam^2 L^2 \sigma^2.
    }

\end{lemma}
\emph{Proof.}
Clearly,  
\eq{
&\Et{\gF{t+1} - \gF{t}} \le 3\Et{\gFa{t+1} - \gFa{t}} \\
&+3\Et{\gF{t+1} - \gFa{t+1}} + 3\Et{\gF{t} - \gFa{t}}.
}

It follows from Assumption~A.1 that
\eql{\label{inequal1combine1}}{
&\Et{\gF{t+1} - \gF{t}} \\
\leq& 3L^2\pr{\Et{\Xa{t+1}-\Xa{t}}+\Et{\Con\Xl{t+1}} + \Et{\Con\Xl{t}}}.
}

Notice that by induction,  $\sum_{i=1}^{n} \yl{t}=\sum_{i=1}^{n} \gl{t}$. Recalling (\ref{dsgta}) gives
  \eql{\label{eq:xat+1-xat11}}{
    &\xa{t+1} - \xa{t} = \frac{\gam }{n}\sum_{i=1}^{n} \yl{t} =  \frac{\gam  }{n}\sum_{i=1}^{n} \gl{t}  \\
    =&   \frac{\gam }{n}\sum_{i=1}^{n} \Big[\nabla f_i\big(\xa{t}\big) + \pr{\nabla f_i\big(\xl{t}_i\big) - \nabla f_i\big(\xa{t}\big)} + \pr{\gl{t}_i - \nabla f_i\big(\xl{t}_i\big) }  \Big]\\
    =& \gam \nabla f\big(\xa{t}\big) +  \frac{\gam }{n}\sum_{i=1}^{n} \Big[\pr{\nabla f_i\big(\xl{t}_i\big) - \nabla f_i\big(\xa{t}\big)} + \pr{\gl{t}_i - \nabla f_i\big(\xl{t}_i\big) }\Big].  
  }

  By Assumptions~A.1 and A.2, we derive
  \eql{\label{eq:Entxat+1-xat1}}{
    &\E\br{\nt{\xa{t+1} - \xa{t}}^2 } \\
     =& \gam^2\E\br{\nt{\nabla f\big(\xa{t}\big) +  \frac{1}{n}\sum_{i=1}^{n} \pr{\nabla f_i\big(\xl{t}_i\big) - \nabla f_i\big(\xa{t}\big)}}^2 }+\frac{\gam^2\sigma^2}{n}\\
     \leq& 2\gam^2 \E\br{\nt{\nabla f\big(\xa{t}\big)}^2} + \frac{2\gam^2 L^2}{n}\Et{\Con\Xl{t}} + 
     \frac{\gam^2\sigma^2}{n}.
  }
  Due to $\gam \leq \frac{1}{L}$, we have
  \eql{\label{eq:xat+1diff1}}{
    \E\br{\nt{\xa{t+1} - \xa{t}}^2 } \leq  2\gam^2\E\br{\nt{\nabla f\big(\xa{t}\big) }^2}+\frac{2}{n}\Et{\Con\Xl{t}}+\frac{\gam^2\sigma^2}{n}.
  }

  Then the conclusion follows~\eqref{inequal1combine1},~\eqref{eq:xat+1diff1} and  $\nf{\Xa{t+1} - \Xa{t}}^2 = n \nt{\xa{t+1} - \xa{t}}^2$.  
$\hfill\square$

\begin{lemma}\label{lem:sumEtConXlt1}
Consider the DSGT~\eqref{dsgta}.
Let Assumptions A.1 and A.2 hold.
Suppose that $\dr{\Wl{t}}_{t \geq 0}$ have consensus rate $\beta$ and $\y \in \R^n$.
If $ \frac{48 \pb \gam^2 L^2}{\pr{1 - \bet}^4 } \leq \frac{1}{2} $, then
\eq{
    \sum_{t=0 }^{T}\Et{\Con\Xl{t}} \leq& 2 \Ly{0}  +  \frac{48 \pb n \gam^4 L^2}{\pr{1 - \bet}^4 } \sum_{t=1}^{T} \E\br{\nt{\nabla f\big(\xa{t - 1}\big)}^2} + \frac{8\pb \gam^2}{\pr{1 - \bet}^4}\nf{\gF{0}}^2 \\
     & + \frac{24\pb \gam^4 L^2 }{\pr{1 - \bet}^4} \pr{T+1} \sigma^2 + \frac{16 \pa n \gam^2 }{\pr{1 - \bet }^3}   \pr{T+1 } \sigma^2.
}

\end{lemma}
\emph{Proof.}
  By the definition of $\Ly{t}$ in~\eqref{eq:Lyz1}, we have that for $t\geq 0 $,
  \eql{\label{eq:Lyt1}}{
    \Ly{t+1} \leq \pr{\frac{1 + \bet }{2}  } \Ly{t} + \frac{4\pb \gam^2}{\pr{1 - \bet}^4}\Et{\gF{t} - \gF{t-1}} + \frac{8 \pa n \gam^2  }{\pr{1 - \bet}^3}\sigma^2.
  }
  Then, for $t\geq 0$, by Lemma~\ref{lem:ConXlt1Lyz!1} and~\eqref{eq:Lyt1},  we have
  \eq{
    &\Et{\Con\Xl{t}}
    \leq {\frac{1 - \bet }{2}}\Ly{t} \\
    \leq& \Ly{t} - \Ly{t+1} + \frac{4\pb \gam^2}{\pr{1 - \bet}^4}\Et{\gF{t} - \gF{t-1}} + \frac{8 \pa n \gam^2  }{\pr{1 - \bet}^3}\sigma^2.  
  }
  For $t \geq 1$, by Lemma~\ref{lem:gFt+1-gFt},
  we derive
  \eql{\label{eq:ConXltLyz1}}{
    &\Et{\Con\Xl{t}}   
    \leq \Ly{t} - \Ly{t+1}  +  \frac{24\pb n \gam^4 L^2}{\pr{1 - \bet}^4 }  \E\br{\nt{\nabla f\big(\xa{t - 1}
    \big)}^2} \\
     & + \frac{36\pb \gam^2 L^2}{\pr{1 - \bet}^4 } \Et{\Con\Xl{t - 1}}  + \frac{12\pb  \gam^2L^2}{\pr{1 - \bet}^4} \Et{\Con\Xl{t}} \\
      & + \frac{12\pb \gam^4 L^2 }{\pr{1 - \bet}^4} \sigma^2 + \frac{8 \pa n \gam^2 }{\pr{1 - \bet }^3} \sigma^2.  
  }

  It follows by the definition of $\Ly{t}$ that
  \eq{
    \Ly{1} \leq \pr{\frac{1 + \bet }{2}}\Ly{0} + \frac{4\pb \gam^2}{\pr{1 - \bet}^4}\nf{\gF{0}}^2 + \frac{8 \pa n  \gam^2  }{\pr{1 - \bet}^3}\sigma^2.
  }
  By Lemma~\ref{lem:ConXlt1Lyz!1}, we obtain
  \eql{\label{eq:EtXl0Lyz1!1}}{
    \nf{\Con\Xl{0}}^2 \leq \frac{1 - \bet}{2}\Ly{0}  \leq \Ly{0} - \Ly{1} + \frac{4\pb \gam^2}{\pr{1 - \bet}^4}\nf{\gF{0}}^2 + \frac{8 \pa n  \gam^2  }{\pr{1 - \bet}^3}\sigma^2.  
  }

  Taking sum on both sides of~\eqref{eq:ConXltLyz1} and noting that $\frac{48 \pb \gam^2 L^2}{\pr{1 - \bet}^4 } \leq \frac{1}{2} $, $\Ly{T+1} \geq 0 $, the lemma is proved.
$\hfill\square$

Lemma~\ref{lem:gFL1} is standard in the analysis of gradient tracking methods.
We attach its proof for completeness.
\begin{lemma}\label{lem:gFL1}
Consider the DSGT~\eqref{dsgta}.  
Let Assumptions A.1 and A.2 hold.
If $\gam \leq \frac{1}{4L}$, then
\eq{
\E\br{f\big(\xa{t+1}\big)} \leq \E\br{f\big(\xa{t}\big)}  - \frac{\gam }{4} \E\br{\nt{\nabla f\big(\xa{t}\big)}^2} 
+\frac{\gam L^2}{n} \Eb{\nf{\Con\Xl{t}}^2}  
    +\frac{\gam^2 L}{2n} \sigma^2.
}

\end{lemma}
\emph{Proof.}
  By Assumption~A.1, we have
  \eq{
    \Eb{f\big(\xa{t+1}\big)} \leq \Eb{f\big(\xa{t}\big)} - 
    \Eb{\jr{\nabla f\big(\xa{t}\big), \xa{t+1} - \xa{t}}} + \frac{L}{2}\Eb{\nt{\xa{t+1} - \xa{t}}^2}.
  }
  It follows from (\ref{dsgta}) that
  \eq{
    &\Eb{\jr{\nabla f\pr{\xa{t}}, \xa{t+1} - \xa{t} }} \\
    =& \gam \Eb{\jr{\nabla f\pr{\xa{t}}, \nabla f\big(\xa{t}\big) + \frac{1}{n}\sum_{i=1}^{n}\pr{\nabla f_i\big(\xl{t}_i\big) - \nabla f_i\big(\xa{t}\big) }}}  \\
    \geq& \gam \Eb{\nt{\nabla f\big(\xa{t}\big)}^2} - \frac{\gam}{2} \Eb{\nt{\gf{\xa{t}}}^2} - \frac{\gam  }{2n}\sum_{i=1}^{n} \Eb{\nt{\nabla f_i\big(\xl{t}_i\big) - \nabla f_i\big(\xa{t})}}^2  \\
    \geq& \frac{\gam }{2} \Eb{\nt{\nabla f\big(\xa{t}\big)}^2} - \frac{\gam L^2}{2n} \Eb{\nf{\Con\Xl{t}}}^2,
  }
  where the first equality is by Assumption~A.2 and~\eqref{eq:xat+1-xat11};
  the second inequality is by Assumption~A.1.

  Recalling \eqref{eq:Entxat+1-xat1} yields 
\eq{
    &\Eb{f\big(\xa{t+1}\big)} \\
    \leq& \Eb{f\big(\xa{t}\big)}
    -\frac{\gam}{2}\pr{1 - 2\gam L}
    \Eb{\nt{\nabla
    f\big(\xa{t}\big)}^2}
    + \frac{\gam L^2}{2n}\pr{1+2\gam L} \Eb{\nf{\Con\Xl{t}}^2}
    +\frac{\gam^2\sigma^2L}{2n} \\
    \leq&\Eb{f\big(\xa{t}\big)}
    - \frac{\gam}{4} \Eb{\nt{\nabla
    f\big(\xa{t}\big)}^2}
    +\frac{\gam L^2}{n} \Eb{\nf{\Con\Xl{t}}^2}
    +\frac{\gam^2 \sigma^2 L}{2n}.
}
The lemma is proved.
$\hfill\square$

Referring to Lemma $26$ of \cite{koloskova2021improved}, we have the following result.
\begin{lemma}\label{lem:choosegam1}
Let $A, B, C, T$ and $\alpha$ be positive constants. 
Define
\eq{
g(\gam) = \frac{A}{\gam T} + B\gam + C\gam^2.
}
Then
\eq{
\inf_{\gam \in (0, \alpha]} g\pr{\gam}
\leq 2\pr{\frac{AB}{T}}^{\frac{1}{2}} + 2C^{\frac{1}{3}}\pr{\frac{A}{T}}^{\frac{2}{3}} + \frac{A  }{\alpha  T}.
}
\end{lemma}

\textbf{Proof of Theorem \ref{the:dsgt}}
Define $f^* = \inf_x f(x)$, $\dffFinit = f\pr{\xa{0}} - f^*$,
$\ConXinit = \nf{\Con\Xl{0}}^2$
and $\dFinit = \sum_{i=1}^{n}\nt{\nabla f_i\pr{\xl{0}_i}}^2$.
We show that
    \eql{\label{eq:DSGTdetail1}}{
        &\frac{1}{T + 1  }\sum_{t=0 }^{T} \Et{\nabla f\pr{\xa{t}}} \\
        \leq&
        \Ot{\sqrt{\frac{\dffFinit L\sigma^2  }{n T}} + \frac{1}{1 - \bet}\pr{\frac{\dffFinit L \sigma}{T}}^{\frac{2}{3}} + \frac{\dffFinit}{\pr{1 - \bet }^2 T}  + \frac{L^2 \ConXinit}{\pr{1-\bet}^2 n T} + \frac{\dFinit}{nT}   },  
    }

    Let $\gam\leq\frac{\pr{1-\bet}^2}{50 L}$ to satisfy the conditions in Lemmas~\ref{lem:sumEtConXlt1} and~\ref{lem:gFL1}.
    Then, we have
    \eq{
        &\frac{1}{T + 1  }\sum_{t=0}^{T} \Et{\nabla f\big(\xa{t}\big)}   \\
        \leq& \frac{4}{\gam \pr{T + 1} }\pr{f\big(\xa{0}\big) - f^*} + \frac{4L^2}{n \pr{T + 1}}\sum_{t=0}^{T} \Et{\Con\Xl{0}} + \frac{2\gam L }{n}\sigma^2  \\
        \leq& \frac{4}{\gam T }\pr{f\big(\xa{0}\big) - f^*} + \frac{2\gam L }{n}\sigma^2 
        + \frac{8 L^2}{n T } \Ly{0} 
        + \frac{32\pb\gam^2 L^2}{\pr{1-\bet}^4 n T}\nf{\gF{0}}^2 \\
        &+  \frac{192\pb  \gam^4 L^4}{\pr{1 - \bet}^4  \pr{T+1} } \sum_{t=1}^{T} \E\br{\nt{\nabla f\big(\xa{t - 1} \big)}^2} 
         + \frac{96\pb \gam^4 L^4 }{\pr{1 - \bet}^4 n} \sigma^2 + \frac{64\pa  \gam^2 L^2 }{\pr{1 - \bet}^3 } \sigma^2.   
    }
    If $\gam \leq \frac{{1-\bet}}{  10 L} $, then $\frac{192\pb \gam^4 L^4}{\pr{1-\bet}^4  } \leq \frac{1}{2} $.   Then,
    we have  
    \eql{\label{eq:sumgftbeforerearrange1}}{
        \frac{1}{T + 1}\sum_{t=0}^{T} \Et{\nabla f\big(\xa{t}\big)}  
        \leq  \frac{8}{\gam T }\pr{f\big(\xa{0}\big) - f^*} + \frac{4\gam L }{n}\sigma^2 
        + \frac{192 \pb \gam^4 L^4 }{\pr{1 - \bet}^4 n} \sigma^2& \\
         + \frac{16 L^2}{n T } \Ly{0}
          + \frac{64\pb \gam^2 L^2}{\pr{1 - \bet}^4 n T}\nf{\gF{0}}^2 
          + \frac{128\pa  \gam^2 L^2 }{\pr{1 - \bet}^3}\sigma^2 &.  
    }
    By~\eqref{eq:Ly0!1lzy} and  $\pa = 4 $, $\pb = 16 $ defined in Lemma~\ref{lem:t^2exp1}, if $\gam \leq \frac{\pr{1 - \bet}^2}{50 L}$ and $T \geq \frac{1}{1 - \bet} $,   we have
    \eql{\label{eq:sumEnfgfxa1}}{
        &\frac{1}{T + 1}\sum_{t=0}^{T} \Et{\nabla f\big(\xa{t}\big)}  \\
        \leq& \frac{8}{\gam T} \pr{f\big(\xa{0}\big) - f^*} + \frac{4\gam L }{n}\sigma^2 \\
        & + \frac{1}{T}\pr{\frac{64 L^2}{\pr{1 - \bet}^2 n  } \nf{\Con\Xl{0}}^2 + \frac{  1024  \gam^2 L^2 }{\pr{1 - \bet}^4 n }\nf{\gF{0}}^2 }  \\
        & + \frac{   2048  \gam^2L^2 }{\pr{1 - \bet}^4 T}\sigma^2 + \frac{3072\gam^4L^4}{\pr{1 - \bet}^4n}\sigma^2 + \frac{512\gam^2L^2}{\pr{1 - \bet}^3}\sigma^2  \\ 
        \leq& \frac{8}{\gam T }\pr{f\big(\xa{0}\big) - f^*} + \frac{5\gam L }{n}\sigma^2  + \frac{2560 \gam^2L^2}{\pr{1 - \bet}^3}\sigma^2 \\
        & + \frac{1}{n T}\pr{\frac{64 L^2}{\pr{1 - \bet}^2   } \nf{\Con\Xl{0}}^2 + \nf{\gF{0}}^2 }.  
    }

    To meet the conditions of Lemma~\ref{lem:sumEtConXlt1}, Lemma~\ref{lem:gFL1},~\eqref{eq:sumgftbeforerearrange1},~\eqref{eq:sumEnfgfxa1}, it suffices to let $\gam \leq \frac{\pr{1 - \bet}^2}{  50  L  }$.
    Then,~\eqref{eq:DSGTdetail1} follows by setting $g\pr{\gam}$ to be the RHS of~\eqref{eq:sumEnfgfxa1}, $A = 8 \pr{f\pr{\xa{0}} - f^*} $, $B = \frac{5 L}{n}\sigma^2 $, $C = \frac{2560   L^2}{\pr{1 - \bet}^3}\sigma^2 $ and $\alpha = \frac{\pr{1 - \bet}^2}{50 L} $
        in Lemma~\ref{lem:choosegam1}.
$\hfill\square$

\section{Numerical Experiments}
\label{Ap:Exp}

\subsection{Network-size independent consensus rate}
\begin{figure}[h!]
	\centering
	\includegraphics[scale=0.26]{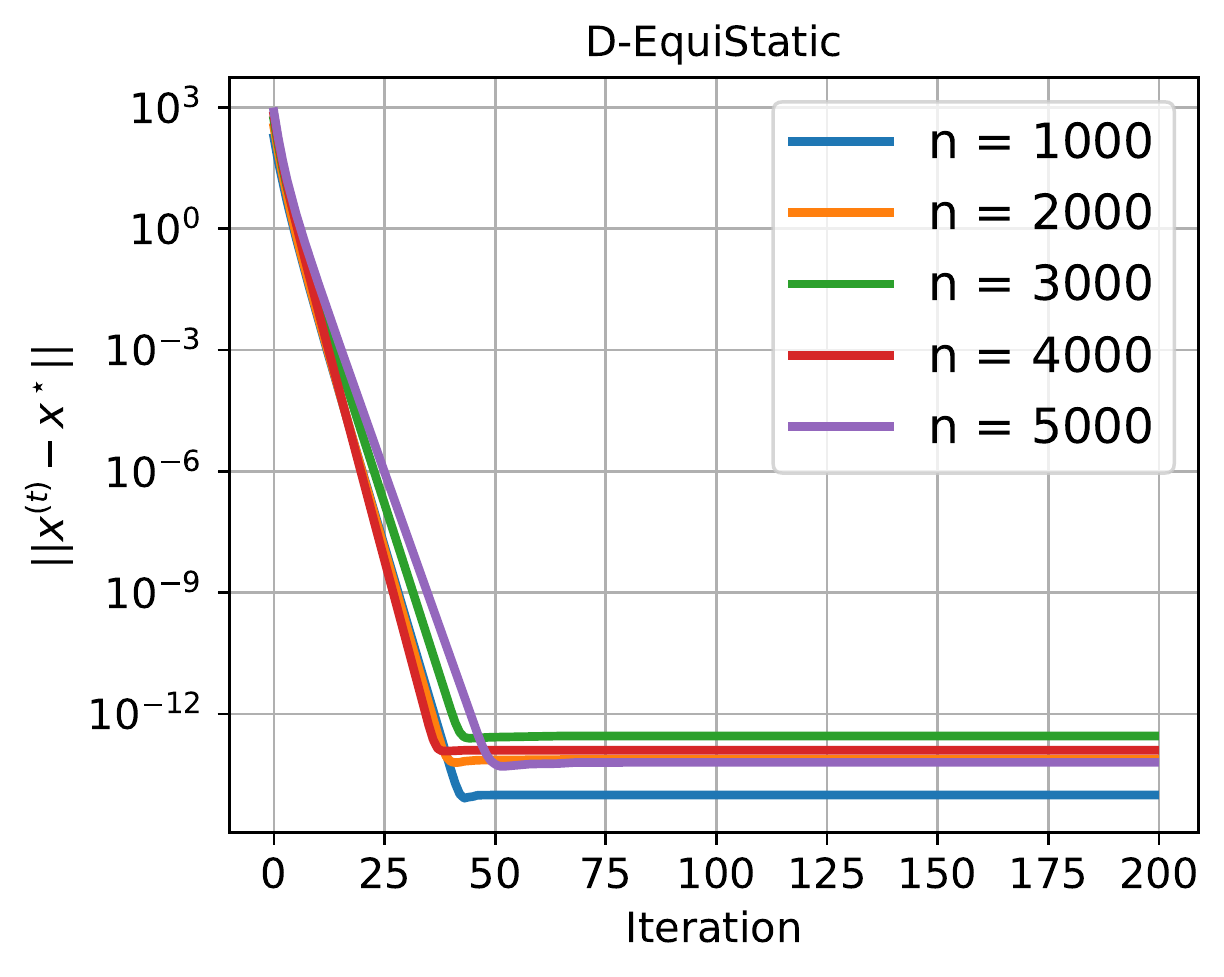}
	\includegraphics[scale=0.26]{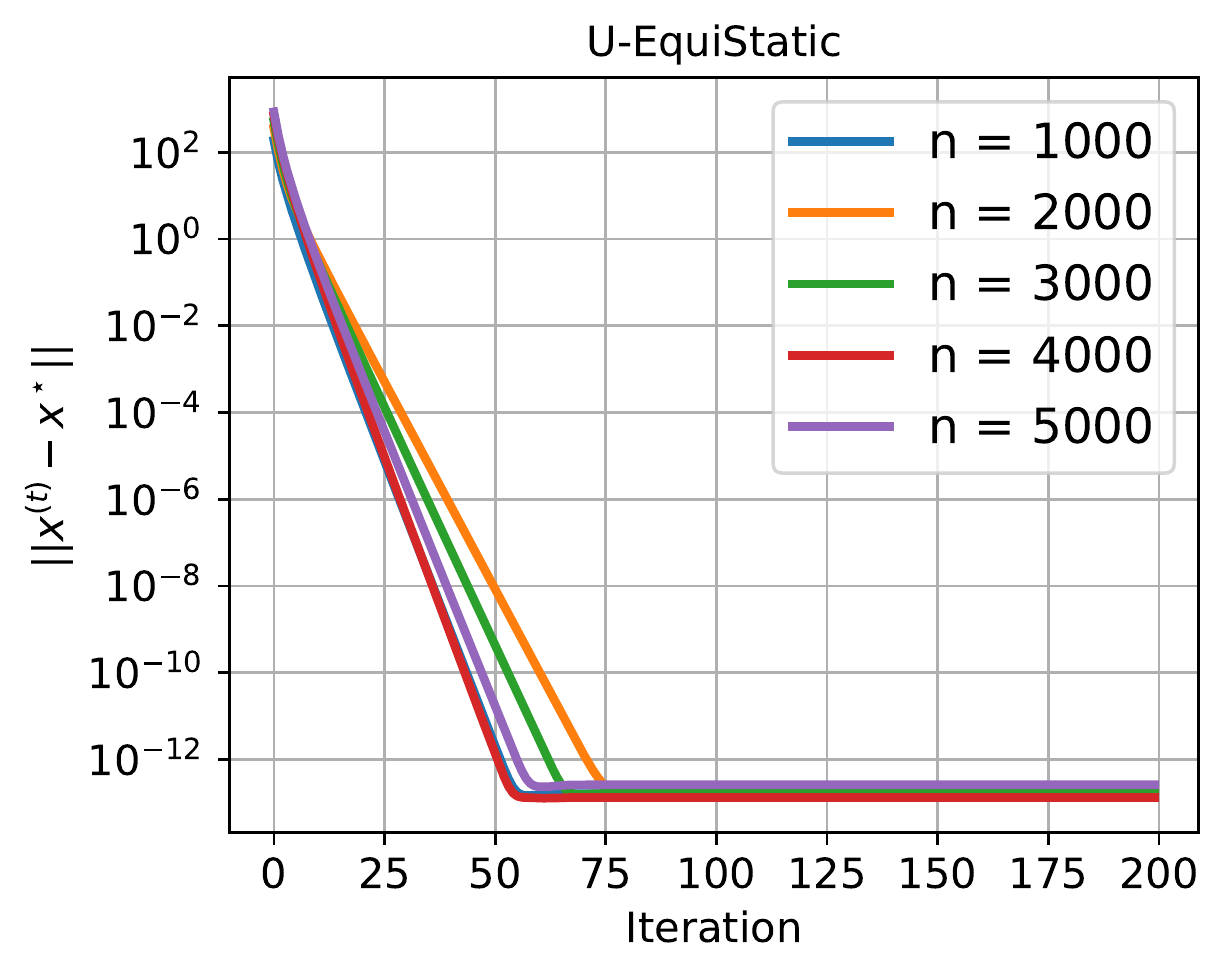}
	\includegraphics[scale=0.26]{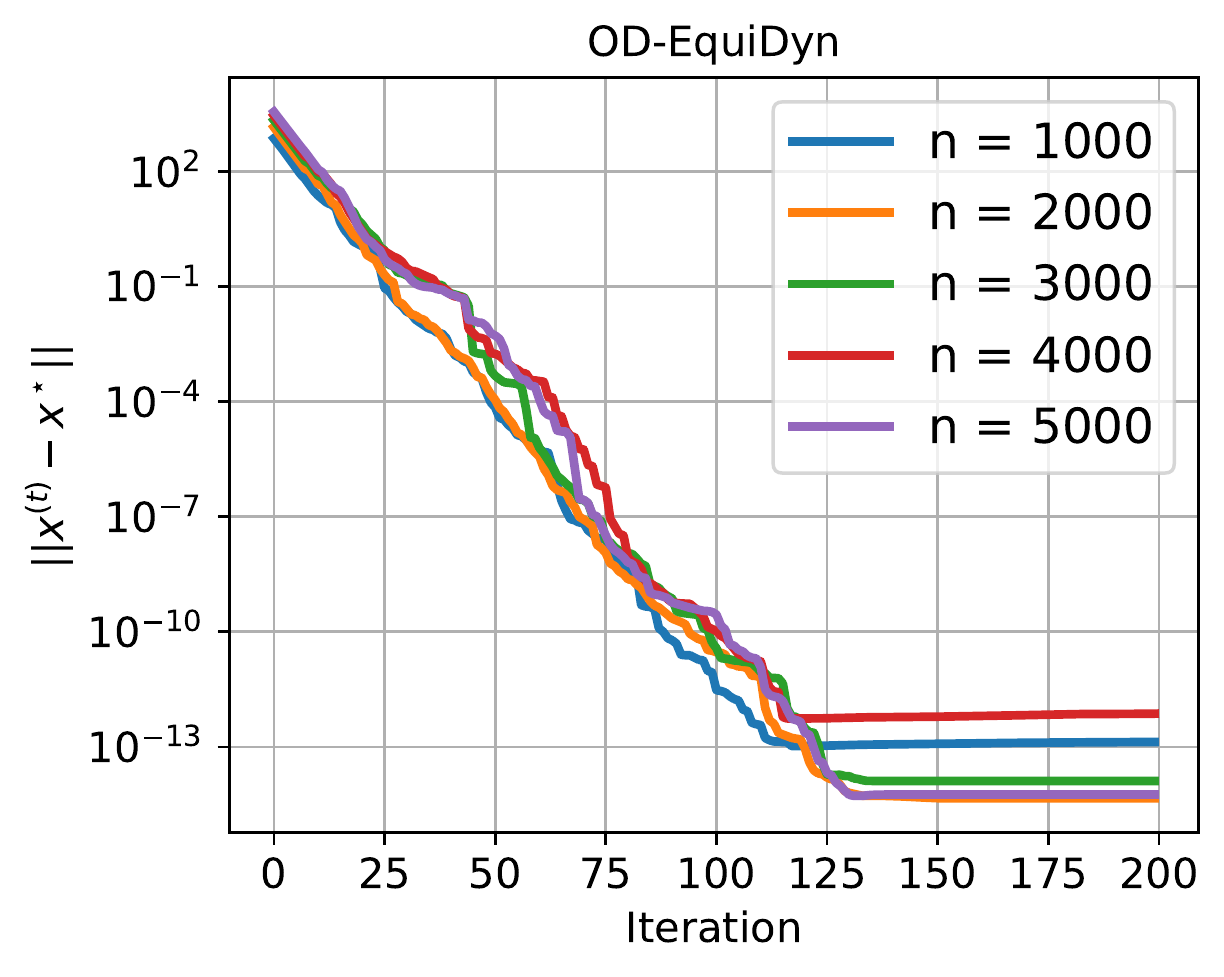}
	\includegraphics[scale=0.26]{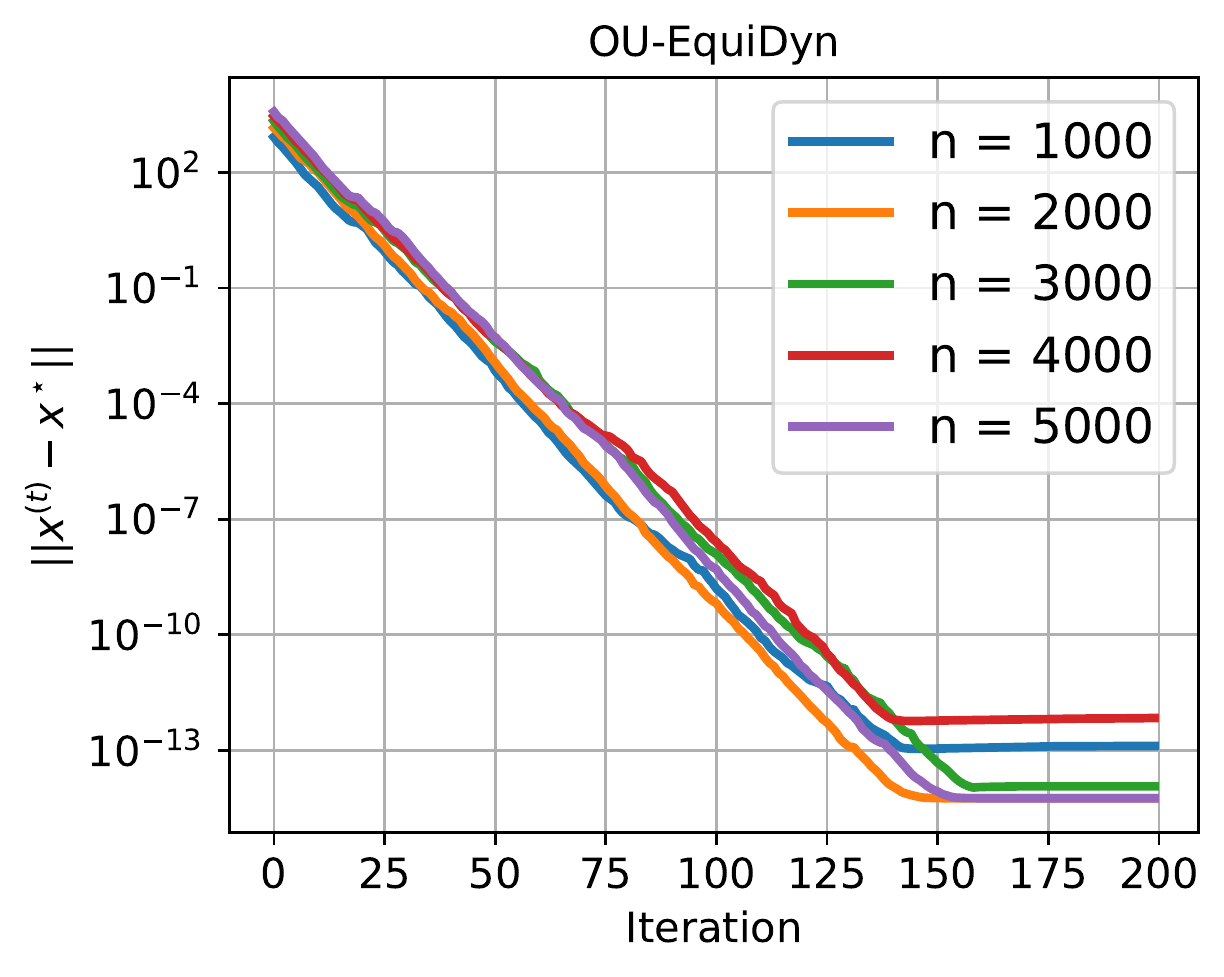}
	\caption{\small The EquiTopo graphs can achieve network-size independent consensus rates.}
	\label{fig:more-n-independent}
\end{figure}

In this experiment, we set $M=5\ln(n)$ for D-EquiStatic and $M=2\ln(n)$ for U-EquiStatic, which is consistent with Theorems \ref{lm-staconstr} and~\ref{thm-u-equirand}. For OD-EquiDyn and OU-EquiDyn, we set $M=5\ln(n)$ and $\eta=0.5$. Fig.~\ref{fig:more-n-independent} shows that the consensus rate is independent of the network size for all EquiTopo graphs. The results are obtained by averaging over $3$ independent random experiments.

\subsection{Comparison with other topologies}
We compare the consensus rate between topologies with one-peer or $\Theta(\ln(n))$ neighbors on network-size $n=300$ and $n=4900$. In the one-peer case, each topology has exactly one neighbor. For OD-EquiDyn and OU-EquiDyn, we set $M=n-1$ and $\eta=0.5$. In the  $\Theta(\ln(n))$ neighbors case, we set $M=9$ and $M=13$ for $n=300$ and $n=4900$, respectively, so that the number of neighbors is identical to the static exponential graph for a fair comparison. The results are obtained by averaging over $10$ and $3$ independent random experiments for $n=300$ and $n=4900$, respectively. 

\begin{figure}[h!]
	\centering
	\includegraphics[scale=0.4]{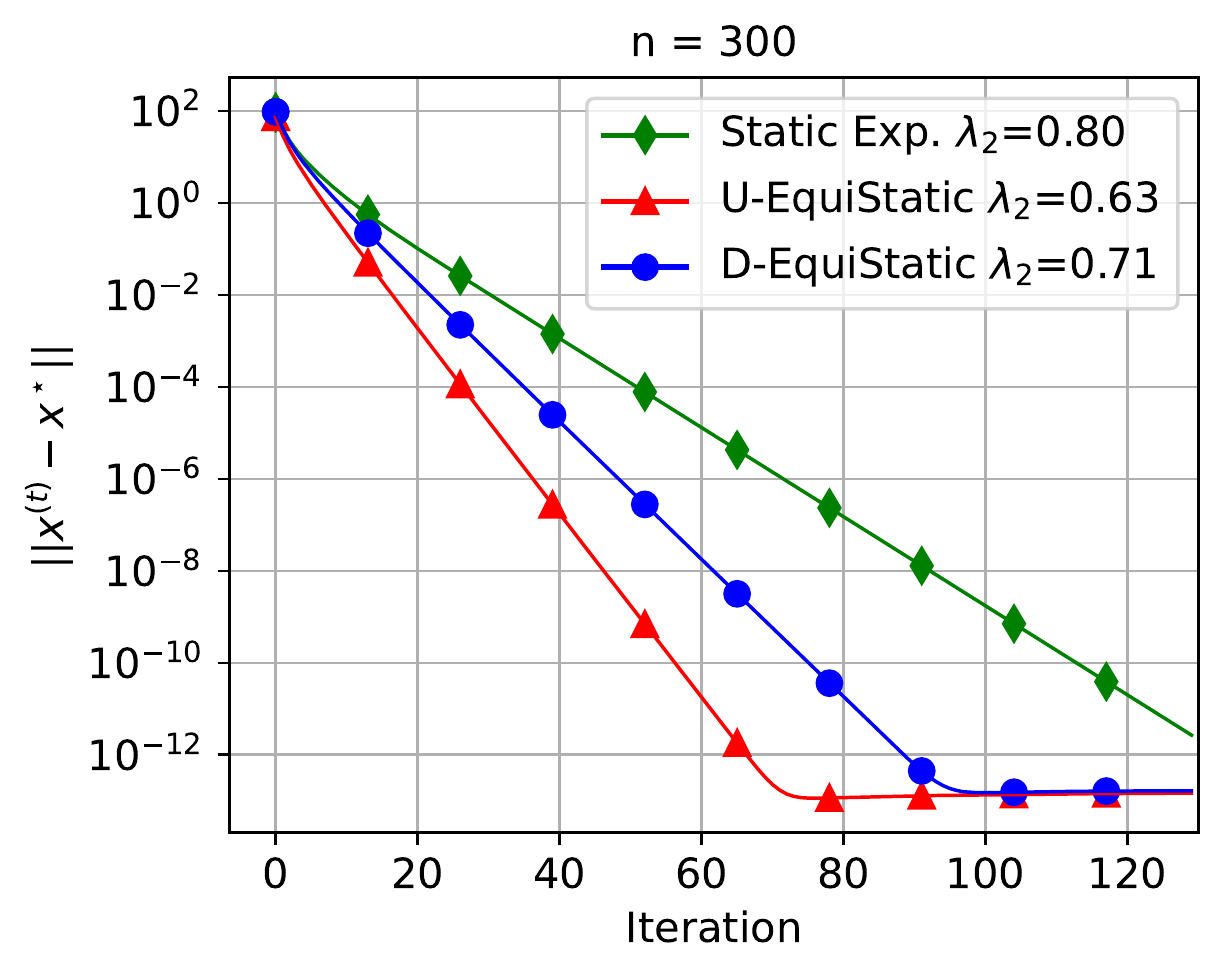}
	\includegraphics[scale=0.4]{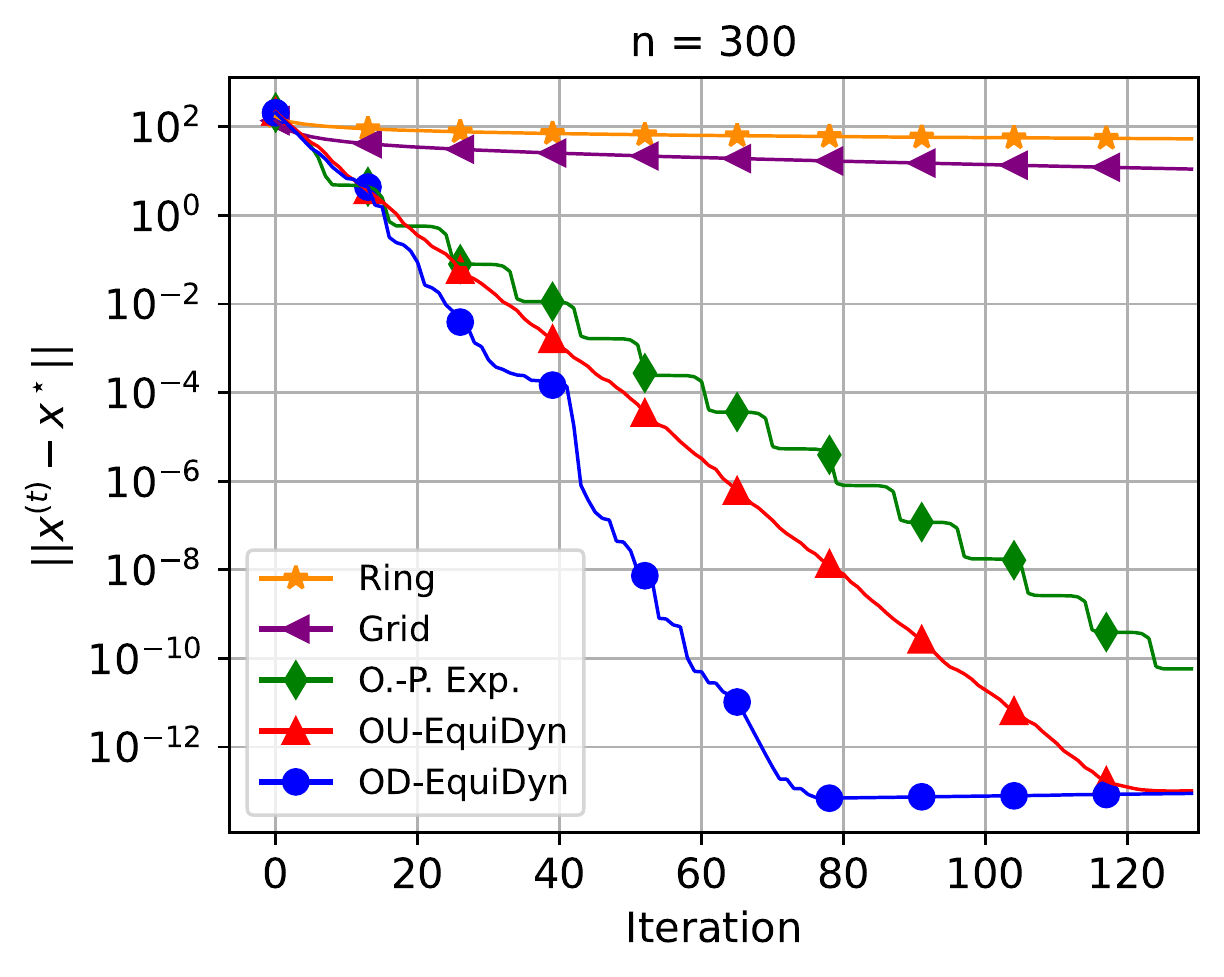}
	\includegraphics[scale=0.4]{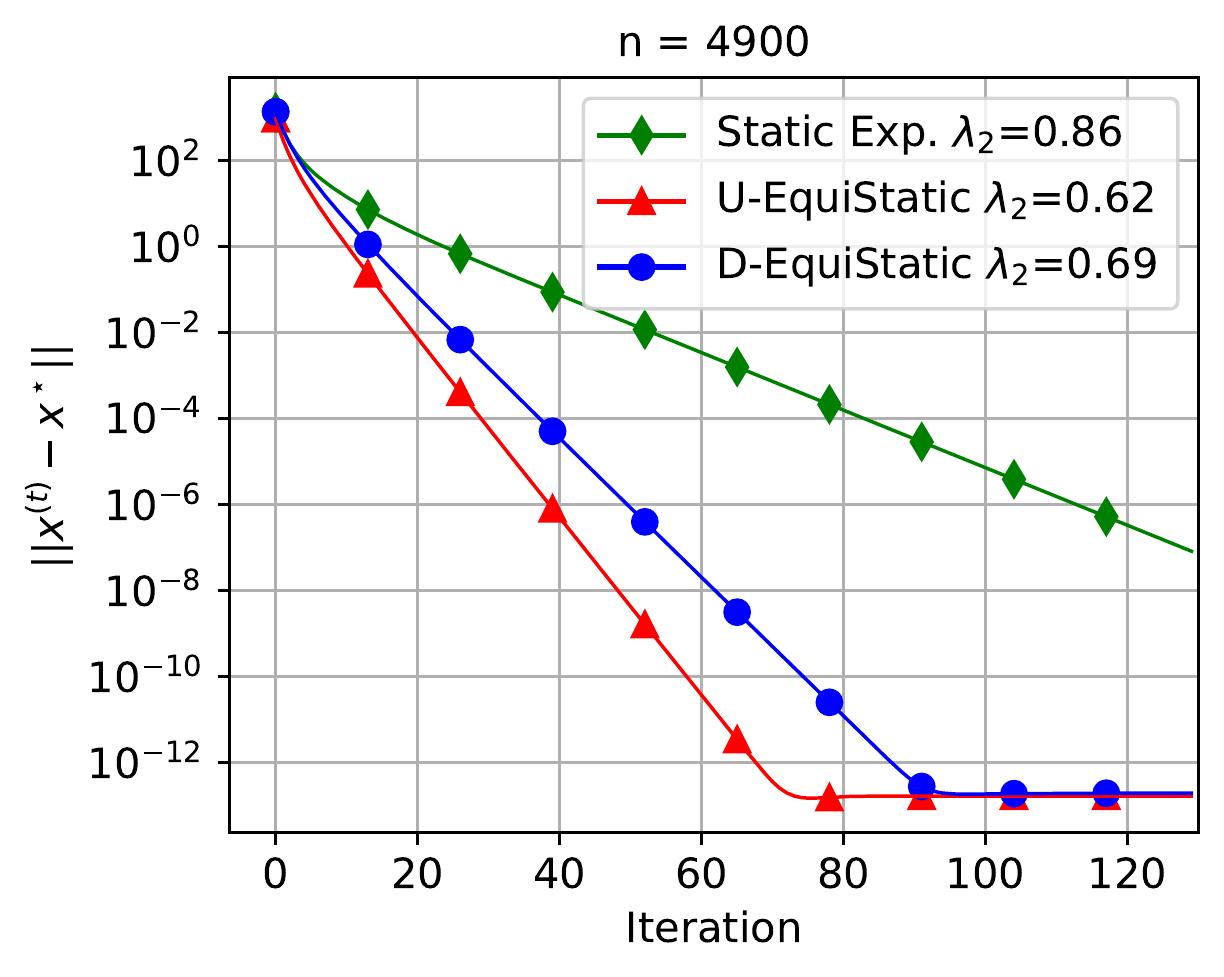}
	\includegraphics[scale=0.4]{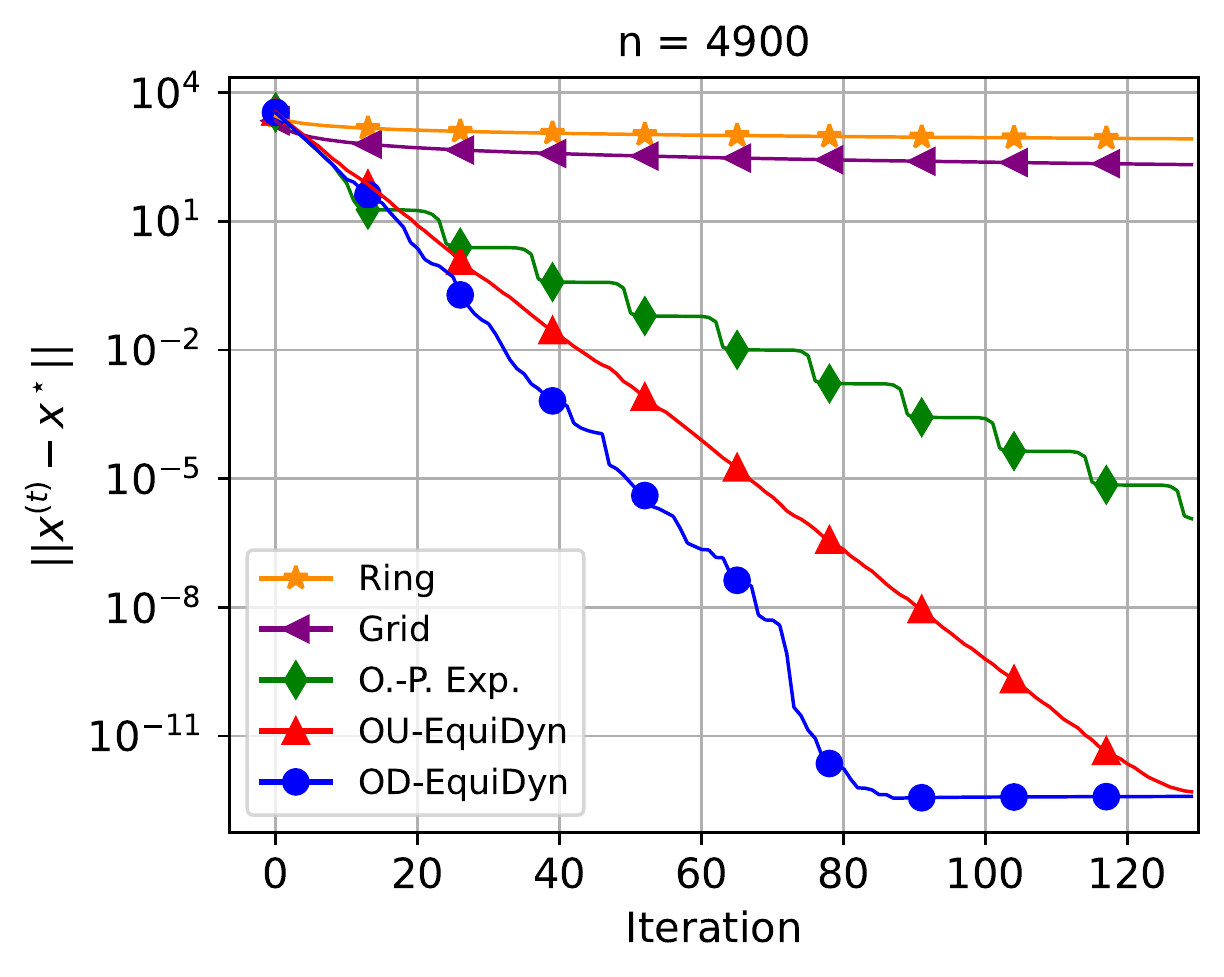}
	\caption{\small Consensus rate comparison among different network topologies in average consensus problem. Left: all graphs are with $\Theta(\ln(n))$ degree. Right: all graphs are with $\Theta(1)$ degree.}
	\label{fig:more-average-consensus}
\end{figure}

\subsection{DSGD with EquiTopo}

\textbf{Least-square}
The distributed least square problems are defined with $f_i(\x) = \|\A_i \x - \bb_i\|^2$, in which $\x\in \R^d$ and $\A_i \in \R^{K\times d}$. In the simulation, we let $d=10$ and $K=50$. At node $i$, we generate each element in $\A_i$ following standard normal distribution. Measurement $\bb_i$ is generated by $\bb_i = \A_i \x^\star + \bs_i$ with a given arbitrary $\x^\star \in \R^d$ where $\bs_i \sim \mathcal{N}(0,\sigma_s^2 \I)$ is some white noise. At each iteration $t$, each node will generate a stochastic gradient via $\widehat{\nabla f}_i(\x) = {\nabla f}_i(\x) + \bn_i$ where $\bn_i \sim \mathcal{N}(0,\sigma_n^2 \I)$ is a white gradient noise. By adjusting constant $\sigma_n$, we can control the noise variance. In this experiment, we set $\sigma_s=0.1$ and $\sigma_n=1$. 
The network size $n$ is $300$, and we set $M=9$ so that D/U-EquiStatic has the same degree as the static exponential graph. After fixing $M$, we sample the basis until the second-largest eigenvalue of the gossip matrix is small enough. The initial learning rate is $0.037$ and decays by $1.4$ every $40$ iterations. The results are obtained by averaging over $10$ independent random experiments.

\textbf{Deep learning}

\begin{figure}[h!]
	\centering
	\includegraphics[scale=0.4]{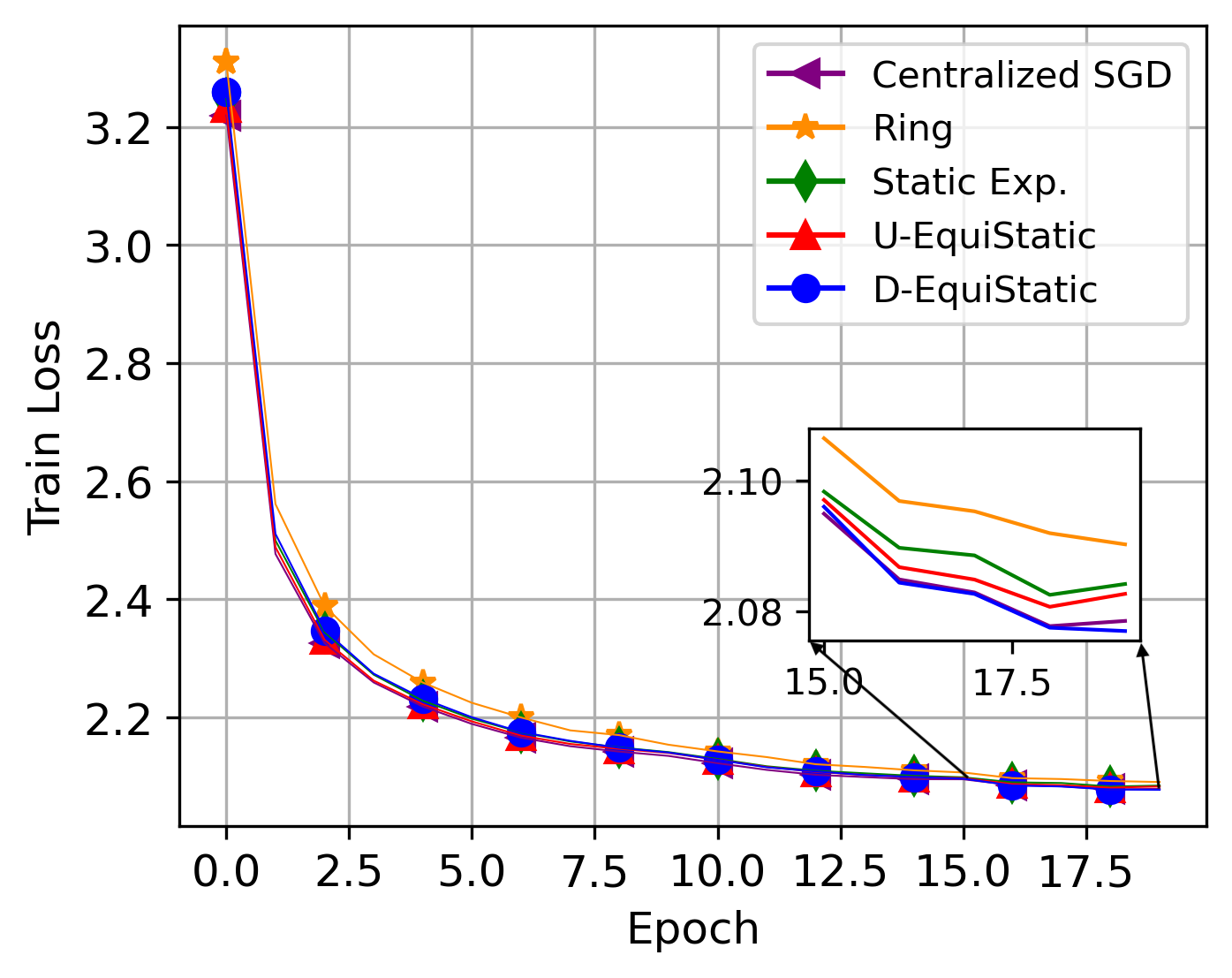}
	\includegraphics[scale=0.4]{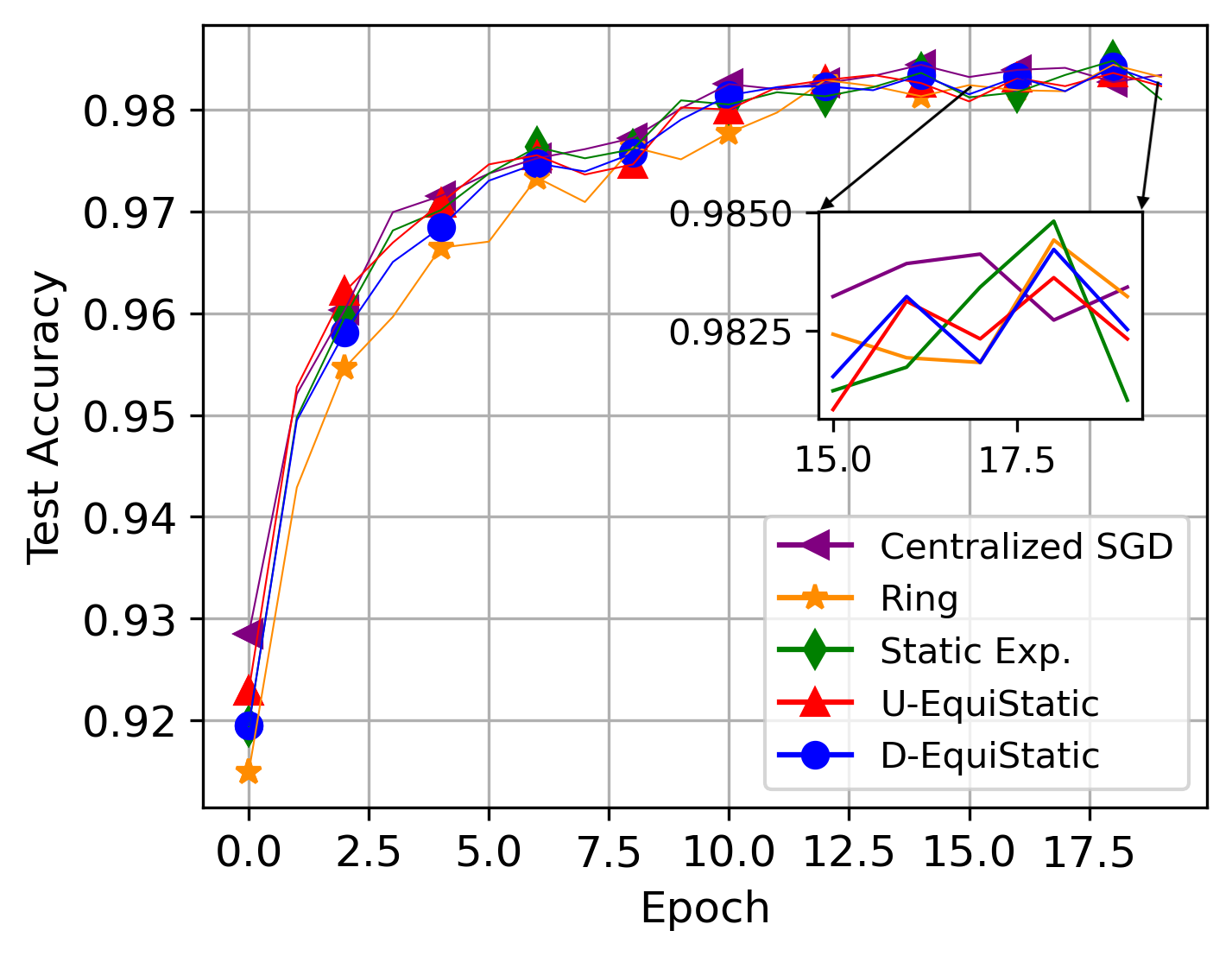}
	\includegraphics[scale=0.4]{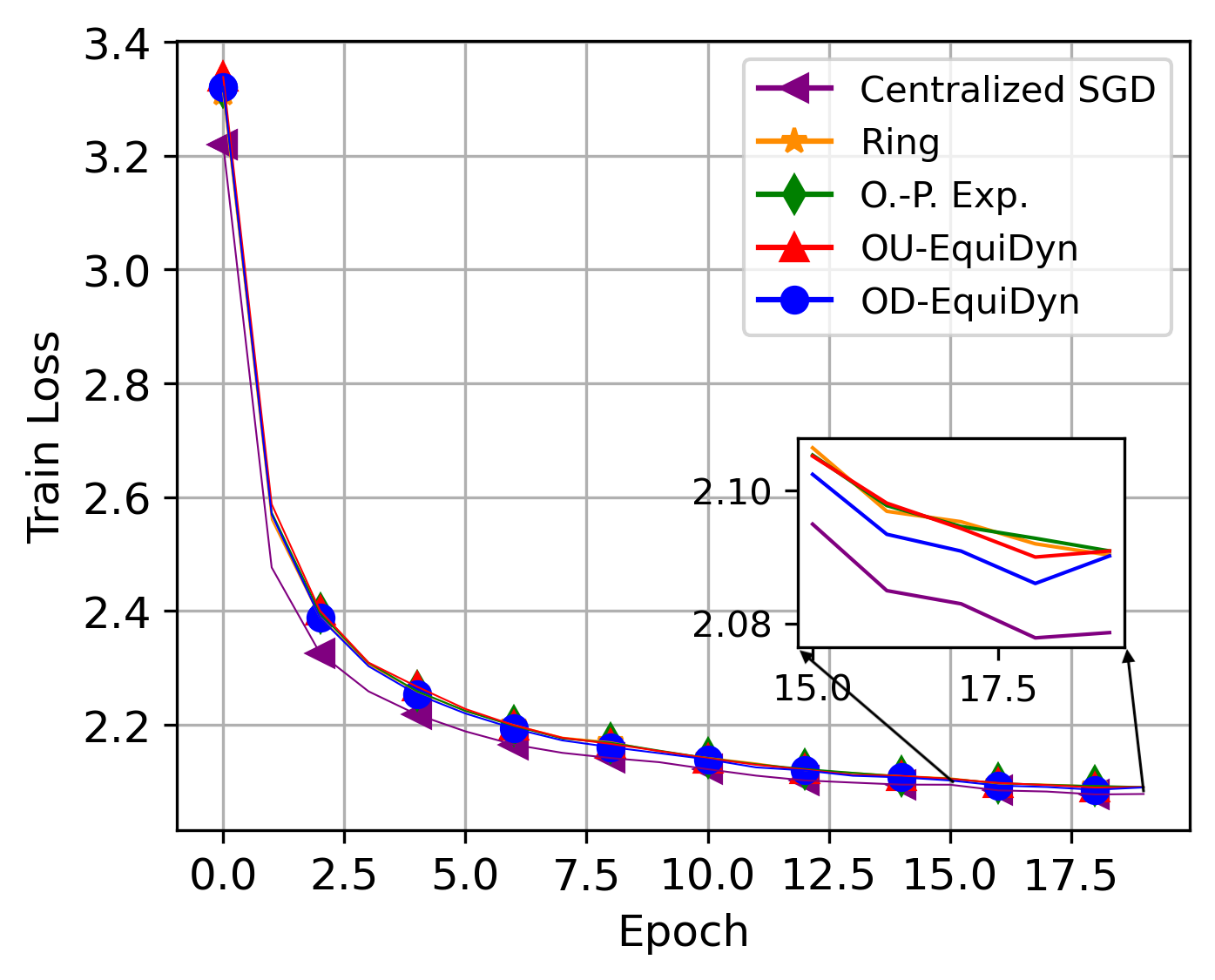}
	\includegraphics[scale=0.4]{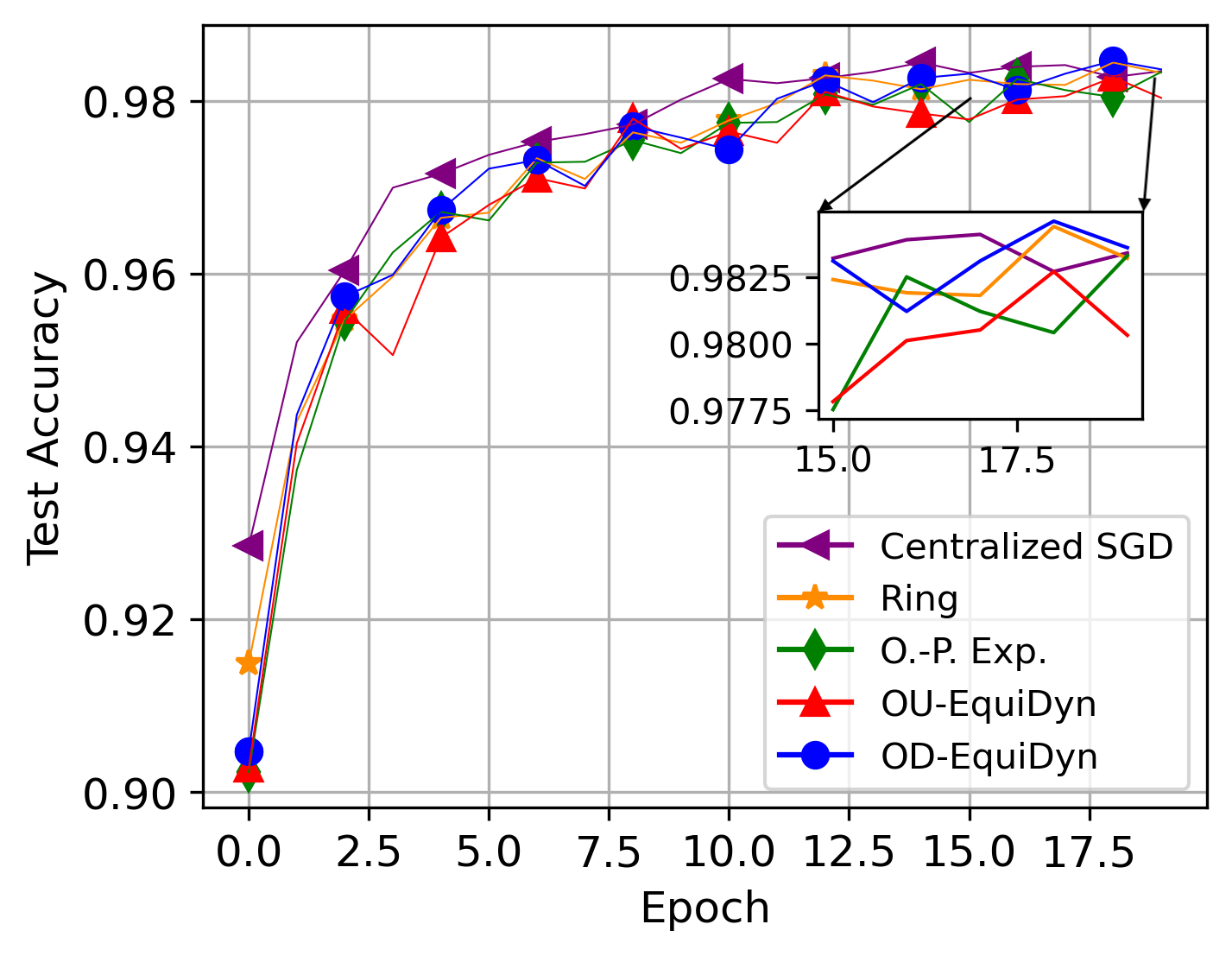}
	\caption{\small Train loss and test accuracy comparisons among different topologies for CNN on MNIST.}
	\label{fig:mnist}
\end{figure}

\textbf{MNIST.} We utilize EquiTopo graphs in DSGD to solve the image classification task with CNN over MNIST dataset \cite{lecun2010mnist}. Like the CIFAR-10 experiment, we utilize BlueFog \cite{ying2021bluefog} to support decentralized communication and topology settings in a cluster of 17 Tesla P100 GPUs. The network architecture is defined by a two-layer convolutional neural network with kernel size 5 followed by two feed-forward layers. Each convolutional layer contains a max pooling layer and a Rectified Linear Unit (ReLu). We generate D/U-EquiStaic with $M=4$ and sample OD/OU-EquiDyn with $M=16$ and $\eta=0.53$. The local batch size is $64$, momentum is $0.5$, the learning rate is $0.01$, and we train for $20$ epochs. Centralized SGD and Ring are included for comparison. See Fig.~\ref{fig:mnist} for the training loss and test accuracy of D/U-EquiStaic and OD/OU-EquiDyn graphs. See Table~\ref{Table:deep_learning_test_acc} for the test accuracy calculated by averaging over last $3$ epochs. EquiTopo graphs achieve competitive train loss and test accuracy to centralized SGD.

\begin{figure}[h!]
	\centering
	\includegraphics[scale=0.4]{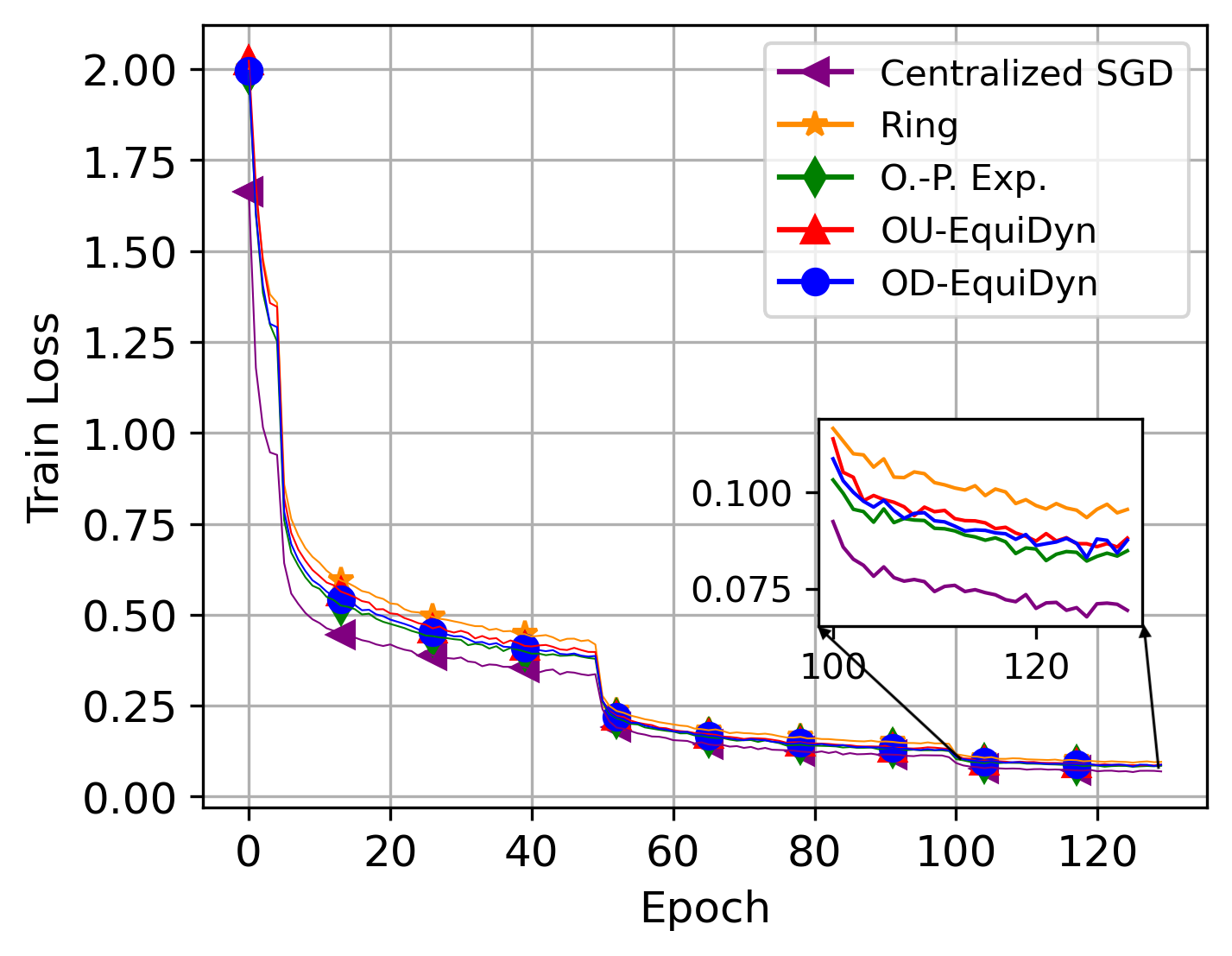}
	\includegraphics[scale=0.4]{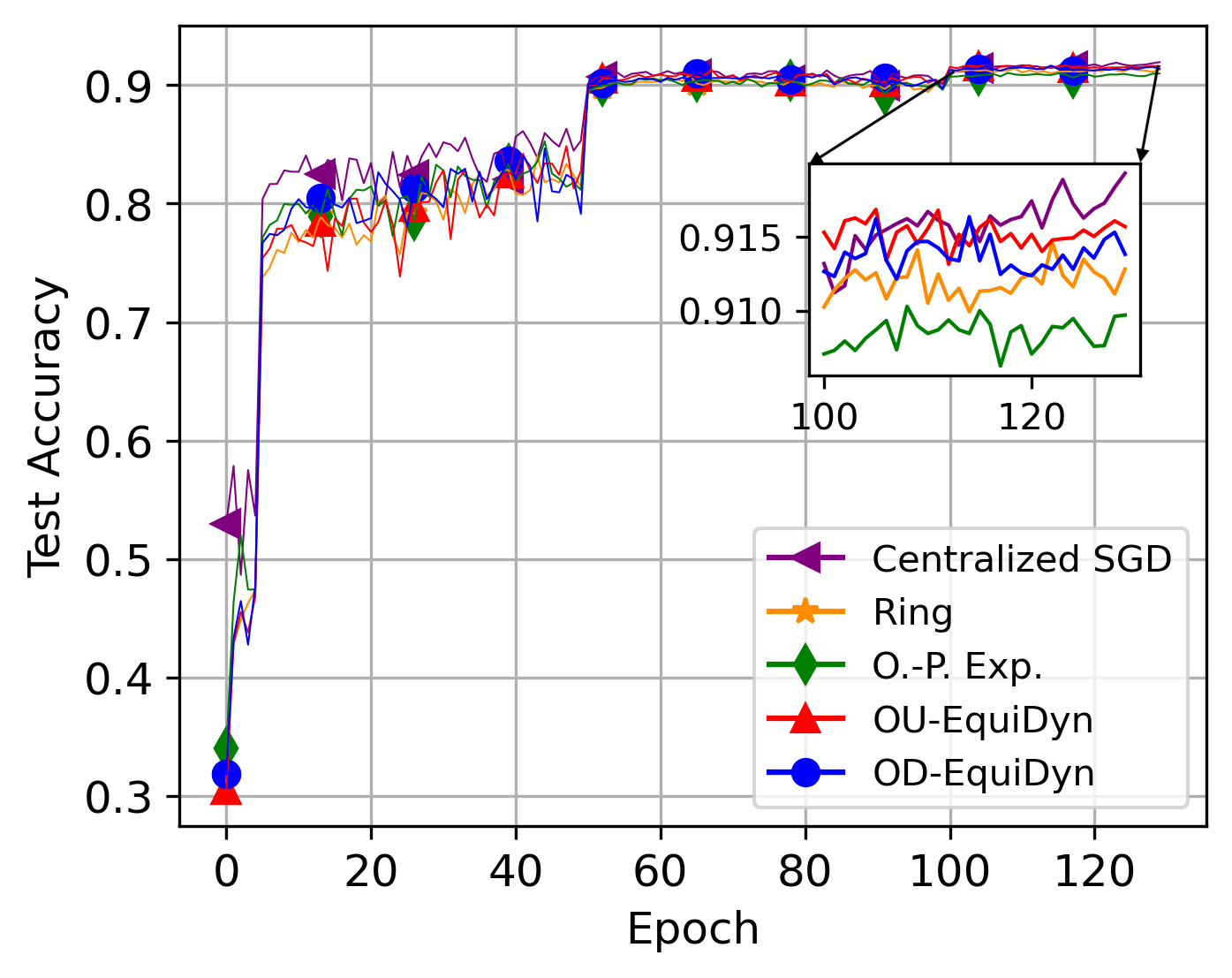}
	\caption{\small Train loss and test accuracy comparisons among different topologies for ResNet-20 on CIFAR-10.}
	\label{fig:cifar10-onepeer}
\end{figure}

\textbf{CIFAR-10.} We use the ResNet-20 model implemented by \cite{Idelbayev18a}. In this experiment, we train for 130 epochs with local batch size $8$, momentum $0.9$, weight decay $10^{-4}$, and the initial learning rate $0.01$, which is divided by $10$ at 50th, 100th, and 120th epochs. We follow the data augmentation from \cite{krizhevsky2009learning}, a $4\times4$ padding followed by a random horizontal flip and a $32\times32$ random crop. We generate D/U-EquiStaic with $M=5$
and sample OD/OU-EquiDyn with $M=16$ and $\eta=0.53$. See Fig.~\ref{fig:cifar10-onepeer} for the training loss and test accuracy of OD/OU-EquiDyn. See Table~\ref{Table:deep_learning_test_acc} for the test accuracy calculated by averaging over last $5$ epochs.

\begin{table}[h!]
    \centering 
    \caption{\small Comparison of test accuracy(\%) with different topologies over MNIST and CIFAR-10 datasets.}
    \vspace{2mm}
	\begin{tabular}{rccllc}
		\toprule
		&\textbf{Topology} &  \textbf{MNIST Acc.}  & \textbf{CIFAR-10 Acc.} &\\ \midrule
		&Centralized SGD &   98.34   & \quad\quad 91.76  &     \\
		&Ring        &    98.32      & \quad\quad 91.25  &   \\ 
		&Static Exp.  &    98.31     &\quad\quad  91.48    &      \\ 
		&O.-P. Exp.    &     98.17      & \quad\quad 90.86   &        \\ 
		&D-EquiStatic   &    98.29   &  \quad\quad\textbf{92.01}& \\ 
		&U-EquiStatic   &     98.26       &  \quad\quad91.74    &    \\ 
		&OD-EquiDyn  &     \textbf{98.39} &   \quad\quad91.44  & \\ 
		&OU-EquiDyn   &    98.12       &  \quad\quad 91.56      &    \\ 
		\bottomrule
	\end{tabular}
	
	\label{Table:deep_learning_test_acc}
\end{table}

\subsection{DSGT with EquiTopo}

In addition to the DSGD experiments, we apply the OD/OU-EquiDyn graphs to the DSGT algorithm when solving logistic regression with non-convex regularizations, i.e., $f_i(\x) = \frac{1}{L}\sum_{\ell=1}^L \ln(1 + \exp(-y_{i,\ell}\bh_{i,\ell}^T\x)) + R\sum_{j=1}^d x_{[j]}^2/(1 + x_{[j]}^2)$ where $x_{[j]}$ is the $j$-th element of $\x$, and $\{\bh_{i,\ell}, y_{i,\ell}\}_{\ell=1}^L$ is the data kept by node $i$. 
Data heterogeneity exists when local data $\xi_i$ follows different distributions $\mathcal{D}_i$ in problem \eqref{form}. To control data heterogeneity across the nodes, we first let each node $i$ be associated with a local solution $\x^\star_{i}$, and such $\x^\star_i$ is generated by $\x^\star_i = \x^\star + \bv_i$ where $\x^\star\sim \mathcal{N}(0, \I_d)$ is a randomly generated vector while $\bv_i \sim \mathcal{N}(0, \sigma^2_h \I_d)$ controls the similarity between each local solution. Generally speaking, a large $\sigma^2_h$ results in local solutions  $\{\x_i^\star\}$ that are vastly different from each other. With $\x_i^\star$ at hand, we can generate local data that follows distinct distributions. At node $i$, we generate each feature vector $\bh_{i,\ell} \sim \mathcal{N}(0, \I_d)$. To produce
the corresponding label $y_{i,\ell}$, we generate a random variable $z_{i,\ell} \sim \mathcal{U}(0,1)$. If $z_{i,\ell} \le 1 + \exp(-y_{i,\ell}\bh_{i,\ell}^T \x_i^\star)$, we set $y_{i,\ell} = 1$; otherwise $y_{i,\ell} = -1$. Clearly, solution $\x_i^\star$ controls the distribution of the labels. This way, we can easily control data heterogeneity by adjusting $\sigma^2_h$. Furthermore, to easily control the influence of gradient noise, we will achieve the stochastic gradient by imposing a Gaussian noise to the real gradient, i.e., $\widehat{\nabla f}_i(\x) = {\nabla f}_i(\x) + \bs_i$ in which $\bs_i\sim \mathcal{N}(0, \sigma^2_{n} \I_d)$. We can control the magnitude of the gradient noise by adjusting $\sigma^2_n$.

\begin{figure}[h!]
	\centering
	\includegraphics[scale=0.35]{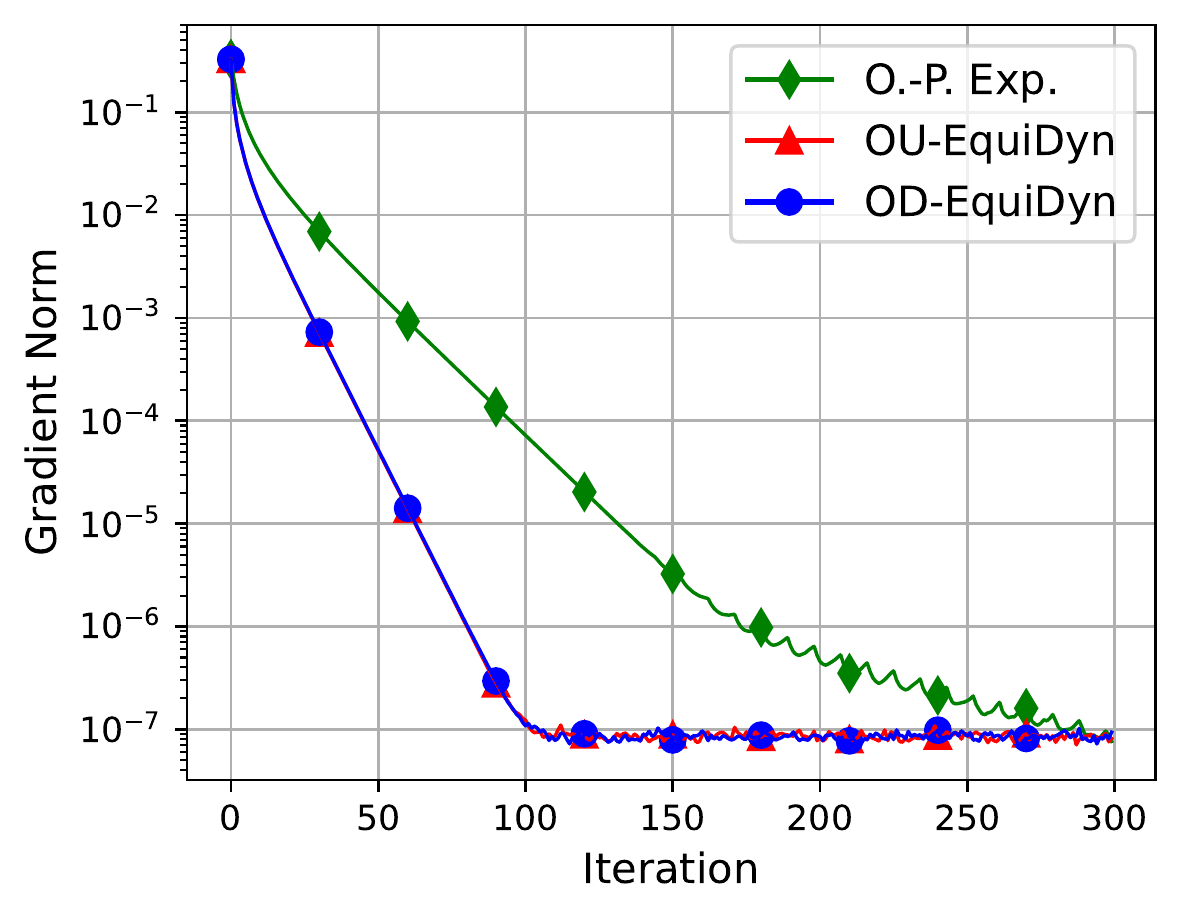}
	\includegraphics[scale=0.35]{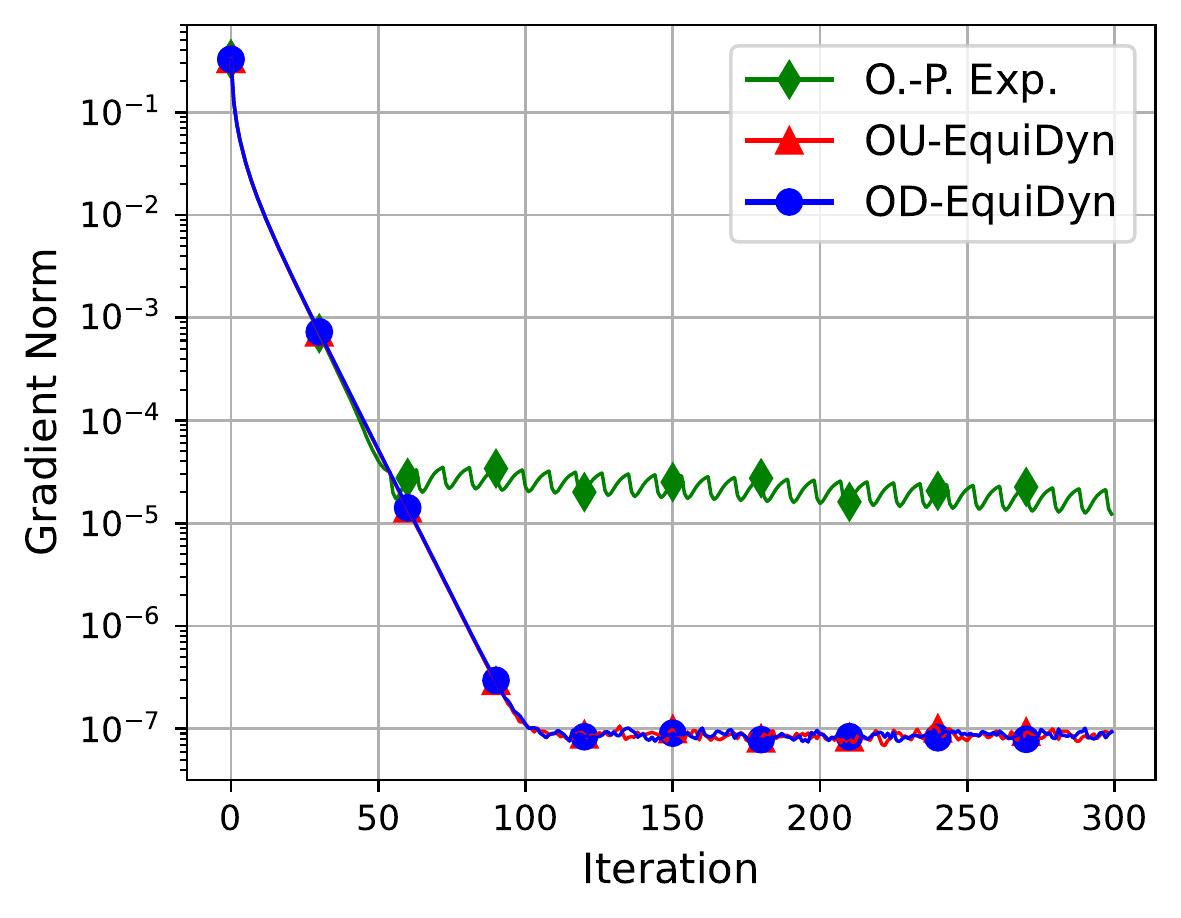}
	\caption{\small OD/OU-EquiDyn in DSGT. Left: The learning rates for O.-P. Exp. and OU/OD-EquiDyn are $1.6$ and $3$, respectively, so that all algorithms achieve the same accuracy. Right: The learning rates for all algorithms are $3$ so that they share the same convergence rate in the initial stage.}
	\label{fig:gt}
\end{figure}

We let $d=10$, $L=1000$, $n=300$, $R=0.001$, and $\sigma_h = 0.2$ in the simulation.
For OD/OU-EquiDyn, we set $M=n-1$ and $\eta=0.5$. The learning rate for OD/OU-EquiDyn and O.-P. Exp. is $3$ and $1.5$, respectively
so that all of them converge to the same level of accuracy. The left plot in Fig.~\ref{fig:gt} depicts the performance of different one-peer graphs in DSGT. The right plot in Fig.~\ref{fig:gt} illustrates how O.-P. Exp. behaves if it has the same learning rate $3$ as OU/OD-EquiDyn. The gradient norm is used as a metric to gauge the convergence performance. The results are calculated by averaging over $10$ independent random experiments. It is observed that OD/OU-EquiDyn converges faster than a one-peer exponential graph.

\end{document}